\documentclass[12pt]{amsart}
\usepackage{amsmath,amssymb,amsthm,amscd}
\usepackage{color}
\usepackage[dvips]{graphicx}
\usepackage{tabmac}

\usepackage[dvipdfmx]{hyperref}
\hypersetup{pagebackref,  
pdfpagemode=FullScreen,  
colorlinks=true}

\definecolor{gray}{gray}{.5}
\definecolor{gray2}{gray}{.8}
\newcommand{\Cell}{\color{gray2}\rule{8.5mm}{5.5mm}}
\newcommand{\cell}{\color{gray}\rule{3mm}{3mm}}
\newcommand{\smallbox}{\color{gray}\rule{2.4mm}{2.4mm}}
\newcommand{\scell}{\color{gray}\rule{1.7mm}{1.7mm}}

\newcommand{\Tableau}[2][sY]{{\text{\tableau[#1]{#2}}}}

\newcommand{\D}{\mathbb{D}}

\newcommand{\A}{\mathcal{A}}
\newcommand{\X}{\mathcal{X}}

\newcommand{\ep}{\varepsilon}

\newcommand{\la}{\lambda}

\newcommand{\R}{\mathcal{R}}
\newcommand{\C}{\mathbb{C}}

\newcommand{\Z}{\mathbb{Z}}

\newcommand{\N}{\mathbb{N}}

\newcommand{\sgn}{\mathrm{sgn}}

\renewcommand{\max}{\mathrm{max}}
\renewcommand{\min}{\mathrm{min}}
\newcommand{\T}{\mathcal{T}}
\newtheorem{thm}{Theorem}[section]
\newtheorem{lem}{Lemma}[section]

\newtheorem{prop}{Proposition}[section]
\newtheorem{cor}{Corollary}[section]
\newtheorem{Def}{Definition}[section]
\newtheorem{conj}{Conjecture}[section]

\theoremstyle{remark}
\newtheorem{example}{Example}[section]
\newtheorem{rmk}{Remark}[section]

\begin{document}

\title{$K$-theoretic analogues of factorial\\
Schur $P$- and $Q$- functions}

\begin{abstract}
We introduce two families of %are non-homogeneous symmetric
symmetric functions generalizing the % Ivanov's 
factorial
Schur $P$- and $Q$- functions due to Ivanov.
We call them $K$-theoretic analogues of factorial Schur $P$- and $Q$- functions.
We prove  
various combinatorial expressions for 
these functions, e.g. as a ratio of Pfaffians,
and a sum over excited Young diagrams.
As a geometric application, 
we show that these functions represent 
the Schubert classes 
in the $K$-theory of torus equivariant 
coherent sheaves on the maximal isotropic Grassmannians
of symplectic and orthogonal types.
This generalizes a corresponding result for the equivariant cohomology 
given by the authors.
We also discuss a remarkable property
enjoyed by these functions,
which we call the $K$-theoretic $Q$-cancellation property.
We prove that the $K$-theoretic $P$-functions 
form a (formal) basis of the ring of 
functions with the $K$-theoretic $Q$-cancellation property.
\end{abstract}

\author{Takeshi Ikeda}

\author{Hiroshi Naruse}
%\date{\today}
\date{30 April, 2013 (to appear in Adv. Math.)}


\thanks{T.I. was partially supported by Grant-in-Aid for
Scientific Research (C) 20540053}

\maketitle

\section{Introduction}
In \cite{Sc}, Schur introduced a family of symmetric polynomials, now called
Schur $Q$-functions, in order to describe the irreducible 
the characters of 
projective representations of symmetric groups.
It plays distinguished parts 
not only in representation theory but also 
in combinatorics
and geometry
(see e.g. \cite{HH}, \cite{Mac}, \cite{Pr}, \cite{St} and references therein).
In \cite{Iv}, Ivanov
introduced 
a multi-parameter deformation
of the $Q$-functions\footnote{For a case of special values of 
deformation parameter the definition is 
due to A.\:Yu.\:Okounkov. 
%See remark \ref{rem:Okounkov}
.}, which we call 
the {\it factorial} $P$- and $Q$-functions,
and proved  
various combinatorial formulas for them analogous 
to the ones for the original $Q$-functions.
In this paper, we introduce a
``$K$-theoretic'' analogue of Ivanov's functions,
and study combinatorial properties of these functions.
As an application, we show that these functions represent
the structure sheaves of the Schubert varieties  
in the torus equivariant $K$-theory of 
the maximal isotropic Grassmannians
of symplectic or orthogonal types.

Let $n$ be a positive integer. 
Let  
$\lambda=(\lambda_1,\ldots,\lambda_r)$
be a {\it strict partition\/} of length $r\leq n,$ i.e. 
a sequence of positive integers 
such that $\lambda_1>\cdots>\lambda_r.$ 
We define polynomials 
$GP_\lambda(x_1,\ldots,x_n|b)$ and $GQ_\lambda(x_1,\ldots,x_n|b)$,
which we call {\it $K$-theoretic factorial $P$- and $Q$-functions}, 
depending on a parameter $\beta$ 
and also on ``equivariant parameters'' $b=(b_1,b_2,\ldots).$
Although these polynomials depend on $\beta$, 
we do not indicate this dependence in the notation.
These are symmetric polynomials
in $x_1,\ldots,x_n$ with coefficients in
the polynomial ring $\Z[\beta,b].$
There are two important specializations.
First letting $\beta=0,$ 
our functions become
the factorial $P$- and $Q$-functions as studied in \cite{Iv}.
Second, 
let $GP_\la(x_1,\ldots,x_n)$ and 
$GQ_\la(x_1,\ldots,x_n)$ denote the functions 
obtained 
by the specialization $b_i=0$ for all $i.$
We simply call these functions 
the $K$-theoretic $P$- and $Q$-functions.
Note that even these ``non-equivariant''
versions have not appeared in the literature 
before an announcement paper \cite{INN}.
If we perform both specializations, then the functions become
the classical $P$- and $Q$-functions. 

As the name suggests,
our motivation 
of this work 
comes from geometry. 
In our previous papers \cite{Ik},\cite{EYD},\cite{rims}, 
the factorial $P$- and $Q$-functions 
of Ivanov were interpreted as 
the Schubert classes in equivariant
cohomology of the maximal isotropic Grassmannians of 
symplectic or orthogonal types. 
These spaces are generalized flag varieties $G/P$ 
of type $C,$ and $B$ or $D$ (see \S \ref{sec:geometry}).
Our original aim is to extend these results to equivariant
$K$-theory. 
Thus the primary goal of this study
is to  
prove that 
the polynomials $GP_\lambda$ and $GQ_\lambda$  
 represent the structure sheaves of Schubert varieties
in the $K$-theory of torus equivariant coherent sheaves on the maximal
isotropic Grassmannians 
of symplectic and orthogonal types (Thm. \ref{thm:main-geom}).

%Regardless of geometric applications,
%the polynomials have independent interest.
%We expect that they will play some significant role in 
%combinatorics, representation theory and 
%invariant theory. 
Another goal of this paper 
is to establish some fundamental 
 properties of the polynomials,
as Ivanov did in \cite{Iv}. 
We obtained the following different expressions for $GP_\lambda$ and $GQ_\lambda$:
\begin{itemize}
\item Hall-Littlewood-type formula (Def. \ref{def:GQGP})
\item Nimmo-type formula  
as a ratio of 
Pfaffians (Prop. \ref{thm:Nimmo})
\item Combinatorial formula using {\it shifted set-valued tableaux} (Thm. \ref{thm:SVT})
\item Combinatorial formula using {\it excited Young diagrams} (Thm. \ref{thm:EYD}).
\end{itemize}
The coincidence of these expressions is the main result about the 
combinatorics of the $K$-theoretic factorial $P$- and $Q$-functions.

\bigskip

In the course of this study,
a remarkable notion 
of a cancellation property for multi-variable polynomials arises.
%We call it
%the $K$-theoretic $Q$-cancellation property.  
Let $\beta$ be a formal parameter.

\begin{Def}\label{def:Ksuper} 
Let $G\Gamma_n$ denote the set of polynomials
$f(x_1,\ldots,x_n)$ in $\Z[\beta][x_1,\ldots,x_n],$
%is called $K$-{\it supersymmetric},
such that $f(x_1,\ldots,x_n)$ is symmetric
in variables $x_1,\ldots,x_n$, and moreover, if for all $1\leq i<j\leq n,$ the rational function  
$$f(x_1,\ldots,x_{i-1},t,x_{i+1},\ldots,x_{j-1},\ominus t,x_{j+1},\ldots,x_n)
$$
does not depend on $t$, where
$\ominus t=-t/(1+\beta t).$
\end{Def}
Obviously the set $G\Gamma_n$ 
%polynomials enjoying this property 
 forms a ring.
%which we denote by $G\Gamma_n.$
A fundamental result about this ring 
is the basis theorem, which asserts that $GP_\lambda(x_1,\ldots,x_n)$'s  
form a $\Z[\beta]$-basis of the ring $G\Gamma_n$ %of $K$-supersymmetric polynomials 
(see Thm. \ref{thm:basis} for a precise statement). 
This result generalizes work of Pragacz \cite{Pr}
in case of $\beta=0,$
where the corresponding property was called the 
$Q$-cancellation property.
So we propose to call the above property the 
{\it $K$-theoretic $Q$-cancellation property}.
This terminology might be a bit confusing in the sense that  
$GP_\lambda$'s rather than 
$GQ_\lambda$'s form a basis.
In fact, we also define a subring $G\Gamma_{n,+}$ of 
$G\Gamma_n$ (\S \ref{Gamma+}). Then $GQ_\lambda$'s form a basis
of $G\Gamma_{n,+}.$
Our excuse is that, since the term ``$Q$-cancellation property''
has been already used in the literature,
 ``$Q$'' rather than ``$P$'' will be better for reader's convenience. 
\bigskip

Let us explain key ideas of our construction in more detail.
For simplicity of notation, we restrict attention to case of type $D$ (the 
maximal isotropic Grassmannians of even orthogonal type) in 
the introduction. 
Let $\R$ be the localization of the polynomial ring $\Z[\beta][b_1,b_2,\ldots]$
by the multiplicative system formed by the products of $1+\beta b_i.$ 
The ring $\R$ is essentially the ring of Laurent polynomials of infinitely many variables.
In fact, when we specialize $\beta=-1$, $\R$ is isomorphic to  
$\Z[e^{\pm t_1},e^{\pm t_2},\ldots]$ where $b_i$ is identified with $1-e^{t_i}.$
Let $G\Gamma$ denote the ring %of $K$-supersymmetric functions 
defined to be
$\varprojlim G\Gamma_n.$ 
For each strict partition $\la$, we can 
construct a certain limit of $GP_\la(x_1,\ldots,x_n|b)$,
an element $GP_\la(x|b)$,  
in the ring $\R\otimes_{\Z[\beta]}G\Gamma.$
The functions $GP_\lambda(x|b)$ when $\la$ runs for the set 
$\mathcal{SP}$ of all strict partitions, form a `basis' of 
the ring $\R\otimes_{\Z[\beta]}G\Gamma$ over $\Z[\beta].$

Let $\Psi$ be the ring defined to be the subalgebra
of $\mathrm{Fun}(\mathcal{SP},\R)$ characterized by 
a $K$-theoretic analogue of 
the so-called Goresky-Kottwitz-MacPherson conditions (cf. \cite{GKM}, see \S \ref{ssec:GKM} for the definition).
For each $\mu\in \mathcal{SP}$, we define a sequence
$b_\mu$ of elements in $\R$ (\S \ref{ssec:vanishing}).
The localization map is defined as follows:
$$\Phi: 
\R\otimes_{\Z[\beta]}G\Gamma{\longrightarrow}
\Psi,\quad
F(x|b)\mapsto (\mathcal{SP}\ni \mu\mapsto F(b_\mu|b)\in \R)
$$
given by substitution of $x$ to $b_\mu$ 
for all $\mu\in \mathcal{SP}.$
The image $\{\Phi(GP_\la(x|b))\}_{\la\in \mathcal{SP}}$ is 
shown to be a {\it family of Schubert classes\/} (Def.\:\ref{def:Schubert}). 
It gives a distinguished `basis' of $\Psi$
as an $\R$ module (Thm.\:\ref{thm:main}). Key facts in the proof 
of the identification of `basis' is the divided difference equations (Thm.\:\ref{thm:piGQ}),
and also vanishing property of $GP_\la(x|b)$ (Prop.\:\ref{prop:vanishing}).
This technique of identification is also used in combinatorial
arguments in \S \ref{sec:comb_proof}.

The localization map can be interpreted as
the restriction map to 
torus fixed points.
Let $\mathcal{G}_n$ be the 
maximal isotropic Grassmannian (of even orthogonal type).   
An algebraic torus $T\cong (\C^{\times})^{n+1}$ acts on  $\mathcal{G}_n$ 
with finitely many fixed points naturally identified with 
the set $\mathcal{SP}(n)$ of strict partitions 
$\la$ such that $\la_1\leq n.$
The same set parametrize the Schubert varieties $\Omega_\la$ in $\mathcal{G}_n.$ 
The Grothendieck group $K_T(\mathcal{G}_n)$ of 
$T$-equivariant coherent sheaves of $\mathcal{G}_n$ 
has an $R(T)$-algebra structure, where $R(T)$ is 
the representation ring of $T$;
$R(T)$ is identified with $\Z[e^{\pm t_1},\ldots,e^{\pm t_{n+1}}].$  
As an $R(T)$-module, $K_T(\mathcal{G}_n)$ has 
a free basis formed by the classes  $[\mathcal{O}_{\Omega_\la}]_T$ of the structure sheaves
of the Schubert varieties. 
Let $i: \mathcal{G}_n^T
\hookrightarrow \mathcal{G}_n$ be the inclusion
map. Then we have a pullback 
$$
i^*: 
K_T(\mathcal{G}_n)
\longrightarrow 
K_T(\mathcal{G}_n^T)\cong 
\mathrm{Fun}(\mathcal{SP}(n),R(T)),
$$
which is known to be injective.
If we denote by $i_\mu^*$ the component of $i^*$
at the $T$-fixed point corresponding to $\mu\in\mathcal{SP}(n),$ 
Thm.\:\ref{thm:main} (combined with Thm.\:\ref{thm:main-geom})
 can be restated in the following 
geometric form.
\begin{thm}
For all $\la,\mu\in \mathcal{SP}(n)$, we have 
$$
i_\mu^*\left([\mathcal{O}_{\Omega_\la}]_T\right)=
GP_\la(b_\mu|b),
$$
where $GP_\la(b_\mu|b)$ is naturally considered to be 
an element in $R(T)$ by specialization $\beta=-1.$
\end{thm}

Therefore, via the localization map, we have constructed 
canonical surjective maps
$$\R\otimes_{\Z[\beta]}G\Gamma\longrightarrow K_T(\mathcal{G}_n)$$
for all positive integers $n$. This map sends $GP_\la(x|b)$ to 
$[\mathcal{O}_{\Omega_\la}]_T$ if $\la\in \mathcal{SP}(n)$
and to zero otherwise (see Thm.\:\ref{thm:main-geom} for a more precise statement). 
One may think of the ring 
$\R\otimes_{\Z[\beta]}G\Gamma$ 
as an equivariant $K$-ring of an infinite dimensional Grassmannian ``$\mathcal{G}_\infty$''.
This point of view will be discussed elsewhere.

\subsection{Organization}
In Section \ref{sec:GQGP}, we define
the $K$-theoretic factorial $P$- and $Q$-polynomials 
$GQ_\la(x_1,\ldots,x_n|b)$ and 
$GP_\la(x_1,\ldots,x_n|b)$
using Hall-Littlewood type summations.
These are symmetric polynomials 
with coefficients in $\Z[\beta,b].$
We show a Nimmo type formula expressing 
these functions as a ratio of Pfaffians (this result is not used 
in the rest of this paper).
This section includes a preliminary 
result on the factorial Grothendieck
polynomials (\cite{Mc}) used in the next section.
Section \ref{sec:beta-ss} establishes some basic
results for functions with $K$-theoretic $Q$-cancellation property. Within this section, we mainly consider the case
of $b=0.$ 
We prove the basis theorem using a factorization 
property of $GP_\lambda(x_1,\ldots,x_n).$
We also discuss the span of $GQ_\lambda(x_1,\ldots,x_n).$
In Section \ref{sec:Weyl}, we set up notation for Weyl groups and root systems.
In Section \ref{sec:GKM}, we introduce 
GKM ring $\Psi$ in a combinatorial
manner. We introduce 
a family of Schubert classes $\{\psi_\lambda\}_\la$
and discuss its basic properties.
These elements (if exist)  
give a natural basis for $\Psi.$
In Section \ref{sec:divdiff}, 
we prove 
the divided difference equation for 
$GQ_\lambda(x|b)$
and $GP_\lambda(x|b).$
In Section \ref{sec:localization}, we define the localization map
which relate the ring spanned by $GP_\la(x|b)$ (or $GQ_\la(x|b)$)
and the GKM ring $\Psi.$
The goal of this section is to identify $\psi_\lambda$ and
$GP_\lambda(x|b)$ (or $GQ_\lambda(x|b)$) via the localization map.
Section \ref{sec:geometry} explains geometric meaning 
of our polynomials and the ring of  
functions with $K$-theoretic $Q$-cancellation property. 
In Section \ref{sec:comb}, we state the 
results of expressing $GP_\lambda(x|b)$ and $GQ_\lambda(x|b)$
in terms of combinatorial objects, shifted set-valued
tableaux and excited Young diagrams. 
Section \ref{sec:comb_proof} is devoted to the proof
the results in Section \ref{sec:comb}.
\subsection{Related works}

For the (non-equivariant) $K$-theory of 
odd maximal orthogonal Grassmannians,
Clifford, Thomas, and Yong \cite{CTY} proved 
a Littlewood-Richardson rule which 
gives an explicit combinatorial description 
for the structure constants for the Schubert basis.
In view of Cor. \ref{thm:non-eq}, these 
constants are equal to the structure constant for
$GP_\la(x_1,\ldots,x_n)$'s.
Note that \cite{CTY} is based on a work of Buch and Ravikumar \cite{BR}. 
Our result in this paper is independent from their result.

For the Lagrangian Grassmannians, 
Ghorpade and Raghavan \cite{GR} 
gave a detailed description of Gr\"obner degeneration of an open 
affine piece of 
%of the coordinate ring of 
a Schubert
variety around a torus fixed point.
Using this result, Kreiman \cite{Kr}, 
derived 
a combinatorial expression for the restriction of Schubert 
to fixed point set of the torus action
class of the equivariant $K$-theory.
Note also that Raghavan and Upadhyay extend \cite{GR}
to orthogonal types.

It will be worth noticing that 
in \cite{DRV}, the center of the
affine BMW-algebra is 
identified with a subalgebra of the ring of symmetric 
Laurent polynomials 
having a cancellation property quite similar to the $K$-theoretic
$Q$-cancellation property. 

It is a natural task to 
extend the $K$-theoretic factorial $P$- and $Q$-functions to 
family of functions relevant for the equivariant $K$-theory of the 
full flag varieties of 
types $B,C,D.$
Kirillov and the second author \cite{KN} introduced such family of functions, the double Grothendieck polynomials
of type $B,C$ and $D$, by
using Id-Coxeter algebra, a $K$-theoretic analogue of the nil-Coxeter algebra.

An outstanding open question is the determination of
the structure constant for $GP_\la(x|b)$ and $GQ_\la(x|b).$
As for the ``equivariant case'', i.e. for the case of arbitrary $b,$
we do not even have a conjecture, except for 
very restricted cases.

In \cite{INN}, we introduced a bumping procedure for the set-valued shifted tableaux.
As an application, we proved a Pieri type rule for $GQ_\lambda(x_1,\ldots,x_n)$,
which is combinatorially equivalent to a problem of 
counting the $KLG$-tableau in the language of \cite{BR},
where they derive the description 
utilizing a more geometric arguments. 
Bumping procedure for 
`orthogonal' type, which is relevant for $GP_\la(x_1,\ldots,x_n)$, is not known now.

\bigskip

After the first version of this paper was completed,
a closely related work by Graham and Kreiman \cite{GrKr}
has appeared. 
They give combinatorial descriptions (using excited Young diagrams and set-valued tableaux)
of the restrictions to $T$-fixed points of the classes of structure sheaves of Schubert 
varieties in the $T$-equivariant $K$-theory of 
the Lagrangian Grassmannians 
and the maximal isotropic  Grassmannians of 
orthogonal types. 
%$\mathcal{G}_n$
%(types $B,C$ and $D$). 
%The method 
%of derivation in \cite{GrKr} is a generalization of our previous result \cite{EYD}
%for $T$-equivariant cohomology. 
Their expressions for the $K$-theoretic restrictions 
differ from the ones obtained from the results in the present paper. 
%In particular, the formula in \cite{GrKr} is {\it positive\/} in
%the sense of \cite{GrKu}.   
See \cite{GrKr}, \S 5.4 for precise comparison.  

\bigskip

{\bf Acknowledgments.} We thank Mark Shimozono for discussions
and his continuous interest in the work.
We also like to thank 
Leonardo Mihalcea and Alexander Yong 
and Hugh Thomas for helpful discussions and 
useful comments
on the draft. 
We are indebted to the anonymous referees for improvements
to the manuscript
(especially Prop. \ref{prop:s_ila} and Thm. \ref{thm:piGQ}).

\bigskip

Note added: In the 5-th MSJ, Seasonal Institute, Schubert calculus 
(Osaka, 2012 July),  
we learned from Alexander Yong and Hugh Thomas
that the result in \cite{CTY}
implies that Conjeture \ref{conj:GP} is true.
%implies Conjecture \ref{conj:GP} 
The details will be presented elsewhere. 
%We thank Alexander Yong and Hugh Thomas
%for the explanation. 

\setcounter{equation}{0}
\section{$K$-theoretic factorial Schur $Q$- and $P$-functions}\label{sec:GQGP}
We define $K$-theoretic analogue of 
factorial $P$- and $Q$-functions.
This generalize a form for $P$- and $Q$-functions
as a particular case of the 
Hall-Littlewood functions with the parameter $t$ equal to $-1$.
These are symmetric polynomials
in $x_1,\ldots,x_n$ with coefficients in $\Z[\beta,b].$
We give another expression for these polynomials 
as ratios of Pfaffians.  
We also give some supplementary discussion on type $A$ case, i.e.
on the factorial Grothendieck polynomials   
defined and studied by McNamara.  
\subsection{The polynomials $GP_\lambda(x_1,\ldots,x_n|b)$ and $GP_\lambda(x_1,\ldots,x_n|b)$}
Let $\beta$ be a parameter. 
Define binary operators $\oplus$ and $\ominus$ by
$$
x\oplus y:=x+y+\beta x y,\quad x\ominus y:=\dfrac{x-y}{1+\beta y}.$$
We also define a deformation of $k$-th powers of $x$ with parameters
$b=(b_1,b_2,...)$ by 
$$
{
[x|b]^k:=(x\oplus b_1)(x\oplus b_2)\cdots (x\oplus b_k)},$$
and its variant by 
$
{[[x|b]]^k:=(x\oplus x)(x\oplus b_1)(x\oplus b_2)\cdots (x\oplus b_{k-1})}.$

Let $\lambda=(\lambda_1,\ldots,\lambda_r)$ be a strictly decreasing 
sequence of positive integers. We call $\lambda$ a {\it strict partition},
and $r$ the {\it length} of $\lambda.$
Let $\mathcal{SP}_n$ denote the set of all strict partitions 
$\lambda$ such that the length $r\leq n.$
We set $[x|b]^{\lambda}=\prod_{j=1}^r[x_j|b]^{\lambda_j}$ and 
$[[x|b]]^{\lambda}=\prod_{j=1}^r[[x_j|b]]^{\lambda_j}.$

\begin{Def}\label{def:GQGP} Let $\lambda=(\lambda_1,\cdots,\lambda_r)$ be a strict partition in $\mathcal{SP}_n.$
We define functions on $n$ variables $x_1,\ldots,x_n$ with parameters
 $\beta$ and $b_1,b_2,\ldots$ as follows:
\begin{eqnarray}
GP_\lambda(x_1,\ldots,x_n|b)&:=&
\displaystyle\frac{1}{(n-r)!}
\sum_{w\in S_n}w\left[
[x|b]^{\lambda}
\prod_{i=1}^r
\prod_{j=i+1}^n \frac{x_{i}\oplus x_{j}}{x_{i}\ominus x_{j}}\right],\label{def:GP}
\\
%%%%%%%%%%%%%%%%%%%%%%%%%%%%%%%%%%%%%%%%
GQ_\lambda(x_1,\ldots,x_n|b)&:=&
\displaystyle\frac{1}{(n-r)!}
\sum_{w\in S_n}w\left[
[[x|b]]^{\lambda}
\prod_{i=1}^r
\prod_{j=i+1}^n \frac{x_{i}\oplus x_{j}}{x_{i}\ominus x_{j}}
\right].\label{def:GQ}
\end{eqnarray}
where $w\in S_n$ acts as a permutation on the variables $x_1,\ldots,x_n$.
\end{Def}

These functions are obviously 
symmetric in the variables $x_1,\ldots,x_n,$ 
and are polynomials 
with coefficients in $\Z[\beta][b_1,b_2,\ldots,b_{m}]$ with $m=\lambda_1$ for
$GP_\lambda$ and $m=\lambda_1-1$ for $GQ_\lambda.$
The fact that these are polynomials follows from 
a standard argument. See \cite{Mac} Chap. III, 1, or \cite{Iv2}, Prop. 1.1 (a) 
for a proof.
We call these polynomials the {\it $K$-theoretic factorial $P$- and $Q$-functions}.

 \begin{rmk}[Comments on terminology]\label{rem:Okounkov}
If we set $\beta=0$ and replace $[x|b]^k$ with the $k$th {\it falling factorial\/}
 $\prod_{i=1}^k(x-i+1)$
in the expression (\ref{def:GP}), the polynomial coincides with
the {\it factorial} $P$-function, whose definition is due to A.\:Yu.\:Okounkov (see \cite{Iv2}, \cite{Iv3}).
 Ivanov \cite{Iv} developed the theory by using 
 $(x|a)^k:=\prod_{i=1}^k(x-a_i)$ with arbitrary sequence
 of parameters
 $a=(a_i)_{i\geq 1}.$
The function 
was first called the {\it generalized  $Q$-($P$-) functions} in \cite{Iv3},
and then the {\it multi-parameter Schur $P$-function}
in \cite{Iv}.
In \cite{Ik}, we simply called Ivanov's function the
{\it factorial Schur $Q$-function}. 
This terminology is 
consistent with a convention in 
type $A$ case (e.g. \cite{MS}).
 \end{rmk}

Both the functions $GP_\lambda(x_1,\ldots,x_n|b)$ and $GQ_\lambda(x_1,\ldots,x_n|b)$ 
can be expanded in the following forms:
$$
GP_\lambda(x_1,\ldots,x_n|b)=
\sum_{k\geq |\la|}
\beta^{k-|\la|}
GP_{\lambda}(x_1,\ldots,x_n|b)_k,$$
$$
GQ_\lambda(x_1,\ldots,x_n|b)
=\sum_{k\geq |\la|}
\beta^{k-|\la|}
GQ_{\lambda}(x_1,\ldots,x_n|b)_k,
$$
where 
 $GP_{\lambda}(x_1,\ldots,x_n|b)_k$ and $GQ_{\lambda}(x_1,\ldots,x_n|b)_k$ are 
symmetric polynomials in $x_1,\ldots,x_n$ with coefficients in
 $\Z[b]=\Z[b_1,b_2,\ldots]$, which are homogeneous of 
degree $k$ when we set $\deg(x_i)=\deg(b_i)=1.$
The lowest homogeneous parts are the factorial 
Schur functions, $P_\lambda(x_1,\ldots,x_n|b)$
and $Q_\lambda(x_1,\ldots,x_n|b)$ respectively. 
It is known that $P_\lambda(x_1,\ldots,x_n|b)$'s (resp.  $Q_\lambda(x_1,\ldots,x_n|b)$'s) 
$,\la \in \mathcal{SP}_n,$ are linearly independent over $\Z[b]$ (cf. \cite{Iv}).
From this fact, it turns out that  $GP_\lambda(x_1,\ldots,x_n|b)$'s (resp. $GP_\lambda(x_1,\ldots,x_n|b)$'s), $\la\in \mathcal{SP}_n,$
are linearly independent over the ring $\Z[\beta,b].$

\begin{rmk}
Here we adopt 
the notation in \cite{DSP} for the factorial $P$- and $Q$-functions,
with the identification $b_i=-t_i.$
This is slightly different from the original 
one in \cite{Iv}. See \cite{DSP}, \S 4.2 for the precise 
correspondence.
\end{rmk}

The specialization obtained by setting $\beta=-1$ is called 
the {\it $K$-theoretic specialization}, because 
the case 
is relevant when we apply 
these functions to $K$-theory.
If we denote by 
$F_\lambda(x_1,\ldots,x_n|b)$ the $K$-theoretic specialization 
of $GP_\lambda(x_1,\ldots,x_n|b)$, then we have  
$$(-\beta)^{-|\lambda|}
F_\lambda(-\beta x_1,\ldots,-\beta x_n|-\beta b_1,-\beta b_2,\ldots)
=GP_\lambda(x_1,\ldots,x_n|b)$$
and similarly for $GQ_\lambda.$ 
Thus 
the case of arbitrary $\beta$ is recovered from 
the $K$-theoretic specialization.
In an application to torus equivariant $K$-theory in \S \ref{sec:geometry}
the parameters $b_i$ are identified with
$1-e^{t_i}$ where $e^{t_i}$ is a standard character of 
the torus.

If all the deformation parameters $b_i\;(i\geq 1)$ are specialized to 
zero, 
then  
the $K$-theoretic factorial $P$- and $Q$- functions 
are simply denoted by $GP_\lambda(x_1,\ldots,x_n)$
and $GQ_\lambda(x_1,\ldots,x_n).$
These are relevant to the non-equivariant 
$K$-theory of the isotropic Grassmannians.

\subsection{Some useful identities.}

The next result will be used in the proofs of Thm.\:\ref{thm:piGQ}.

\begin{lem}\label{lem:Id_C} We have
\begin{equation}
\sum_{i=1}^m\frac{x_i\oplus x_i}{x_i\ominus t}
\prod_{j\neq i}\frac{x_i\oplus x_j}{x_i\ominus x_j}
+\prod_{i=1}^m
\frac{t\oplus x_i}{t\ominus x_i}=1.\label{eq:C}
\end{equation}
\end{lem}
{\it Proof.} 
We prove this by 
induction on $m.$
We can write the left-hand side as 
$$
\frac{F(x_m)}{(x_m-t)\prod_{i=1}^{m-1}(x_m-x_i)}
$$
where $F(x_m)$ is a polynomial in $x_m.$
We claim that $\deg F(x_m)\leq m.$ 
It is easy to see that 
$\deg F(x_m)\leq m+1.$ 
By letting $x_i$ to $\ominus x_i$ 
and $t$ to $\ominus t$ in (\ref{eq:C}) for $m-1$, we have
$$
\sum_{i=1}^{m-1}\frac{x_i\oplus x_i}{x_i-t}
\prod_{j\neq i}\frac{x_i\oplus x_j}{x_i-x_j}
+\prod_{i=1}^{m-1}
\frac{t\oplus x_i}{t-x_i}=
\prod_{i=1}^{m-1}(1+\beta x_i).
$$
This equation implies that the coefficient of degree $m+1$ of $F(x_m)$ vanishes.
Therefore 
it suffices to show the equation 
$F(x_m)=(x_m-t)\prod_{i=1}^{m-1}(x_m-x_i)$
for $m+1$ different values of $x_m$.
In fact, by letting $x_m=\ominus x_i\;(1\leq i\leq m-1)$ and $x_m=0,\ominus t$, 
we can check the equality by using inductive hypothesis. $\qed$

\bigskip

The next result will be used in the proofs of Thm.\:\ref{thm:piGQ} and also 
Prop.\:\ref{prop:A10}.
Analogous equation for the case of $\beta=0$ 
can be found in \cite{Iv3}, Prop.\:2.4.

\begin{lem}\label{var0-Iv}
For $k=0,1,2$, we have
\begin{equation}
\sum_{i=1}^{n}
(1+\beta x_i)^k
\prod_{j\neq i}\frac{x_i\oplus x_j}{x_i\ominus x_j}
=\begin{cases}
\prod_{i=1}^n(1+\beta x_i)^k-\prod_{i=1}^n(1+\beta x_i)
&\;\mbox{if}\; n\;\mbox{is even}\\
\prod_{i=1}^n(1+\beta x_i)^k&\;\mbox{if}\; n\;\mbox{is odd}.
\end{cases}
\end{equation}
\end{lem}

{\it Proof.} 
Let $E_{n,k}$ denote the equation.
We prove $E_{n,k}\;(k=0,1,2)$ by induction on $n.$
For $n=1,2$, the verification is straightforward.
Suppose $n$ is odd and $n\geq 3.$ 
Assume $E_{m,k}$ hold for $m<n,k=0,1,2.$ 
We write the left-hand side of $E_{n,0}$ as 
\begin{equation}
\frac{F(x_n)}{\prod_{i=1}^{n-1}(x_n-x_i)}\label{eq:frac}
\end{equation}
where $F(x_n)$ is a polynomial in $x_n.$
It is easy to see that $\deg F(x_n)\leq n.$
Using $E_{n-1,1}$, we see that the coefficient of $x_n^{n}$ in 
$F(x_n)$ is zero.
Hence we have $\deg F(x_n)\leq n-1.$
Now in order to prove $E_{n,0}$, 
we only have to show this for $n$ distinct values of $x_n.$
Let us put $x_n=\ominus x_i\;(1\leq i\leq n-1).$
Then $E_{n,0}$ is reduced to $E_{n-2,0}.$
If we put $x_n=0$ then 
$E_{n,0}$ is reduced to $E_{n-1,0}.$
Thus we have $E_{n,0}.$ Note that 
$E_{n,2}$ is obtained from $E_{n,0}$ by $x_i\mapsto \ominus x_i
\;(1\leq i\leq n).$

Next we prove 
\begin{equation}
\sum_{i=1}^{n}
(1+\beta x_i)
\prod_{j\neq i}\frac{x_i\oplus x_j}{x_i\ominus x_j}
-\prod_{i=1}^n(1+\beta x_i)=0
\end{equation}
which is equivalent to $E_{n,1}.$ 
Write the left-hand side as in (\ref{eq:frac}).
Then by using $E_{n-1,2}$, we see that $\deg F(x_n)\leq n-1.$
Now evaluation $x_n=\ominus x_i\;(1\leq i\leq n-1)$ is 
reduced to $E_{n-2,1},$ whereas the specialization $x_n=0$ 
is reduced to $E_{n-1,1}.$
Thus we obtained $E_{n,1}.$

In a similar way, $E_{n+1,0}$ follows from $E_{n,1},\;E_{n,0},$ and $E_{n-1,0}.$
Then we have $E_{n+1,2}$ by switching $x_i\mapsto \ominus x_i\;(1\leq i\leq n).$ 
Finally $E_{n+1,1}$ follows from $E_{n,2},\;E_{n,1},$ and $E_{n-1,1}.$
$\qed$

\bigskip

For positive integers $n,r$ such that $r\leq n$,
we define
$$
\Phi_{n,r}(x)=\prod_{i=1}^r\prod_{j=i+1}^n\frac{x_i\oplus x_j}{x_i\ominus x_j}.
$$
We need the next lemma to prove Prop. \ref{prop:GPexpand}, Thm. \ref{thm:piGQ}.

\begin{lem}\label{lem:big-r} Let $n,r$ be such that $n>r$ and 
$n-r$ is odd.
If the length of $\la$ is less than or equal to $r,$
then we have $$
GP_\la(x_1,\ldots,x_n)
=\frac{1}{(n-r)!}\sum_{w\in S_{n}}w\left(x^\la\Phi_{n,r}(x)\right).
%=\frac{1}{(n-r+1)!}\sum_{w\in S_{n}}w\left(x^\mu\Phi_{n,r+1}(x)\right)
$$
\end{lem}
{\it Proof.} Apply Lemma \ref{var0-Iv}, when the number of variables is odd
and $k=0.$ $\qed$

\subsection{Nimmo type formula for $GQ_\lambda,GP_\lambda$}
Let $A$ be a skew-symmetric matrix of even size. 
We denote by $\mathrm{Pf}(A)$ 
the Pfaffian of $A.$
There is a formula which is an analogue of Schur's 
evaluation of Pfaffian \cite{Sc}.
\begin{lem}\label{Pf-Schur}
We have
\begin{equation}
\underset{1\leq i<j\leq 2m}
{\mathrm{Pf}}\left(\dfrac{x_i-x_j}{x_i\oplus x_j}\right)
=\prod_{1\leq i<j\leq 2m}\dfrac{x_i-x_j}{x_i\oplus x_j}.\label{eq:Pf-Schur}
\end{equation}
\end{lem}
{\it Proof.} 
By Knuth's theorem \cite{Knu}, we only have to show this for 
$m=2$, which is straightforward.  
Indeed, (\ref{eq:Pf-Schur}) is a special case of (4.4) in \cite{Knu}.
$\qed$

\bigskip

 Nimmo \cite{Ni} derived an expression for 
 $Q_\lambda(x_1,\ldots,x_n)$ 
 as a ratio of Pfaffians.
The next result is a generalization of Nimmo's formula 
for $GQ_\lambda(x_1,\ldots,x_n|b)$ and $GP_\lambda(x_1,\ldots,x_n|b).$
\begin{thm}\label{thm:Nimmo} 
Let $\lambda$ be a strict partition of length $r\leq n.$ 
Let $m=r$ if $n-r$ is even, and 
$m=r+1$ if $n-r$ is odd, and set $\la_{r+1}=0.$
Let $A_0=A_0(x_1,\ldots,x_n)$ be a skew symmetric $n\times n$ matrix
with $(i,j)$ entry $(x_i-x_j)/(x_i\oplus x_j),$ and let  
$B_\lambda$ be an $n\times m$ matrix with $(i,j)$ entry
$[[x_i|b]]^{\lambda_{m-j+1}}(1+\beta x_i)^{-j}.$
Let  
\begin{equation}
A_\la(x_1,\ldots,x_n|b)=
\left(\begin{array}{cc}
A_0 & B_\lambda\\
-{}^t\! B_\lambda & 0 
\end{array}\right),
\end{equation}
which is a skew symmetric matrix of $(n+m)\times (n+m).$
Put
$$\mathrm{Pf}_0(x_1,\ldots,x_n)=
\begin{cases}\mathrm{Pf}\,A_0(x_1,\ldots,x_n) &\mbox{if}\; n\;\mbox{is even}\\
\mathrm{Pf}\,A_0(x_1,\ldots,x_{n},0)&\mbox{if}\; n\;\mbox{is odd}\end{cases}.
$$
Then 
\begin{eqnarray*}
GQ_\lambda(x_1,\ldots,x_n|b)
=
\prod_{i=1}^n(1+\beta x_i)^m
\frac{\mathrm{Pf}\,A_\la(x_1,\ldots,x_n|b)}{
\mathrm{Pf}_0(x_1,\ldots,x_n)}.
\label{eq:Nimmo}
\end{eqnarray*}
The similar formula holds for $GP_\lambda$ 
when we use 
$[x|b]^k$ instead of $[[x|b]]^k.$
\end{thm}

Our proof of this theorem goes along the same line as in \cite{Ni}.

\begin{prop}[cf. \cite{Ni}, (A10)]\label{prop:A10} 
Under the same notation as in Thm.\:\ref{thm:Nimmo}, we have
$$
\prod_{j=1}^n(1+\beta x_{j})^{-m}
\prod_{1\leq i<j\leq n}\dfrac{x_{i}- x_{j}}{x_{i}\oplus x_{j}}\times
GQ_\lambda(x_1,\ldots,x_n|b)
$$
$$
=\sum_{w\in \hat{S}_{n,m}}\sgn(w)
\times 
\det
\left(
[[x_{w(i)}|b]]^{\lambda_j}
(1+\beta x_{w(i)})^{j-m-1}\right)\times 
\prod_{m+1\leq i<j\leq n}\frac{x_{w(i)}- x_{w(j)}}{x_{w(i)}\oplus x_{w(j)}},
$$
where $\hat{S}_{n,m}$ denote the set of permutations $w$ 
of $\{1,\ldots,n\}$ 
satisfying $w(1)<\cdots<w(m),\;
w(m+1)<\cdots<w(n).$
The similar formula holds for $GP_\lambda(x_1,\ldots,x_n|b)$ 
when we use 
$[x|b]$ instead of $[[x|b]].$
\end{prop}

{\it Proof.} 
We can rewrite the definition of $GQ_\la(x_1,\ldots,x_n|b)$ as
follows:
\begin{equation}
\frac{1}{(n-m)!}
\sum_{w\in S_{n}}\prod_{i=1}^m
[[x_{w(i)}|b]]^{\lambda_i}
(1+\beta x_{w(i)})^{i-1} 
\prod_{i=1}^m
\prod_{j=i+1}^n
\frac{x_{w(i)}\oplus x_{w(j)}}{x_{w(i)}- x_{w(j)}}.\label{eq:HLtoNimmo}
\end{equation}
In fact,  if $n-r$ is even (so $m=r$), this is immediate
by using $x\ominus y=(x-y)(1+\beta y).$
If $n-r$ is odd and $m=r+1$, then we can deduce this expression by 
using Lemma\:\ref{var0-Iv} (the case of $k=0$ and the number of 
variables is odd).
Now we write the last factor as
$$
\prod_{i=1}^m
\prod_{j=i+1}^n
\frac{x_{w(i)}\oplus x_{w(j)}}{x_{w(i)}- x_{w(j)}}
=\prod_{1\leq i<j\leq n}\frac{x_{w(i)}\oplus x_{w(j)}}{x_{w(i)}- x_{w(j)}}
\times
\prod_{m+1\leq i<j\leq n}\frac{x_{w(i)}- x_{w(j)}}{x_{w(i)}\oplus x_{w(j)}}.
$$
Noting that 
$\prod_{j=1}^n(1+\beta x_{w(j)})^{m}=
\prod_{j=1}^n(1+\beta x_{j})^{m}$, which is invariant
under $S_n$,
and 
$$
\prod_{1\leq i<j\leq n}\frac{x_{w(i)}\oplus x_{w(j)}}{x_{w(i)}- x_{w(j)}}
=\mathrm{sgn}(w)
\prod_{1\leq i<j\leq n}\frac{x_{i}\oplus x_{j}}{x_{i}- x_{j}},
$$
we deduce that (\ref{eq:HLtoNimmo}) equals 
$$
\frac{1}{(n-m)!}\times\prod_{j=1}^n(1+\beta x_{j})^{m}
\prod_{1\leq i<j\leq n}\frac{x_{i}\oplus x_{j}}{x_{i}- x_{j}}
\sum_{w\in S_n}\sgn(w)
\times \prod_{i=1}^m
[[x_{w(i)}|b]]^{\lambda_i}
(1+\beta x_{w(i)})^{i-m-1}$$
\begin{equation}
\times 
\prod_{m+1\leq i<j\leq n}\frac{x_{w(i)}- x_{w(j)}}{x_{w(i)}\oplus x_{w(j)}}.\label{eq:second-factor}
\end{equation}
Now for each ${w}\in \hat{S}_{n,m}$, we perform the summation over 
the subgroup $G_{{w}}$ of $S_n$
stabilizing $\{x_{w(1)},\ldots,x_{w(m)}\}$ and
$\{x_{w(m+1)},\ldots,x_{w(n)}\}.$
Clearly $G_{{w}}\cong S_m\times S_{n-m}.$
We write $\sigma\in G_w$ as $\sigma_1\sigma_2$ 
correspondingly.
Summation over $\sigma_1$ creates 
$\det
\left(
[[x_{w(i)}|b]]^{\lambda_j}
(1+\beta x_{w(i)})^{j-m-1}\right)$, while 
$\sigma_2$ gives 
a factor $(n-m)!$ since 
the last factor in (\ref{eq:second-factor}) is anti-symmetric
in $n-m$ variables $\{x_{w(m+1)},\ldots,x_{w(n)}\}.$
This completes the proof.
$\qed$

\bigskip

We recall the following rule for expanding 
a Pfaffian of a matrix with blocks of zeros.

\begin{lem}%[cf. \cite{Ni}, (A8)]
Let $A$ be an $n\times n$
skew-symmetric matrix.
Let $m$ be a positive integer such that $m\leq n$ and 
$n-m$ is even. Let 
$B=(b_{ij})$ be an $n\times m$ matrix. Then
\begin{equation}
\mathrm{Pf}\left(\begin{array}{cc}
A & B\\
-{}^t\! B & 0 
\end{array}\right)
=(-1)^{\frac{m(m-1)}{2}}
\sum_{w\in\hat{S}_{n,m}}
\mathrm{sgn}(w)
\det(B_{w(1),\ldots,w(m);1,\ldots,m})
\mathrm{Pf}(A_{w(m+1),\ldots,w(n)}),\label{eq:Pf-exp}
\end{equation}
where $B_{w(1),\ldots,w(m);1,\ldots,m}$ denotes an $m\times m$ matrix
with $(i,j)$ entry 
$b_{w(i),j},$ and $A_{w(m+1),\ldots,w(n)}$ denotes
the submatrix of $A$ whose indices of rows and columns 
are $ w(m+1),\ldots,w(n).$
\end{lem}
{\it Proof.}
This formula has appeared in \cite{Ni}, (A8) . A proof can be found in \cite{IW}.
$\qed$

\bigskip

{\it Proof of Thm. \ref{thm:Nimmo}.} 
Since $n-m$ is even, we have by Lemma \ref{Pf-Schur}
$$
\prod_{m+1\leq i<j\leq n}\frac{x_{w(i)}- x_{w(j)}}{x_{w(i)}\oplus x_{w(j)}}
=\underset{1\leq i<j\leq n-m}{\mathrm{Pf}}\left(\frac{x_{w(m+i)}-x_{{w(m+j)}}}{x_{w(m+i)}\oplus x_{w(m+j)}}
\right).
$$
Then the right-hand side of the equation 
of Prop. \ref{prop:A10} becomes 
the form of
(\ref{eq:Pf-exp}).
By permuting  
columns of $B$, we can eliminate the
sign factor $(-1)^{\frac{m(m-1)}{2}}$ to obtain (\ref{eq:Nimmo}).
$\qed$

\subsection{Factorial Grothendieck polynomials}\label{subsec:Gro}

For a positive integer $n$,
let $\mathcal{P}_n$ denote the set 
of partitions of length less than or equal to $n;$ i.e. 
$\mathcal{P}_n$ is a set consisting of 
a sequence $\lambda=(\lambda_1,\ldots,\lambda_n)$ 
of non-negative integers 
such that 
$\lambda_1\geq \cdots\geq \lambda_n\geq 0.$
Let $\lambda\in \mathcal{P}_n.$ Define the 
following symmetric function
\begin{equation}
G_\lambda(x_1,\ldots,x_n|b)
=\frac{\det([x_i|b]^{\lambda_j+n-j}(1+\beta x_i)^{j-1})_{n\times n}}{
\prod_{1\leq i<j\leq n}(x_i-x_j)}.\label{eq:ratio}
\end{equation}
Using column operations to evaluate the determinant, we see that 
$G_\emptyset(x_1,\ldots,x_n|b)=1.$
\begin{rmk}\label{rmk:McNa}
McNamara \cite{Mc} introduced the {\it factorial Grothendieck polynomials\/} in terms of 
set-valued semistandard tableaux. 
It can be proved that (\ref{eq:ratio}) is 
identical to McNamara's function. In fact,
the argument in \S \ref{sec:comb_proof},
modified to type A case suitably, 
is applicable to show this fact. 
Note that we do not use this coincidence in this paper. 
For our present purpose, the expression (\ref{eq:ratio}) is useful. 
\end{rmk}

Consider a partition $\la\in \mathcal{P}_n.$
Define the sequence
$b_\la=(\ominus b_{\lambda_1+n},\ldots,
\ominus b_{\lambda_i+n-i+1},\ldots,\ominus b_{\lambda_n+1}).$
We identify a partition with its Young diagram
as usual (see \cite{Mac}).
\begin{prop}[cf. \cite{Mc}]\label{lem:vanishA}
Let $\lambda,\mu\in \mathcal{P}_n.$ Then
$$G_\lambda(b_\mu|b)=\begin{cases}
0
&\mbox{if}\; \lambda\not\subset \mu\\
\prod_{(i,j)\in \lambda}
(b_{n+j-\lambda_j'}\ominus b_{\lambda_i+n-i+1})
&\mbox{if} \;\lambda=\mu
\end{cases},
$$
where $\la_j'=\#\{i\;|\la_i\geq j\}.$  
\end{prop}
{\it Proof.} It is convenient to write 
the definition of $G_\lambda(x|b)$ in the following form:
\begin{equation}
G_\lambda(x_1,\ldots,x_n|b)=
\displaystyle\sum_{w\in S_n} w
\left[
\frac{[x|b]^{\lambda+\rho_{n-1}}}{\prod_{1\leq i<j\leq n}(x_i\ominus x_j)}
\right],\label{def:Gro}
\end{equation}
where $\rho_{n-1}=(n-1,\ldots,2,1,0).$
Then the proposition is proved by 
straightforward calculations (cf. \cite{Iv2}, \cite{MS}).
$\qed$

\bigskip

Let $\R$ be the localization of $\Z[\beta][b_1,b_2,\ldots]$ 
by the multiplicative system 
formed by products of $1+\beta b_i\;(i\geq 1).$
Note that $G_\lambda(b_\mu|b)$ belongs to $\R.$
The following lemma was proved in \cite{Mc}. 
As our definition of $G_\lambda(x|b)$ is
different from the one in \cite{Mc},  
we give a proof here 
for completeness.

\begin{lem}[cf. \cite{Mc}]\label{G_basis}
$G_\lambda(x_1,\ldots,x_n|b)\;(\lambda\in \mathcal{P}_n)$ form an $\R$-basis
of $\R[x_1,\ldots,x_n]^{S_n}.$
\end{lem}
{\it Proof.} 
Let $m$ be a positive integer.
Let $\mathcal{P}_{n,m}$ be the set of partitions
$\lambda$ in $\mathcal{P}_n$ such that $\lambda_1\leq m.$
Let $L_{n,m}$ be the $\Z[\beta,b]$--span of 
$m_\mu(x)\;(\mu\in \mathcal{P}_{n,m}),$
where $m_\mu(x)$ is the monomial 
symmetric function in variables $x_1,\ldots,x_n$
corresponding to $\mu$ (cf. \cite{Mac}, I, 2).
Note that the highest possible 
power of each $x_i$ 
in $G_\lambda(x_1,\ldots,x_n|b)$ 
is $m$
(this fact
can be seen from (\ref{eq:ratio})).
Thus we can define the
elements $d_{\lambda\mu}(\beta,b)
\in \Z[\beta,b]$ by 
\begin{equation}
G_\lambda(x_1,\ldots,x_n|b)=\sum_{\mu\in\mathcal{P}_{n,m} }
d_{\lambda\mu}(\beta,b)\,
m_\mu(x)\quad(\lambda\in \mathcal{P}_{n,m}).\label{eq:linear}
\end{equation}
We claim that 
$\det(d_{\lambda\mu}(\beta,b))$ is a product of factors 
of the form $1+\beta b_i$ for some $i.$  
Once the claim is proved, we can
invert the system of linear equations (\ref{eq:linear}) over the ring 
$\R=\Z[\beta,b][(1+\beta b_i)^{-1}(i\geq 1)]$
to express $m_\la(x)$' s as $\R$-linear combination of 
$G_\lambda(x_1,\ldots,x_n|b)$'s with $\la\in \mathcal{P}_{n,m}$.
Since each element in $\R[x_1,\ldots,x_n]^{S_n}$
is in $\R\otimes_{\Z[\beta,b]}L_{n,m}$ for some $m$, we 
have the lemma.
 
In order to prove the claim,
we first show that $\det(d_{\lambda\mu}(0,0))=1.$
Now we specialize $\beta$ and $b$ to $0$ then 
$
G_\lambda(x_1,\ldots,x_n|b)$ becomes 
the classical Schur function $s_{\lambda}(x_1,\ldots,x_n).
$
It follows that, $d_{\lambda\mu}(0,0)$ 
is equal to the Kostka number (cf. \cite{Mac}, I. 6) if $|\lambda|=|\mu|.$
Moreover, since $s_{\lambda}(x_1,\ldots,x_n)$ is a homogeneous polynomial
of degree $|\lambda|$, 
we have  $d_{\lambda\mu}(0,0)=0$ 
if $|\lambda|\neq |\mu|.$ 
Therefore the matrix $(d_{\lambda\mu}(0,0))_{\la,\mu}$ is lower triangular (with respect to 
a linear order which is a refinement of the dominance order) and 
all the entries of the main diagonal are $1$. 
This implies $\det(d_{\lambda\mu}(0,0))=1.$ 

Now using Lemma \ref{lem:vanishA},
it can be proved that the only possible irreducible factors of 
$\det(d_{\la\mu}(\beta,b))$ are of the form $1+\beta b_i$ 
(see \cite{Mc}, Lemma 4.7 and an argument after the lemma).
Because we have $\det(d_{\la\mu}(0,0))=1$, 
we can conclude that $\det(d_{\la\mu}(\beta,b))$
is actually a product of $1+\beta b_i$ for some $i.$
\qed

\bigskip

As a by-product of the proof of the preceding lemma, we have
the following.

\begin{cor}\label{cor:basis} 
$G_\lambda(x_1,\ldots,x_n)\;(\lambda\in \mathcal{P}_n)$ form a $\Z[\beta]$-basis
of $\Z[\beta][x_1,\ldots,x_n]^{S_n}.$
\end{cor}
{\it Proof.} Indeed, by letting $b_i=0$ for all $i$, 
the proof of 
Lemma \ref{G_basis} works when we consider $\Z[\beta]$ instead of 
$\R.$ 
$\qed$

\subsection{Factorization formula}\label{sec:Fac}
The next result is an analogue of the factorization formula of Stanley 
(\cite{Cont}, Problem 4 in Problem Session).
\begin{prop}\label{prop:factor}
For every $k\geq 1$ let $\rho_k$ denote the partition $(k,k-1,\ldots,2,1).$

(1) For positive integer $n$, we have 
\begin{eqnarray}
GP_{\rho_{n-1}}(x_1,\ldots,x_n|b)
&=&
GP_{\rho_{n-1}}(x_1,\ldots,x_n)
=\prod_{1\leq i<j\leq n}
(x_i\oplus x_j),\\
GQ_{\rho_n}(x_1,\ldots,x_n|b)
&=&GQ_{\rho_{n}}(x_1,\ldots,x_n)
=\prod_{1\leq i\leq j\leq n}(x_i\oplus x_j).
\label{eq:triangle}
\end{eqnarray}

 (2) For positive integer $n$ and  
$\lambda\in \mathcal{P}_n$, 
we have
\begin{eqnarray}
GP_{\rho_{n-1}+\lambda}
(x_1,\ldots,x_n|b)
&=&GP_{\rho_{n-1}}(x_1,\ldots,x_n)
G_\lambda(x_1,\ldots,x_n|b),\\
GQ_{\rho_n+\lambda}(x_1,\ldots,x_n|b)
&=&GQ_{\rho_n}(x_1,\ldots,x_n)
G_{\lambda}(x_1,\ldots,x_n|b).
\label{eq:fac}
\end{eqnarray}

\end{prop}
{\it Proof.} 
We only prove the $GP_\lambda$ case. $GQ_\lambda$ case is similar.
Note that 
the length of $\rho_{n-1}+\lambda$ is $n-1$ or $n.$
In each cases,
we symmetrize the following rational function
$$
[x|b]^{\lambda+\rho_{n-1}}
\prod_{1\leq  i<j\leq n}
\frac{x_i\oplus x_j}{x_i\ominus x_j}
=
\prod_{1\leq i<j\leq n}
(x_i\oplus x_j)\times
\frac{[x|b]^{\lambda+\rho_{n-1}}}
{\prod_{1\leq i<j\leq n}(x_i\ominus x_j)}
$$
Since 
$
\prod_{1\leq i<j\leq n}
(x_i\oplus x_j)
$ is symmetric, we have (see (\ref{def:Gro}))
\begin{equation}
GP_{\rho_{n-1}+\lambda}
(x_1,\ldots,x_n|b)
=\prod_{1\leq  i<j\leq n}
(x_i\oplus x_j)\times 
G_\lambda(x_1,\ldots,x_n|b).\label{eq:facwa}
\end{equation}
In order to see $GP_{\rho_{n-1}}(x_1,\ldots,x_n|b)=\prod_{1\leq i<j\leq n}
(x_i\oplus x_j),$
we set $\lambda=\emptyset$ 
in (\ref{eq:facwa}). 
\qed

\setcounter{equation}{0}
\section{$K$-theoretic $Q$-cancellation property}\label{sec:beta-ss}
In this section, we 
 prove that the $K$-theoretic $P$-functions 
form a basis of the ring of functions with the $K$-theoretic 
$Q$-cancellation property (Thm. \ref{thm:basis}).
Here we mainly consider ``non-equivariant'' versions
$GQ_\lambda(x_1,\ldots,x_n)$ and
$GP_\lambda(x_1,\ldots,x_n).$

\subsection{Polynomials with $K$-theoretic $Q$-cancellation property}

The polynomials with the $K$-theoretic $Q$-cancellation property (Def. \ref{def:Ksuper})
form a subring in $\Z[\beta][x_1,\ldots,x_n].$
We denote it by 
$G\Gamma_n.$

\begin{prop}
$GP_\lambda(x_1,\ldots,x_n)$ and 
$GQ_\lambda(x_1,\ldots,x_n)$ are in $G\Gamma_n.$
\end{prop}

This is a consequence of the next Lemma which was proved in \cite{Iv2} for the case of $\beta=0.$

 \begin{lem}\label{lem:Iv} Let $r\leq n$, and $f(x_1,\ldots,x_r)$ be a 
polynomial in $\Z[\beta][x_1,\ldots,x_r].$
Then the following function belongs to $G\Gamma_n$:
$$R_n(x_1,\ldots,x_n)=
\sum_{w\in S_n}
f(x_{w(1)},\ldots,x_{w(r)})
\prod_{i=1}^r\prod_{j=i+1}^n \frac{x_{w(i)}\oplus x_{w(j)}}{x_{w(i)}\ominus x_{w(j)}}.$$
If moreover $f$ 
 is divisible by $x_1\cdots x_r$ then 
 $\tilde{R}_n(x_1,\ldots,x_n)=\frac{1}{(n-r)!} R_n(x_1,\ldots,x_n)$
 has stability. i.e.
 $\tilde{R}_{n+1}(x_1,\ldots,x_n,0)
 =\tilde{R}_n(x_1,\ldots,x_n).$
\end{lem} 
{\it Proof.}
Let $x_i=t,\;x_j=\ominus t$, for arbitrary integers $i,j$ such that $1\leq i<j\leq n.$
We claim that each summand corresponding to $w\in S_n$ 
does not depend on $t.$
To see this, set 
 $F_n(x_1,\ldots,x_n)=
 \prod_{i=1}^r\prod_{j=i+1}^n\dfrac{x_i\oplus x_{j}}{x_{i}\ominus x_{j}}.$
In fact, if one of $i,j$ is in $\{w(1),\ldots,w(r)\}$, then
$F_n(x_{w(1)},\ldots,x_{w(n)})$
 vanishes.
If $i,j\in \{w(r+1),\ldots,w(n)\},$ then by
using identity
$$
\frac{x\oplus t}{x\ominus t}
\frac{x\ominus t}{x\oplus t}=1,
$$
we see that $F_n(x_{w(1)},\ldots,x_{w(n)})$ does not depend on $t$.
Also $f(x_{w(1)},\ldots,x_{w(r)})$ is obviously independent of $t.$
Hence the claim follows. 

Next we prove the second assertion.
Let $x_{n+1}=0.$ 
If $n+1$ is in $\{w(1),\ldots,w(r)\}$, then the corresponding term vanishes because 
$f(x_{w(1)},\ldots,x_{w(r)})=0$ by the assumption. 
If $n+1\in \{w(r+1),\ldots,w(n)\},$ then 
$F_{n+1}(x_{w(1)},\ldots,x_{w(n+1)})=
F_n(x_{w'(1)},\ldots,x_{w'(n)})$
where $w'$ is a permutation obtained by multiplying the transposition $(w(n+1),n+1)$
from the left to $w.$
This implies the desired result. 
$\qed$

\begin{thm}[Basis theorem]\label{thm:basis} The polynomials 
$GP_\lambda(x_1,\ldots,x_n)\,(\lambda\in \mathcal{SP}_n)$
form a $\Z[\beta]$ basis of $G\Gamma_n.$
\end{thm}

\subsection{Proof of Basis Theorem}\label{ssec:PfBasis}

We use the same strategy of the proof as in \cite{Pr},
with the aid of the following lemma.
\begin{lem}\label{lem:2var}
Let $A$ be any unique factorization
domain,  
$\beta,x,y$  
be independent variables over $A.$
If $f(\beta,x,y)\in A[\beta,x,y]$ vanishes when we made the 
substitution $x=\ominus y,$
then $f(\beta,x,y)$ is divisible by $x\oplus y.$ 
\end{lem}
{\it Proof.} 
We can prove this lemma by applying the division algorithm 
for $f(\beta,x,y)$ as a polynomial in $x.$
$\qed$

\bigskip

{\it Proof of Thm. \ref{thm:basis}.}
We first consider the case when $n$ is an even integer.
We use induction on $n$.
Let $n=2$ and let 
$f(x_1,x_2)\in G\Gamma_2$.
We may assume the constant term
of $f$ is zero
so that $f(0,0)=0.$
Since $f(x_1,x_2)\in G\Gamma_2$, 
we have $f(t,\ominus t)=f(0,0)=0.$ 
Then by Lemma \ref{lem:2var},
$f(x_1,x_2)$ is divisible by 
$(x_1\oplus x_2).$
Thus $f$ can be 
written as $f(x_1,x_2)=(x_1\oplus x_2)g(x_1,x_2).$
Since $g(x_1,x_2)$ is a symmetric polynomial, 
we can write $g(x_1,x_2)=\sum_{\lambda\in \mathcal{P}_2}c_\lambda
G_{\lambda}(x_1,x_2),\; c_\lambda\in \Z[\beta].$
By the factorization formula, we have
$$
f(x_1,x_2)=\sum_{\lambda\in \mathcal{P}_2}c_\lambda
GP_{\lambda+\rho_2}(x_1,x_2).
$$
Thus $f(x_1,x_2)$ is a $\Z[\beta]$-linear combination of
$GP_{\lambda}(x_1,x_2),\;\lambda\in\mathcal{SP}_2.$

If $n\geq 4.$ We proceed as follows.
$f(x_1,\ldots,x_{n-2},0,0)\in G\Gamma_{n-2}$.
By induction,
$f'=f(x_1,\ldots,x_{n-2},0,0)$
is a linear combination of
$GP_\lambda(x_1,\ldots,x_{n-2})$'s where
$\lambda\in \mathcal{SP}_{n-2}.$ 
Let 
$$
f'
=\sum_{\lambda} c_\lambda \;GP_\lambda(x_1,\ldots,x_{n-2}),\quad
c_\lambda\in \Z[\beta].
$$
Consider the polynomial
$g(x)=\sum_{\lambda} c_\lambda \;GP_\lambda(x_1,\ldots,x_n).$
Let $h(x)=f(x)-g(x).$
We have
$
h(x_1,\ldots,x_{n-2},t,\ominus t)
=h(x_1,\ldots,x_{n-2},0,0)=0.$
This implies that 
$x_{n-1}\oplus x_n$ divides $h(x).$
Since $h(x)$ is symmetric, $h(x)$
is a multiple of 
$GP_{\rho_{n-1}}(x_1,\ldots,x_n)=\prod_{1\leq i<j\leq n}(x_{i}\oplus x_{j}).$ 
Thus
$$
f=g+GP_{\rho_{n-1}}(x_1,\ldots,x_n)s(x)
$$
where $s(x)\in \Z[\beta][x_1,\ldots,x_n]^{S_n}.$
Write by Cor. \ref{cor:basis} 
$$
s(x)=\sum_{\lambda \in \mathcal{P}(n)}c_\lambda'
G_\lambda(x_1,\ldots,x_n),\quad
c'_\lambda\in \Z[\beta].
$$
Hence we have
\begin{equation}
f(x)=
\sum_{\lambda\in \mathcal{SP}_{n-2}} c_\lambda \;GP_\lambda(x_1,\ldots,x_n)
+\sum_{\lambda \in \mathcal{P}_n}c_\lambda'
GP_{\rho_{n-1}+\lambda}(x_1,\ldots,x_n).\label{f(x)}
\end{equation}
 Note that if $\lambda\in \mathcal{P}_n$ then 
$\rho_{n-1}+\lambda\in \mathcal{SP}_n.$ 
This completes the proof.

Next consider the case when $n$ is odd. 
Let $n=1.$
This case is obvious since 
$G\Gamma_1=\Z[\beta][x_1]$ by definition
and $GP_m(x_1)=x_1^m.$
If $n\geq 3$ the proof is the same 
as the even case. 
$\qed$

\subsection{Characterization of subring spanned by $GQ_\lambda(x)$}\label{Gamma+}
Let $G\Gamma_{n,+}$ be the set consisting 
of $F\in G\Gamma_n$ such that  
$F(t,x_2,\ldots,x_n)-F(0,x_2,\ldots,x_n)$ is divisible by $t\oplus t.$
Obviously $G\Gamma_{n,+}$ is a subring of  $G\Gamma_{n}.$

\begin{prop}[Basis theorem for $G\Gamma_{n,+}$]\label{prop:basisGQ}
 $G\Gamma_{n,+}=\bigoplus_{\lambda\in \mathcal{SP}_n}\Z[\beta]
 GQ_\lambda(x_1,\ldots,x_n).$ 
\end{prop}
{\it Proof.} We know that $GQ_\lambda(x_1,\ldots,x_n)$ belongs to $G\Gamma_{n}$. By definition, one see that 
$GQ_\lambda(t,x_2,\ldots,x_n)-GQ_\lambda(0,x_2,\ldots,x_n)$
is divisible by $t\oplus t.$ So $GQ_\lambda(x_1,\ldots,x_n)$ is 
an element of $G\Gamma_{n,+}.$
Let $F\in G\Gamma_{n,+}.$ 
We will show that $F$ is a $\Z[\beta]$-linear combination of 
$GQ_\lambda(x_1,\ldots,x_n)\;(\lambda\in \mathcal{SP}_n)$
by induction on $n.$
We can proceed 
in the same way as in the proof for Thm.\:\ref{thm:basis}, 
by using $GQ_\lambda(x_1,\ldots,x_n)$ instead of 
$GP_\lambda(x_1,\ldots,x_n).$
 $\qed$

\subsection{Inverse limit}
As with the other types of symmetric functions, 
our polynomials have the following stability property.
 
  \begin{prop}\label{prop:stable} Let $\lambda$ be a strict partition 
  of length $r$ 
  in $\mathcal{SP}_n.$ Then

(1) $GQ_\lambda(x_1,\ldots,x_{n-1},0)
 =GQ_\lambda(x_1,\ldots,x_{n-1}),$
 
(2) $GP_\lambda(x_1,\ldots,x_{n-1},0)
 =GP_\lambda(x_1,\ldots,x_{n-1}),$

where the right-hand sides are zero if $r=n.$
 \end{prop}
{\it Proof.} We can apply 
Lemma \ref{lem:Iv}.
$\qed$
   
\begin{rmk}\label{rmk:stab}
By a similar argument as
in the proof of Lemma \ref{lem:Iv}, one sees that $$GP_\lambda(x_1,\ldots,x_{n-2},0,0|b)
 =GP_\lambda(x_1,\ldots,x_{n-2}|b).$$
However, the property 
$GP_\lambda(x_1,\ldots,x_{n-1},0|b)
 =GP_\lambda(x_1,\ldots,x_{n-1}|b)$
 does not hold in general; for example, we have $GP_{1}(x_1,x_2)=x_1\oplus x_2,$
 whereas  
 $GP_1(x_1)=x_1\oplus b_1.$
In contrast to this fact, we have
 $GQ_\lambda(x_1,\ldots,x_{n-1},0|b)
 =GQ_\lambda(x_1,\ldots,x_{n-1}|b),$
 which holds since $[[x|b]]^\la$ is divisible by $x_1\cdots x_r$, whereas 
$[x|b]^\la$ is not (see Lemma \ref{lem:Iv}).
\end{rmk}

Let $\varphi_{n+1}: G\Gamma_{n+1}\rightarrow G\Gamma_n$
be the morphism of the $\Z[\beta]$-algebras
given by the specialization $x_{n+1}=0.$ 
Then 
$\{G\Gamma_n,\varphi_n\}$ form an inverse system.
Let $G\Gamma$ denote the inverse limit
$\varprojlim G\Gamma_n.$
%We call this the {\it ring of $K$-supersymmetric functions}.
Then, by the stability property of $GP_\lambda(x_1,\ldots,x_n)$ (Prop. \ref{prop:stable}), 
we have $GP_\lambda(x):=\varprojlim GP_\lambda(x_1,\ldots,x_n)
\in G\Gamma.$ 
Similarly we define $G\Gamma_{+}=\varprojlim G\Gamma_{n,+},$
which is a subring of $G\Gamma.$
Then $GQ_\lambda(x):=\varprojlim GQ_\lambda(x_1,\ldots,x_n)$ in
$G\Gamma_+.$ 
\bigskip

\begin{rmk}\label{rmk:beta-exp}
Both the functions $GP_\lambda(x)$ and $GQ_\lambda(x)$ 
can be expressed in the following form:
\begin{equation}
F_{|\lambda|}(x)+\beta F_{|\lambda|+1}(x)
+\beta^2 F_{|\lambda|+2}(x)+\cdots,\quad
F_{k}(x) \in \Lambda^{k},
\end{equation}
where $\Lambda^k\;(k\geq 0)$ is the space of homogeneous 
symmetric functions of degree $k$ (see \cite{Mac}).
The initial terms $F_{|\lambda|}(x)$ are $P_\lambda(x)$ and $Q_\lambda(x)$
respectively. 
In particular, we see that the functions $GP_\lambda(x)$ (resp. $GP_\lambda(x)$) $,\la\in \mathcal{SP},$ are
linearly independent over $\Z[\beta].$
\end{rmk}

\begin{rmk}
From combinatorial results proved in \S \ref{sec:comb}, 
one can see that 
all the homogeneous parts $F_{k}(x)$ are actually non-zero.
Moreover, each $F_{k}(x)$ is 
a non-negative linear combination of monomial symmetric functions
 $m_\mu(x)$ such that $|\mu|=k.$
\end{rmk}

\bigskip

Here is the `basis' theorem for $G\Gamma.$
\begin{prop}\label{prop:BasisGP}
Any $f(x)\in G\Gamma$  
can be expressed uniquely as 
a possibly infinite $\Z[\beta]$-linear combination of $GP_\lambda$'s, 
\begin{equation}
f(x)=\sum_{\lambda\in \mathcal{SP}}c_\lambda\cdot GP_\lambda(x),\quad
c_\lambda\in \Z[\beta],
\end{equation}
such that for all positive integer $n$, the set 
$\{\lambda\in \mathcal{SP}_n\;|\;c_\lambda\neq 0\}$
is finite.
\end{prop}

{\it Proof.} This is a direct consequence of Thm.\:\ref{thm:basis}. 
$\qed$

\bigskip

We can introduce $c_{\lambda\mu}^\nu\in \Z[\beta]$
by the following expansion with possibly infinitely many 
terms: 
$$
GP_\lambda(x)\cdot GP_\mu(x)
=\sum_{\nu}c_{\lambda\mu}^\nu
GP_\mu(x).
$$

We expect that the right-hand side is actually a finite sum. 
For type $A$, the corresponding statement is true
as a consequence of the 
explicit description of LR-coefficients
given there (\cite{Buch}, Cor. 5.5).
\begin{conj}\label{conj:GP}
$\bigoplus_{\lambda\in\mathcal{SP}}\Z[\beta]GP_\lambda(x)$ is 
a subring of $G\Gamma.$
\end{conj}

\begin{rmk}
Our result (Thm. \ref{thm:non-eq}) shows that 
the constants $c_{\lambda\mu}^\nu$, when
specialized to $\beta=-1$, are 
equal to $K$-theory Littlewood-Richardson coefficients
of the maximal orthogonal Grassmannians.
Clifford, Thomas and Yong \cite{CTY}
have given a combinatorial interpretation for the LR-coefficients.
%See the note at the end of \S 1. 
We do not know the description in \cite{CTY}
implies Conjecture \ref{conj:GP}.
\end{rmk}

\bigskip

The algebra $G\Gamma$ has a natural decreasing  
filtration $G\Gamma=F^0\supset F^1\supset \cdots \supset F^k\supset \cdots,$ 
defined by 
$$
F^k=F^k G\Gamma
=\left.\left\{\sum_{\lambda\in \mathcal{SP}}c_\lambda \cdot GP_\lambda(x)\in G\Gamma\;\right|\;
c_\lambda=0\;\mbox{if}\;|\lambda|<k\right\}.
$$
By proposition \ref{prop:BasisGP} and a consideration of degree in $x_i$, one sees that 
$F^k$ is actually an ideal and $F^k\cdot F^l \subset F^{k+l}.$ 
The associated graded ring $\mathrm{gr}\,G\Gamma
:=\bigoplus_{k=0}^\infty F^k/F^{k+1}
$ is naturally isomorphic to 
$\Z[\beta]\otimes_\Z\Gamma,$
where $\Gamma$ denote the ring of Schur $P$-functions,
the $\Z$-span of the $P_\lambda(x)\;(\lambda\in \mathcal{SP}).$
This fact follows from Remark \ref{rmk:beta-exp}.

\bigskip

We have a similar result for $G\Gamma_+.$ This is 
a consequence of Prop. \ref{prop:basisGQ}.  
\begin{prop}\label{prop:BasisGQ}
Any $f(x)\in G\Gamma_+$  
can be expressed uniquely as 
a possibly infinite $\Z[\beta]$-linear combination of $GQ_\lambda$'s, 
\begin{equation}
f(x)=\sum_{\lambda\in \mathcal{SP}}c_\lambda\cdot GQ_\lambda(x),\quad
c_\lambda\in \Z[\beta],
\end{equation}
such that for all positive integers $n$, the set 
$\{\lambda\in \mathcal{SP}_n\;|\;c_\lambda\neq 0\}$
is finite.
\end{prop}

\begin{conj}
$\bigoplus_{\lambda\in\mathcal{SP}}\Z[\beta]GQ_\lambda(x)$ is 
a subring of $G\Gamma_+.$
\end{conj}

\begin{rmk}
The Pieri type rule in \cite{INN} implies that 
$GQ_\la(x)GQ_k(x)$ is a finite $\Z[\beta]$ linear 
combination of $GQ_\mu(x)$'s. 
\end{rmk}

\subsection{Expansion of $GP_\la(x|b)$ in terms of $GP_\la(x)$}\label{ssec:ExpGP}
We define that the equivariant version of $GP_\la(x)$ is defined 
by
$$
GP_\la(x|b)=\varprojlim GP_{\la}(x_1,\ldots,x_{2m}|b).
$$
Here we use only the polynomials with even number of variables (see Remark \ref{rmk:stab}).
We have the next results.
\begin{prop}\label{prop:GPexpand}
Let $r$ be the length of $\la.$
Let 
$\la\in\mathcal{SP},$ and $n$ be the smallest even integer
such that $r\leq n.$
Let $d_{\la\mu}\in \Z[\beta,b]$ be defined by the following equation:
\begin{equation}
GP_\la(x_1,\ldots,x_n|b)
=\sum_{\mu\in \mathcal{SP}_n}d_{\la\mu}GP_\mu(x_1,\ldots,x_n).\label{eq:exGP}\end{equation}
Then the following equation holds 
in $\Z[b]\otimes G\Gamma:$
\begin{equation}
GP_\la(x|b)
=\sum_{\mu\in \mathcal{SP}_n}d_{\la\mu}GP_\mu(x).\label{eq:GP_expand}
\end{equation}
\end{prop}
{\it Proof.} We first claim that $d_{\la\mu}\neq 0$ implies that 
the length of $\mu$ is $r$, if $r$ is odd, and $r$ or $r-1$, if $r$ is even.
In fact, by the factorization theorem, in either cases,  we have % if $r$ is odd, then $n=r+1$ and 
$GP_\la(x_1,\ldots,x_n|b)=GP_{\rho_{n-1}}(x_1,\ldots,x_{n-1})
G_\mu(x_1,\ldots,x_n|b)$
with $\mu\in \mathcal{P}_r.$
Since $G_\mu(x_1,\ldots,x_n|b)$ is a $\Z[\beta,b]$-linear combination 
of $G_\nu(x_1,\ldots,x_n)$'s such that $\mu\in \mathcal{P}_r,$
the claim follows, again by the factorization theorem.

%If $r$ is even, then $n=r$ and 
%$GP_\la(x_1,\ldots,x_n|b)
%=GP_{\rho_{n-1}}(x_1,\ldots,x_{n-1})
%G_\mu(x_1,\ldots,x_n|b)$ in (\ref{}) is 
%\begin{lem}
%Let 
%$\la\in\mathcal{SP}$ and $n$ be the smallest even integer
%such that $r\leq n.$
%$GP_\la(x_1,\ldots,x_n|b)
%=\sum_{\mu}d_{\la\mu}GP_\mu(x_1,\ldots,x_n).$
%
%\end{lem}
%{\it Proof.} 
%Let $n$ be the smallest even integer greater than $r,$
%the length of $\la.$
%By assumption and the above claim, 
%the polynomial $$GP_\la(x_1,\ldots,x_n|b)
%=\frac{1}{(n-r)!}\sum_{w\in S_{n}}w\left([x|b]^\la\Phi_{n,r}(x)\right)
%$$ is equal to 
%

By Lemma \ref{lem:big-r} and the claim,
we see that each polynomial $GP_\mu(x_1,\ldots,x_n)$
appearing in the right hand side of (\ref{eq:exGP})
is equal to 
\begin{equation}
\frac{1}{(n-r)!}\sum_{w\in S_{n}}w\left([x|b]^\mu\Phi_{n,r}(x)\right).\label{eq:r-sum}
\end{equation}
%$$
%%\frac{1}{(n-r)!}\sum_{w\in S_{n}}w\left([x|b]^\la\Phi_{r,n}(x)\right)
%\frac{1}{(n-r)!}
%\sum_{\mu}d_{\la\mu}
%\sum_{w\in S_{n}}w\left(x^\mu\Phi_{n,r}(x)\right).
%%+\frac{1}{(n-r+1)!}
%%\sum_{l(\mu)=r-1}d_{\la\mu}
%%\sum_{w\in S_{n}}w\left(x^\mu\Phi_{n,r-1}(x)\right).
%$$
Now we will calculate 
$$GP_\la(x_1,\ldots,x_{n+2}|b)
=\frac{1}{(n+2-r)!}\sum_{w\in S_{n+2}}w\left([x|b]^\la\Phi_{n+2,r}(x)\right).
$$
Since $$
\Phi_{n+2,r}(x)=
\Phi_{n,r}(x)
\prod_{i=1}^r\frac{x_i\oplus x_{n+1}}{x_i\ominus x_{n+1}}
\frac{x_i\oplus x_{n+2}}{x_i\ominus x_{n+2}},
$$
%we can deduce 
%$$
%\frac{1}{(n+2-r)!}
%\sum_{w\in S_{n+2}}w\left([x|b]^\la\Phi_{n+2,r}(x)\right)
%\prod_{i=1}^r\frac{x_i\oplus x_{w(n+1)}}{x_i\ominus x_{w(n+1)}}
%\frac{x_i\oplus x_{w(n+2)}}{x_i\ominus x_{w(n+2)}},
%$$
%$$
%=\frac{1}{(n+2-r)!}
%\sum_{\mu}d_{\la\mu}
%\sum_{w}w\left(x^\mu\Phi_{n,r}(x)\right)
%\prod_{i=1}^r\frac{x_i\oplus x_{w(n+1)}}{x_i\ominus x_{w(n+1)}}
%\frac{x_i\oplus x_{w(n+2)}}{x_i\ominus x_{w(n+2)}},
%$$
%where $w$ permutes $\{1,\ldots,n\}-\{w(n+1),w(n+2)\}.$
%%$$
%%GP_\la(x_{w(1)},\ldots,x_{w(n)})
%%\prod_{i=1}^r\frac{x_i\oplus x_{w(n+1)}}{x_i\ominus x_{w(n+1)}}
%%\frac{x_i\oplus x_{w(n+2)}}{x_i\ominus x_{w(n+2)}},
%%$$
%$$
%=\sum_{\mu}d_{\la\mu}
%\frac{1}{(n+2-r)!}\sum_{w\in S_{n+2}}w\left(x^\mu\Phi_{n+2,r}(x)\right).
%$$
we can 
easily derive the following equation from (\ref{eq:exGP}) and (\ref{eq:r-sum})
$$GP_{\mu}(x_1,\ldots,x_{n+2}|b)
=\sum_{\mu}d_{\la\mu}
GP_{\mu}(x_1,\ldots,x_{n+2}).
$$
Hence, by passing to the limit, we have (\ref{eq:GP_expand}).$\qed$

\setcounter{equation}{0}
\section{Preliminaries on Weyl groups and root systems}\label{sec:Weyl}

\subsection{Weyl groups}
Let $X=B,C,D$ and $W=W(X_\infty)$ be the Weyl group of type $X_\infty.$
This is a Coxeter group whose generators $\{s_i\}_{i\in I}$ are indexed by   
$I=\{0,1,2,\ldots\}$ for $B_\infty,\;C_\infty$ and
$I=\{\hat{1},1,2,\ldots\}$ for $D_\infty.$
If $W$ is type $C_\infty$ ($B_\infty$) then the defining relations are  
$s_i^2=1\;(i=0,1,2\ldots)$ and 
\begin{eqnarray}
s_is_{i+1}s_i&=&s_{i+1}s_is_{i+1}\quad (i\geq 1),\label{WA1}\\
s_is_j&=&s_js_i\quad (i,j\geq 1,\;|i-j|\geq 2),\label{WA2}\\
s_0s_1s_0s_1&=&s_1s_0s_1s_0,\quad 
s_0s_i=s_is_0\quad (i\geq 1).\label{WA3}
\end{eqnarray}
If $W$ is type $D_\infty$ then the defining relations are  
$s_i^2=1\;(i=\hat{1},1,2\ldots)$, (\ref{WA1}), (\ref{WA2}), and 
\begin{eqnarray}
s_{\hat{1}}s_2s_{\hat{1}}&=&s_2s_{\hat{1}}s_2,\quad
s_{\hat{1}}s_j=s_j s_{\hat{1}}\quad (j\neq 2).
\end{eqnarray}
The corresponding graphs are given as follows:
\bigskip
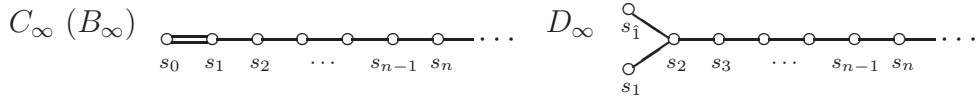
\setlength{\unitlength}{0.4mm}
\begin{center}
  \begin{picture}(150,25)
  \thicklines
  \put(-50,20){${C}_{\infty}$ (${B}_{\infty}$)}
  \put(0,10)
  {\put(0,5){$\circ$}
  \put(4,6.5){\line(1,0){12}}
  \put(4,8.5){\line(1,0){12}}
  \multiput(15,5)(15,0){5}{
  \put(0,0){$\circ$}
  \put(4,2.4){\line(1,0){12}}}
  \put(90,5){$\circ$}}
  \put(0,8){\tiny{$s_0$}}
  \put(15,8){\tiny{$s_1$}}
  \put(30,8){\tiny{$s_2$}}
  \put(70,8){\tiny{$s_{n-1}$}}
  \put(90,8){\tiny{$s_{n}$}}
  \put(50,8){\tiny{$\cdots$}}
  \put(94.5,17.5){\line(1,0){10}}
  \put(106,15){$\cdots$}
  \end{picture}
 \begin{picture}(150,20)
  \thicklines
  \put(-25,20){${D}_{\infty}$}
   \put(0,0)
  {\put(0,25){$\circ$}
  \put(0,5){$\circ$}
  \put(4,26){\line(3,-2){12}}
  \put(4,8.5){\line(3,2){12}}
  \multiput(15,15)(15,0){4}{
  \put(0,0){$\circ$}
  \put(4,2.4){\line(1,0){12}}}
  \put(75,15){$\circ$}
  \put(120,10)}
 \put(0,20){\tiny{$s_{\hat{1}}$}}
  \put(0,0){\tiny{$s_1$}}
  \put(15,8){\tiny{$s_2$}}
  \put(30,8){\tiny{$s_3$}}
  \put(70,8){\tiny{$s_{n-1}$}}
  \put(90,8){\tiny{$s_{n}$}}
  \put(50,8){\tiny{$\cdots$}}
 \put(79,17.5){\line(1,0){12}}
  \put(90,15){$\circ$}
    \put(94.5,17.5){\line(1,0){10}}
  \put(106,15){$\cdots$}
  \end{picture}
\end{center}

Let $I_n=\{0,1,\ldots,n-1\}$ for $X=B,C$ 
and $\{\hat{1},1,2,\ldots,n-1\}$ for $X=D.$
Let $W(X_n)$
denote the subgroups of $W(X_\infty)$
generated by $s_i\;(i\in I_n).$ 

\subsection{Weyl groups and signed permutations}
Let $\mathbb{N}$ be the set of positive integers $\{1,2,\ldots\}.$
Denote by $\overline{\mathbb{N}}$ a `negative' 
copy $\overline{\mathbb{N}}=  \{\bar{1},\bar{2},\ldots\}$
of $\mathbb{N}.$ A {\it signed permutation\/} $w$ of $\mathbb{N}$, by definition, is
a bijection on the set $\N\cup \overline{\N}$ such that 
$\overline{w(i)}=w(\,\bar{i}\,)$ for all $i\in \N$ and 
$w(i)=i$ for $i>n$ for some $n\geq 1.$
Denote by $\overline{S}_\infty$ the group of all signed permutations.
We often denote $w\in \overline{S}_\infty$ by one-line form $w(1)w(2)\cdots w(n)\cdots.$ 
Note that we only have to specify $w(i)$ for positive $i.$ 
For example $w=\bar{5}2\bar{1}\bar{3}4\cdots$
is a signed permutation, where dots mean the part 
that $w(i)=i$ holds. 

Define the signed permutations $s_i\;(i\geq 0)$ by
\begin{eqnarray}
s_0(1)=\overline{1},\quad
s_0(i)=i\;(i\in \N, \;i\neq 1),
\\
s_i(i)=i+1,\quad s_i(i+1)=i, \quad s_i(j)=j\;(j\in \N, \;j\neq i).
\end{eqnarray}
Then (\ref{WA1}), (\ref{WA2}), and (\ref{WA3}) hold,
and $W(C_\infty)=W(B_\infty)$ can be identified with
 $\overline{S}_\infty.$   
Let $s_{\hat{1}}=s_0s_1s_0.$ Explicitly, we have 
\begin{equation}
s_{\hat{1}}(1)=\overline{2},\quad
s_{\hat{1}}(2)=\overline{1},\quad
s_{\hat{1}}(i)=i\;(i\in \N, \;i\neq 1,2).
\end{equation}
Then $W(D_\infty)$ can be identified with the subgroup
of $\overline{S}_\infty$ %the group of all signed permutations
generated by $s_{\hat{1}},\;s_i\;(i\geq 1),$
which we denote this subgroup by $\overline{S}_{\infty,+}.$

\subsection{Strict partitions and Grassmannian elements} \label{ssec:SPGr}
Let $\mathcal{SP}$ denote the set of all strict partitions,
i.e. $ \mathcal{SP}=\bigcup_{n\geq 0}\mathcal{SP}_n.$
For each $w\in \overline{S}_\infty, $ 
define $\lambda(w)=\overline{\N}\cap w (\N).$ 
For example if $w=\bar{5}2\bar{1}\bar{3}4\cdots$ then $\lambda(w)=\{\bar{5},\bar{3},\bar{1}\}.$
Such $\lambda(w)$ can be considered as a strict partition;  for example the last $\lambda(w)$
corresponds to $(5,3,1)\in \mathcal{SP}.$
Thus we have a surjective map $\lambda: \overline{S}_\infty\rightarrow
\mathcal{SP}, w\mapsto \lambda(w).$ 
This map gives a bijection 
$\overline{S}_\infty/S_\infty \cong \mathcal{SP},$ where $S_\infty$ is the 
subgroup of $\overline{S}_\infty$ consisting of `ordinary' permutations, i.e. those 
$w\in \overline{S}_\infty$ such that $w(\N)\subset \N.$

An element $w\in \overline{S}_\infty$ is called {\it Grassmannian}, if 
$w(1)<w(2)<\cdots$, where the elements are ordered
as $\cdots <\bar{2}<\bar{1}<1<2<\cdots.$ 
For example $w=\bar{5}\bar{3}\bar{1}24\cdots$ is
Grassmannian. Let $\overline{S}_\infty^{\,0}$ denote the set of Grassmannian elements in $\overline{S}_\infty.$
Note that $\overline{S}_\infty^{\,0}$ is equal to the set of elements 
$w\in \overline{S}_\infty$ such that 
$\ell(ws_i)>\ell(w)$ for all $i$ other than $i=0,$
where $\ell$ is the length function of $W(B_\infty)=W(C_\infty)\cong \overline{S}_\infty.$

The next fact is well-known.  
See for example \cite{BL}, \cite{EYD} (note that our convention here is different from there). 

\begin{prop}\label{prop:SP-BC}
The set  $\overline{S}_\infty^{\,0}$ forms a set of  
coset representatives for $\overline{S}_\infty/S_\infty.$
Let $\lambda$ denote the  %order preserving 
resulting bijection
$\overline{S}_\infty^{\,0}\cong \mathcal{SP}.$
Then for $w,v\in  \overline{S}_\infty^{\,0}$, we have $|\lambda(w)|=\ell(w),$ and 
$\lambda(w)\subset \lambda(v)\Longleftrightarrow
w\leq v$ (Bruhat-Chevalley order).
\end{prop}

Next consider the type $D$ case, i.e. $\overline{S}_{\infty,+}.$
The group $\overline{S}_{\infty,+}$ consists of the  
signed permutations such that the cardinality of $\lambda(w)$ is even.
The image $\lambda(\overline{S}_{\infty,+})\subset \mathcal{SP}$
is a subset of $\lambda$ having even length.
This image 
can be also identified with $\mathcal{SP}$ by 
``removing all the diagonal boxes'' as illustrated by the following:

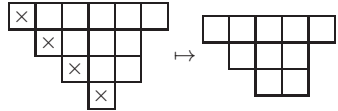
\begin{center}
\begin{tiny}
$\Tableau{{\times} & {}& {}&{}&{}&{}\\ ~&{\times}&{}&{}&{}\\~&~&{\times}&{}&{}\\~&~&~&{\times}}\mapsto
\Tableau{{} & {}& {}& {}& {}\\ ~&{}&{}& {}\\~&~&{}& {}}$
\end{tiny}
\end{center}
Let $\lambda_+: \overline{S}_{\infty,+}\rightarrow \mathcal{SP}$ denote the obtained map.
For example if $w=\bar{6}\bar{4}\bar{3}\bar{1}25\cdots$
then $\lambda_+(w)=(5,3,2)\in \mathcal{SP}.$
Thus $\lambda_+$ gives a bijection $\overline{S}_{\infty,+}/S_\infty\cong \mathcal{SP}.$
The set of Grassmannian elements in $\overline{S}_{\infty,+}$ is denoted by 
$\overline{S}_{\infty,+}^{\,\hat{1}}.$
As in the previous case of type $B,C$, the set $\overline{S}_{\infty,+}^{\,\hat{1}}$ 
is equal to the set of elements 
$w\in \overline{S}_{\infty,+}$ such that 
$\ell(ws_i)>\ell(w)$ for all $i\in I-\{\hat{1}\}.$

The corresponding result (also well-known) for type $D$ is the following.
\begin{prop}\label{prop:SP-D}
The set $\overline{S}_{\infty,+}^{\,\hat{1}}$ forms a set of %standard 
coset representatives for $\overline{S}_{\infty,+}/S_\infty.$ 
Let $\lambda_+$ denote the resulting bijection
$\overline{S}_{\infty,+}^{\,\hat{1}}\cong \mathcal{SP}.$
Then for $w,v\in  \overline{S}_{\infty,+}^{\,\hat{1}}$, we have $|\lambda_+(w)|=\ell(w),$ and 
$\lambda(w)\subset \lambda(v)\Longleftrightarrow
w\leq v$ (Bruhat-Chevalley order).
\end{prop}

Note that the length function $\ell$ and the Bruhat order are
those of type $D_\infty.$

\bigskip

Through the bijections given in 
Prop. \ref{prop:SP-BC} and Prop. \ref{prop:SP-D}, 
Weyl group acts naturally on $\mathcal{SP}.$
We will describe the action explicitly. 
Define $c(\alpha)\in I$ of $\alpha \in \mathbb{D}(\lambda)$ to be 
$j-i$ for type $B$ and $C$ cases. For type $D$ case, we define 
$c(i,j)=j-i+1\;(i<j)$ and $c(i,i)=\hat{1}$ of $i$ is even 
and $c(i,i)=1$ of $i$ is odd. For each type, a strict partition 
$\lambda$ is $i$-{\it removable} 
if there is a box $\alpha\in \mathbb{D}(\lambda)$ with $c(\alpha)=i$ such that 
$\mathbb{D}(\lambda)-\{\alpha\}$ is a shifted diagram.
Then we denote the diagram $\mathbb{D}(\lambda)-\{\alpha\}$ by $\lambda^{(i)}.$ 
Conversely $\lambda$ is $i$-{\it addable\/} 
if there is $\mu$ such that $\mu$ is $i$-removable 
and $\mu^{(i)}=\lambda.$
Then we denote $\mu=\lambda^{+(i)}.$
%Let $\lambda\in \mathcal{SP}$ and 
%$w$ be the corresponding element in $\overline{S}^0_\infty.$
%$s_iw$
%For each type, there is an action of Weyl group on 
% $\mathcal{SP}$ via the corresponding bijection of Prop. \ref{prop:SP-BC} or Prop. \ref{prop:SP-D}.
 Note that the bijection gives an ordering on $\mathcal{SP}$
 induced by the Bruhat-Chevalley order of the corresponding type.
 We denote $\la<\mu$ when $\la\leq \mu$ and $\la\neq \mu.$
 \begin{prop}\label{prop:s_ila} Let $\lambda\in \mathcal{SP}.$
Then 

\begin{itemize}
\item $s_i\lambda <\lambda\Longleftrightarrow $ 
$\lambda$ is $i$-removable, 

\item $s_i\lambda >\lambda\Longleftrightarrow$
$\lambda$ is $i$-addable,

\item $s_i\lambda=\lambda
\Longleftrightarrow$
$\lambda$ is neither $i$-removable nor
$i$-addable. 

\end{itemize}
Furthermore, if $s_i\lambda <\lambda$, then $s_i \lambda=\lambda^{(i)},$ %,
and if $s_i\lambda >\lambda$, then $s_i \lambda=\lambda^{+(i)}.$
\end{prop}
{\it Proof.} We first consider the case $\overline{S}_\infty.$
For each $j\geq 1$, define $\rho_j=s_{j-1}\cdots s_1s_0.$
It is easy to check the following relations:
\begin{equation}
s_i\rho_j=\rho_js_{i+1}\;(1\leq i\leq j-2),\quad
s_i\rho_j=\rho_js_i\;(i\geq j+1),
\label{eq:s_i_comm}
\end{equation}
\begin{equation}
s_{j-1}\rho_j=\rho_{j-1}\;(\mbox{with}\;\rho_0=1),\quad
s_j\rho_{j}=\rho_{j+1}.\label{eq:srho}
\end{equation}
If $\la=(\la_1,\ldots,\la_l)\in \mathcal{SP}$, then the corresponding 
element in $\overline{S}_\infty^0$ is given by 
$w_\la=\rho_{\la_l}\cdots \rho_{\la_1}$
(see \cite{EYD}, \S 7).
We first consider the case $i\geq 1.$
If $\la\in \mathcal{SP}$ is $i$-removable, then there is $k\geq 1$ such that $\la_k=i+1$ and $\la_{k+1}\leq i-1.$
By using the relations (\ref{eq:s_i_comm}), (\ref{eq:srho}), we have $s_iw_\la=w_{\la^{(i)}}.$ 
If 
$\la\in \mathcal{SP}$ is 
$i$-addable, then there is $k\geq 1$ such that 
 $\la_{k+1}=i$ and $\la_{k}\geq i+2.$
Then we have $s_iw_\la=w_{\la^{+(i)}}$ 
also by using the relations (\ref{eq:s_i_comm}), (\ref{eq:srho}).
Next we consider the case that 
$\lambda$ is neither $i$-removable nor
$i$-addable. 
We need one more important relation:
\begin{equation}
s_i\rho_i\rho_{i+1}=\rho_i\rho_{i+1}s_1,\label{i,i,i+1}
\end{equation}
which can be checked straightforwardly.
If there is $k\geq 1$ such that $\la_k=i+1$ and $\la_{k+1}=i$, then 
$
w_\la=\rho_{\la_l}\cdots \rho_{i}\rho_{i+1}\cdots \rho_{\la_1},
$
so by using (\ref{eq:s_i_comm}), (\ref{i,i,i+1}), we have 
$s_iw_\la=w_\la s_{k}. $
%
%$\la_k=i$ and $\la_{k+1}=i+1$ then 
%$s_iw_\la=w_\la s_{k+1}$
The remaining case is that all $\la_j$ are different from $i$ and $i+1$. Then there is $k$ such that 
 $\la_k\geq i+2$ and $\la_{k+1}\leq i-1.$ 
Then we have $s_iw_\la=w_\la s_{i+k}$
by using (\ref{eq:s_i_comm}).
Next consider the case $i=0$. The strict partition 
$\la$ is $0$-removable if and only if $\la_l=1,$
otherwise $\la$ is $0$-addable. Then the result is 
obvious from the form of the reduced expression 
$w_\la.$

Now we will consider type $D$.
We define $\rho^{-}_j=s_j\cdots s_2s_{{1}}$ for $j\geq 1$, and 
$$
\rho^{+}_j=s_j\cdots s_2s_{\hat{1}}\quad(j\geq 2),\quad
\rho^+_1=s_{\hat{1}}.
$$
Let $\la\in\mathcal{SP}.$
The corresponding element in
$\overline{S}^{\hat{1}}_{\infty,+}$ 
is given by 
$
w_\la=\rho_{\la_l}^{\varepsilon}
\cdots
\rho_{\la_2}^{-}\rho_{\la_1}^{+}$, where $\varepsilon=+$, if $l$ is odd
and $\varepsilon=-$, if $l$ is even.
We have the following easily checked relations:
$$
s_i\rho_j^{\pm}=\rho_j^{\pm}s_{i+1}\;(2\leq i\leq j-1),\quad
s_i\rho^{\pm}_j=\rho^{\pm}_js_{i}\;(i\geq j+2),
$$
$$
s_{j+1}\rho_j^\pm=\rho_{j+1}^\pm,\quad
s_j\rho_j^\pm=\rho_{j-1}^\pm,\quad
s_{{1}}\rho^-_j=\rho_j^-s_2,\quad
s_{\hat{1}}\rho^+_j=\rho_j^+s_2\;(j\geq 2).$$
The results for $s_i\;(i\geq 2)$ are proved as in the proof for $s_i\;(i\geq 1)$ of 
type $C$ case.
Let us consider $i=1$ (the case $i=\hat{1}$ is similar). 
If $\la$ is $1$-removable or 
$1$-addable, the desired relation is obvious. 
Suppose $\la$ is neither $1$-removable nor 
$1$-addable.
This occurs when (i) $l$ is odd and $\la_{l}=1$,
or (ii) $l$ is even and $\la_{l}\geq 2.$
If (i) holds, then we have $s_1w_\la=w_\la s_l.$
If (ii) holds, then we have $s_1w_\la=w_\la s_{l+1}.$
$\qed$

\begin{example}

If $\lambda=(3,1)$ then for type $B,C$, we have the following examples:

\setlength{\unitlength}{0.4mm}
  \begin{picture}(200,35)
  \put(0,18){\tiny{$\lambda=$}}
  \put(15,25){\line(1,0){30}}
  \put(15,15){\line(1,0){30}}
  \put(25,5){\line(1,0){10}}
  \put(15,15){\line(0,1){10}}
  \put(25,5){\line(0,1){20}}
  \put(35,5){\line(0,1){20}}
 \put(45,15){\line(0,1){10}}
  \put(27,8){\tiny{$0$}}
  \put(17,18){\tiny{$0$}}
  \put(27,18){\tiny{$1$}}
    \put(37,18){\tiny{$2$}}
  \put(70,18){\tiny{$s_0\lambda=$}}
  \put(95,25){\line(1,0){30}}
  \put(95,15){\line(1,0){30}}
%  \put(105,25){\line(1,0){20}}
  \put(95,15){\line(0,1){10}}
  \put(105,15){\line(0,1){10}}
  \put(115,15){\line(0,1){10}}
  \put(125,15){\line(0,1){10}}
%  \put(135,35){\line(0,1){10}}
%  \put(107,28){\tiny{$0$}}
  \put(97,18){\tiny{$0$}}
  \put(107,18){\tiny{$1$}}
%  \put(117,28){\tiny{$1$}}
  \put(117,18){\tiny{$2$}}
%  \put(127,38){\tiny{$3$}}
  %
  \put(140,18){\tiny{$s_1\lambda=$}}
  \put(165,25){\line(1,0){30}}
  \put(165,15){\line(1,0){30}}
  \put(175,5){\line(1,0){20}}
  \put(165,15){\line(0,1){10}}
  \put(175,5){\line(0,1){20}}
  \put(185,5){\line(0,1){20}}
  \put(195,5){\line(0,1){20}}
%  \put(205,35){\line(0,1){10}}
  \put(177,8){\tiny{$0$}}
  \put(167,18){\tiny{$0$}}
  \put(177,18){\tiny{$1$}}
  \put(187,8){\tiny{$1$}}
  \put(187,18){\tiny{$2$}}
%  \put(197,38){\tiny{$3$}}
  %
  \put(210,18){\tiny{$s_2\lambda=$}}
  \put(235,25){\line(1,0){20}}
  \put(235,15){\line(1,0){20}}
  \put(245,5){\line(1,0){10}}
  \put(235,15){\line(0,1){10}}
  \put(245,5){\line(0,1){20}}
  \put(255,5){\line(0,1){20}}
%  \put(265,15){\line(0,1){10}}
%  \put(205,35){\line(0,1){10}}
  \put(247,8){\tiny{$0$}}
  \put(237,18){\tiny{$0$}}
  \put(247,18){\tiny{$1$}}
%  \put(257,28){\tiny{$1$}}
%  \put(257,18){\tiny{$2$}}
%  \put(197,38){\tiny{$3$}}
%
  \put(280,18){\tiny{$s_3\lambda=$}}
  \put(305,25){\line(1,0){40}}
  \put(305,15){\line(1,0){40}}
  \put(315,5){\line(1,0){10}}
  \put(305,15){\line(0,1){10}}
  \put(315,5){\line(0,1){20}}
  \put(325,5){\line(0,1){20}}
  \put(335,15){\line(0,1){10}}
  \put(345,15){\line(0,1){10}}
  \put(317,8){\tiny{$0$}}
  \put(307,18){\tiny{$0$}}
  \put(317,18){\tiny{$1$}}
%  \put(257,28){\tiny{$1$}}
  \put(327,18){\tiny{$2$}}
  \put(337,18){\tiny{$3$}}
   \end{picture}
   
\noindent where we filled $c(\alpha)$ in each box $\alpha$.
For the same $\lambda$ in type $D$, we have

\setlength{\unitlength}{0.4mm}
  \begin{picture}(200,35)
  \put(0,18){\tiny{$\lambda=$}}
  \put(15,25){\line(1,0){30}}
  \put(15,15){\line(1,0){30}}
  \put(25,5){\line(1,0){10}}
  \put(15,15){\line(0,1){10}}
  \put(25,5){\line(0,1){20}}
  \put(35,5){\line(0,1){20}}
 \put(45,15){\line(0,1){10}}
  \put(27,8){\tiny{$1$}}
  \put(17,18){\tiny{$\hat{1}$}}
  \put(27,18){\tiny{$2$}}
    \put(37,18){\tiny{$3$}}
  \put(70,18){\tiny{$s_{\hat{1}}\lambda=$}}
  \put(95,25){\line(1,0){30}}
  \put(95,15){\line(1,0){30}}
  \put(105,5){\line(1,0){10}}
  \put(95,15){\line(0,1){10}}
  \put(105,5){\line(0,1){20}}
  \put(115,5){\line(0,1){20}}
  \put(125,15){\line(0,1){10}}
%  \put(135,35){\line(0,1){10}}
%  \put(107,28){\tiny{$0$}}
  \put(97,18){\tiny{$\hat{1}$}}
  \put(107,18){\tiny{$2$}}
 \put(107,8){\tiny{$1$}}
  \put(117,18){\tiny{$3$}}
%  \put(127,38){\tiny{$3$}}
  %
  \put(140,18){\tiny{$s_1\lambda=$}}
  \put(165,25){\line(1,0){30}}
  \put(165,15){\line(1,0){30}}
%  \put(175,25){\line(1,0){20}}
  \put(165,15){\line(0,1){10}}
  \put(175,15){\line(0,1){10}}
  \put(185,15){\line(0,1){10}}
  \put(195,15){\line(0,1){10}}
%  \put(205,35){\line(0,1){10}}
%  \put(177,28){\tiny{${1}$}}
  \put(167,18){\tiny{$\hat{1}$}}
  \put(177,18){\tiny{$2$}}
%  \put(187,28){\tiny{$3$}}
  \put(187,18){\tiny{$3$}}
%  \put(197,38){\tiny{$3$}}
  %
  \put(210,18){\tiny{$s_2\lambda=$}}
  \put(235,25){\line(1,0){30}}
  \put(235,15){\line(1,0){30}}
  \put(245,5){\line(1,0){20}}
  \put(235,15){\line(0,1){10}}
  \put(245,5){\line(0,1){20}}
  \put(255,5){\line(0,1){20}}
  \put(265,5){\line(0,1){20}}
%  \put(205,35){\line(0,1){10}}
  \put(247,8){\tiny{$1$}}
  \put(237,18){\tiny{$\hat{1}$}}
  \put(247,18){\tiny{$2$}}
  \put(257,8){\tiny{$2$}}
  \put(257,18){\tiny{$3$}}
%  \put(197,38){\tiny{$3$}}
%
  \put(280,18){\tiny{$s_3\lambda=$}}
  \put(305,25){\line(1,0){20}}
  \put(305,15){\line(1,0){20}}
  \put(315,5){\line(1,0){10}}
  \put(305,15){\line(0,1){10}}
  \put(315,5){\line(0,1){20}}
  \put(325,5){\line(0,1){20}}
%  \put(335,15){\line(0,1){10}}
%  \put(345,35){\line(0,1){10}}
  \put(317,8){\tiny{$1$}}
  \put(307,18){\tiny{$\hat{1}$}}
  \put(317,18){\tiny{$2$}}
%  \put(257,28){\tiny{$1$}}
%  \put(327,18){\tiny{$3$}}
%  \put(337,38){\tiny{$3$}}
   \end{picture}
\noindent etc.
%  \end{center}
\end{example}

{{}
\subsection{Grassmannian elements in the finite rank case}
Let  $\mathcal{SP}(n)$ denote the set of strict partitions 
such that $\la\subset \rho_n.$
Consider a subgroup of $W(X_n)$ 
generated by 
$s_1,\ldots,s_{n-1},$ which is
isomorphic to $S_n.$
If we denote it simply by $S_n$, we have the following natural bijection:
$$
\mathcal{SP}(n)\cong
W(C_n)/S_n=
W(B_n)/S_n
\cong
W(D_{n+1})/S_{n+1}.
$$
Indeed, when we identify $W(C_\infty)$ with $\overline{S}_\infty$,
a set of coset representatives
for $W(C_n)/S_n$
is given by $W(C_{n})\cap \overline{S}_\infty^0$, which is 
naturally identified with $\mathcal{SP}(n).$
Similarly, $W(D_{n+1})/S_{n+1}$ can be identified 
with $W(D_{n+1})\cap \overline{S}_{\infty,+}^{\hat{1}}
\cong \mathcal{SP}(n).$
}

\subsection{The ring $\R$ and $W$-action}\label{subsec:WR}

So far, we have worked over the ring of coefficients $\Z[\beta][b_1,b_2,\ldots].$ 
Now we need to work over a larger coefficient ring:
\begin{equation}
\R:=\Z[\beta][b_1,b_2,\ldots,{b}_{\bar{1}},{b}_{\bar{2}},\ldots]
/\langle b_i\oplus {b}_{\bar{i}}\;|\;i\geq 1\rangle
\end{equation}
%As is shown below 
This is essentially the ring of Laurent polynomials
in infinite variables with coefficients in $\Z[\beta].$
One sees that $\R$ is isomorphic to $\Z[\beta][b_1,b_2,\ldots][(1+\beta b_i)^{-1}(i\geq 1)]
$ (see \S \ref{subsec:Gro}) and also to $\Z[\beta][{b}_{\bar{1}},{b}_{\bar{2}},\ldots][(1+\beta {b}_{\bar{i}})^{-1}(i\geq 1)].$
$\R$ has a natural action of the Weyl group defined 
by $w(b_i)=b_{w(i)},\;w(b_{\bar{i}})=b_{\overline{w(i)}}\;(w\in W,\;i=1,2,\ldots)$ 
with $\beta$ fixed.

\subsection{Roots}\label{ssec:roots}
We fix notation about root systems.
Let $L$ denote the free $\Z$-module with basis
$\{t_i\}_{i\geq 1}.$
The {\it positive roots\/} $\Delta^{+}\subset L$ 
(set $\Delta^{-}:=-\Delta^{+}$ the {\it negative roots}) 
are defined by 
\begin{eqnarray*}
\mbox{Type}\; {B}_{\infty}:&
\Delta^{+}&=\{t_{i}\;|\; i\geq 1\}
\cup\{t_{j}\pm t_{i}\;|\;j>i\geq 1\},\\
\mbox{Type}\; {C}_{\infty }:&
\Delta^{+}&=\{2t_{i}\;|\; i\geq 1\}
\cup\{t_{j}\pm t_{i}\;|\;j>i\geq 1\},\\
\mbox{Type}\; {D}_{\infty}:&
\Delta^{+}&=\{t_{j}\pm t_{i}\;|\;j>i\geq 1\}.
\end{eqnarray*}
For a positive integer $n$,
define $\Delta_n^+:=\Delta^+\cap L_n$ where 
$L_n:=\bigoplus_{i=1}^n\Z t_i.$
These are the corresponding positive root system of finite rank $n$
of type 
$B_n,C_n$, and $D_n.$

The following elements of $\Delta^+$ are called the {\it simple roots}:
\begin{eqnarray*}
\mbox{Type}\; {B}_{\infty}:&
\alpha_{0}&=t_{1},
\quad \alpha_{i}=t_{i+1}-t_{i}\quad(i\geq 1) ,\\
\mbox{Type}\; {C}_{\infty}:&
\alpha_{0}&=2 t_{1},
\quad
\alpha_{i}=t_{i+1}-t_{i}\quad(i\geq 1),\\
\mbox{Type}\; {D}_{\infty}:&
\alpha_{\hat{1}}&=t_{1}+t_{2},
\quad \alpha_{i}=t_{i+1}-t_{i}\quad(i\geq 1).
\end{eqnarray*}

We can define a map 
$e: L=\bigoplus_{i=1}^\infty \Z t_i\rightarrow \R$
satisfying 
$
e(t_i)= b_i,\;
e(-t_i)={b}_{\bar{i}},
$
$
e(\alpha+\gamma)=e(\alpha)\oplus e(\gamma),\;
e(\alpha-\gamma)=e(\alpha)\ominus e(\gamma)\;
(\alpha,\,\gamma\in L).
$
Note that 
the map $e$ is compatible with the natural 
action of $W$ on $L$ 
and the action on
$\R$ defined in \S \ref{subsec:WR}.
For simple roots $\alpha_i$, the explicit form of $e(\alpha_i)$ 
written in terms of $b_i$'s are given by
\begin{eqnarray*}
B_\infty&: &e(\alpha_0)=b_1,\quad
e(\alpha_i)=b_{i+1}\ominus b_i\;(i\geq 1),\\
C_\infty&: &e(\alpha_0)=b_1\oplus b_1,\quad
e(\alpha_i)=b_{i+1}\ominus b_i\;(i\geq 1),\\
D_\infty&: &e(\alpha_{\hat{1}})=b_1\oplus b_2,\quad
e(\alpha_i)=b_{i+1}\ominus b_i\;(i\geq 1).
\end{eqnarray*}

\setcounter{equation}{0}
\section{GKM ring and its Schubert basis}\label{sec:GKM}
In this section, we introduce a ring $\Psi$, which is 
defined by the $K$-theoretic GKM condition.
A notion of ``Schubert classes'' is defined in a combinatorial way.

\subsection{GKM ring $\Psi$}\label{ssec:GKM}
Let $\mathrm{Fun}(\mathcal{SP},\R)$ denote the set of 
all maps from $\mathcal{SP}$ to $\R.$
%$\mathrm{Fun}(\mathcal{SP},\R).$
$\mathrm{Fun}(\mathcal{SP},\R)$ is naturally an $\R$-algebra;  
the $\R$-module structure is given by diagonal multiplication,
and multiplication is defined in point wise manner. 
For each $\alpha\in \Delta^+$, we have $\alpha=w(\alpha_i)$ for 
some $i\in I$ and $w\in W.$ Then let $s_\alpha=w s_i w^{-1}.$
\begin{Def}[GKM ring] Let $\Psi$ be the subring of  
$\mathrm{Fun}(\mathcal{SP},\R)$ defined by the following condition:
$$
\psi(s_\alpha \mu)-\psi(\mu)\in e(\alpha)\cdot \R\quad
\mbox{for all}\;
\mu\in \mathcal{SP}\;
\mbox{and}\;
\alpha\in \Delta^+.
$$
\end{Def}

This condition and the associated geometry 
is discussed in \S \ref{ssec:K_T}.

\subsection{Divided difference operators on $\Psi$}

We %can 
now define the {\it divided difference operators\/} $\pi_i\;(i\in I)$ on 
$\Psi$ by the following formula:
\begin{equation}
(\pi_i \psi)(\mu)=\frac{s_i (\psi(s_i\mu))-(1+\beta e(\alpha_i))\psi(\mu)}{e(\alpha_i)}
\quad (\psi\in \Psi).
\end{equation}
This can be rewritten as follows:
\begin{equation}
(\pi_i \psi)(\mu)=\frac{\psi(\mu)-(1+\beta e(-\alpha_i))s_i (\psi(s_i\mu))}{e(-\alpha_i)}.
\end{equation}
If $\psi\in \Psi$ then 
the family $\{(\pi_i \psi)(\mu)\}_\mu$
actually gives an element of $\Psi.$
This fact can be shown by a similar argument as in \cite{KT} (the first Lemma in the Appendix).

The operators $\pi_i$ satisfy the following relation:
if $W$ is of type $C_\infty$ ($B_\infty$) then 
$\pi_i^2=-\beta\pi_i\;(i=0,1,2\ldots)$ and 
\begin{eqnarray}
\pi_i\pi_{i+1}\pi_i&=&\pi_{i+1}\pi_i\pi_{i+1}\quad (i\geq 1),\label{piWA1}\\
\pi_i\pi_j&=&\pi_j\pi_i\quad (i,j\geq 1,\;|i-j|\geq 2),\label{piWA2}\\
\pi_0\pi_1\pi_0\pi_1&=&\pi_1\pi_0\pi_1\pi_0,\quad 
\pi_0\pi_i=\pi_i\pi_0\quad (i\geq 1).
\end{eqnarray}
If $W$ is type $D_\infty$ then   
$\pi_i^2=-\beta\pi_i\;(i=\hat{1},1,2\ldots)$, (\ref{piWA1}), (\ref{piWA2}), and 
\begin{eqnarray}
\pi_{\hat{1}}\pi_2\pi_{\hat{1}}&=&\pi_2\pi_{\hat{1}}\pi_2,\quad
\pi_{\hat{1}}\pi_j=\pi_j \pi_{\hat{1}}\quad (j\neq 2).
\end{eqnarray}

\begin{rmk}
The operators $\pi_i$ are `left' divided difference operators.
These operators do not appear explicitly in \cite{KK},
where `right' divided difference operators are used.
Unfortunately, we cannot
find a reference that clarify a geometric origin of the `left' one.
However, it turns out that the left divided difference operators
are very useful especially for `parabolic' situations. 
See \cite{Kn, KT} for some discussions 
and applications of the left divided difference operators
in equivariant cohomology.
\end{rmk}

\subsection{Schubert classes}

\begin{Def} \label{def:Schubert}
A set of elements
$\{\psi_\lambda|\;\lambda\in \mathcal{SP}\}\subset \Psi$ is 
called a family of {\it Schubert
classes} if the following conditions are satisfied:
\begin{equation}
\pi_i \psi_\lambda=\begin{cases}\psi_{\lambda^{(i)}}
&\mbox{if}\quad s_i\lambda<\lambda\\
-\beta \psi_{\lambda}& \mbox{if}\quad s_i\lambda\geq \lambda
\end{cases},\label{eq:divdiffSchubert}
\end{equation}
\begin{equation}
\psi_\lambda(\emptyset )=\delta_{\lambda,\emptyset}\mbox{ (Kronecker's delta)}.\label{eq:InitialSchubert}
\end{equation}
\end{Def}

A family of Schubert classes exists.
This fact can be proved by a geometric argument in \S \ref{sec:geometry}.
We have a proof of the existence 
as a consequence of Thm. \ref{thm:main}.

\begin{lem}\label{lem:rec}
Let $\{\psi_\lambda\}_\lambda$ be a family of Schubert classes.
Suppose $\mu$ be a strict partition such that $\mu\ne \emptyset.$ 
Let $i\in I$ be such that $s_i \mu<\mu.$ Then
\begin{equation}
\psi_\lambda(\mu)
=\begin{cases}
(1+\beta e(-\alpha_i))s_i\left( \psi_\lambda(s_i\mu)\right)+e(-\alpha_i)s_i\left( 
\psi_{s_i \lambda}(s_i\mu)\right)
& \mbox{if}\quad s_i \lambda<\lambda\\
s_i \left(\psi_\lambda(s_i\mu)\right)&  \mbox{if}\quad s_i \lambda\geq \lambda\\
\end{cases}.\label{eq:rec}
\end{equation}
\end{lem}
{\it Proof.}
This is a direct consequence of 
the divided difference equation.
$\qed$

\begin{rmk}
Recurrence equation (\ref{eq:rec}) is called `left-hand' recurrence in 
\cite{LSS}, Remark 2.3. 
\end{rmk}

\begin{prop}\label{prop:unique}
A family of Schubert classes, if exists, is unique.
\end{prop}
{\it Proof.} Let $\{\psi_\lambda\}_{\lambda\in \mathcal{SP}}$ be 
a family of Schubert classes. 
By definition, we have $\psi_\la(\emptyset)=\delta_{\la,\emptyset}.$  
Suppose $\mu\neq\emptyset.$ 
Then there is $i\in I$ such that $s_i\mu<\mu.$ 
By the recurrence equation (\ref{eq:rec}),
$\{\psi_\lambda(\mu)\;|\;\lambda\in \mathcal{SP}\}$
is uniquely determined from 
$\{\psi_\lambda(s_i\mu)\;|\;\lambda\in \mathcal{SP}\}.$
Hence by using induction on $|\mu|$, we conclude the uniqueness. 
$\qed$

\begin{prop}\label{prop:vanish}
Let  $\{\psi_\lambda\}_\lambda$ be a family of Schubert classes. 

{\rm (1)} If $\lambda\not\subset \mu$ then $\psi_\lambda(\mu)=0.$

{\rm (2)} We have $
\psi_\lambda(\lambda)=\prod_{\alpha\in \mathrm{Inv}(\lambda)}e(-\alpha),\label{eq:euler}$
where $\mathrm{Inv}(\lambda)=\{\alpha\in \Delta^{+}\;|\; s_\alpha \lambda<\lambda\}.$
\end{prop}

{\it Proof.} (1) We proceed by induction on $|\mu|.$
For $\mu=\emptyset$, the vanishing property holds by the initial condition.
Let $\mu\neq \emptyset .$
There exists $i$ such that $s_i \mu<\mu.$ 
Then we have $\lambda\not\subset s_i \mu.$ 
So by inductive hypothesis, 
we have $\psi_{\lambda}(s_i\mu)=0.$
Thus, if $s_i\lambda\geq\lambda$, then
from (\ref{eq:rec}), we have $\psi_\lambda(\mu)=0.$
Next suppose $s_i\lambda<\lambda.$ 
Then since both $\lambda$ and $\mu$ are $i$-removable,
we have $s_i\lambda\not\subset s_i\mu.$
By inductive hypothesis, we have 
$\psi_{\lambda}(s_i\mu)=\psi_{s_i \lambda}(s_i \mu)=0.$
Hence from (\ref{eq:rec}), we have $\psi_\lambda(\mu)=0.$

(2) There exists $i$ such that 
$s_i\lambda<\lambda$. 
By the recurrence equation (\ref{eq:rec}) together 
with $\psi_\lambda(s_i\la)=0$ as a consequence of (1), we have
$
\psi_\lambda(\lambda)=e(-\alpha_i)\,s_i\left( \psi_{s_i\lambda}(s_i\lambda)\right).
$
By this equation
$\psi_\lambda(\lambda)$ are determined inductively.
We see that 
$\prod_{\alpha\in \mathrm{Inv}(\lambda)}e(-\alpha)$ satisfies 
this equation
since we have $\mathrm{Inv}(\lambda)=s_i\, \mathrm{Inv}(s_i \lambda)\cup \{\alpha_i\}.$
$\qed$

\begin{prop}\label{prop:Schubert_basis}
Let  $\{\psi_\lambda\}_\lambda$ be a family of Schubert classes. 
Any element in $\Psi$ can be uniquely expressed as a possibly 
infinite $\R$-linear combination of $\psi_\lambda
\;(\lambda\in \mathcal{SP}).$
\end{prop}
{\it Proof.} Let $\psi\in \Psi.$
Define $\mathrm{Supp}(\psi)=\{\mu\in\mathcal{SP}\;|\;\psi(\mu)\neq 0\}.$ 
Let $\lambda\in \mathrm{Supp}(\psi)$ be
a minimum element in the ordering of strict partitions by inclusion. From the GKM condition 
and Prop. \ref{prop:vanish}, we see that 
$\psi(\lambda)$ is divisible by $e(-\alpha)$ for 
all $\alpha\in \mathrm{Inv}(\lambda).$
We see that the elements $\{e(-\alpha)\;|\;\alpha\in \mathrm{Inv}(\lambda)\}$
are relatively prime, and hence $\psi(\lambda)$ is divisible by their product $\psi_\lambda(\lambda)
=\prod_{\alpha\in \mathrm{Inv}(\lambda)}e(-\alpha).$
Let $\psi'=\psi-\psi(\lambda)/\psi_\lambda(\lambda)\cdot \psi_\lambda.$
By the minimality of $\la$ in $\mathrm{Supp}(\psi)$ and 
Prop. \ref{prop:vanish} (1), 
$\mathrm{Supp}(\psi')\subsetneqq \mathrm{Supp}(\psi).$
By repeating this, we may write $\psi$ as 
a possibly infinite $\R$-linear combination of $\psi_\lambda$'s.
%, 
$\qed$

\setcounter{equation}{0}
\section{Divided difference equation for $GP_\lambda(x|b)$ and $GQ_\lambda(x|b)$}\label{sec:divdiff}
In this section, we prove a divided difference 
equation for 
$GP_\lambda(x|b)$ and $GQ_\lambda(x|b).$

\subsection{Rings $G\Gamma_\R$ and inverse limit of the 
$K$-theoretic factorial
$P$- and $Q$-functions}\label{ssec:GGamma_R}
 Let $X$ be $B,C$ or $D.$
In order to state results for different types
 in a parallel way,
we use the following convention:
we denote by $G\Gamma_\R^X$ the following : 
$G\Gamma_\R^C=\R\otimes_{\Z[\beta]}G\Gamma_{+},$
$G\Gamma_\R^B=G\Gamma_\R^D=\R\otimes_{\Z[\beta]}G\Gamma.$
We often suppress $X$ when there is no fear of 
confusion.

For each $\la\in \mathcal{SP}_n,$ we define 
\begin{eqnarray}
GB_\la^{(n)}(x|b)&=&GP_\la(x_1,\ldots,x_n|0,b_1,b_2,\ldots),\\
GC_\la^{(n)}(x|b)&=&GQ_\la(x_1,\ldots,x_n|b_1,b_2,\ldots).\\
GD_\la^{(n)}(x|b)&=&\begin{cases}GP_\la(x_1,\ldots,x_n|b_1,b_2,\ldots)
&  \mbox{if}\;n\;\mbox{is even}\\
GP_\la(x_1,\ldots,x_n,0|b_1,b_2,\ldots) & \mbox{if}\;n\;\mbox{is odd}.\label{eq:defGXfinite}
\end{cases}
\end{eqnarray}
We define $GX_\lambda(x|b)=\varprojlim GX_\la^{(n)}(x|b)\in G\Gamma_\R^X$
for $X=B,C,$ and $D.$
We also use notation $GQ_\la(x|b)=GC_\la(x|b)$
and $GD_\la(x|b)=GP_\la(x|b)$ (defined in \S \ref{ssec:ExpGP}).

\begin{prop} 
Any $f(x)\in G\Gamma^X_\R$  
can be expressed uniquely as 
a possibly infinite $\R$-linear combination of $GX_\lambda(x|b)$'s, 
\begin{equation}
f(x)=\sum_{\lambda\in \mathcal{SP}}c_\lambda\cdot GX_\lambda(x|b),\quad
c_\lambda\in \R,
\end{equation}
such that for all positive integer $n$, the set 
$\{\lambda\in \mathcal{SP}_n\;|\;c_\lambda\neq 0\}$
is finite.
\end{prop}
{\it Proof.}
Using Lem. \ref{G_basis} in place of Cor. \ref{cor:basis},
we can show this 
by the same argument as in 
\S \ref{ssec:PfBasis} (cf. Prop. \ref{prop:BasisGP} and Prop. \ref{prop:BasisGQ}).
$\qed$

\subsection{Action of $W$ on $G\Gamma_\R$}
We define an action of 
$W(X_\infty)$ on $G\Gamma_\R^X.$
Let $s_i\;(i\geq 1)$ act on $\R\otimes_{\Z[\beta]}G\Gamma$ 
by 
$$
s_i\left(\sum_{\alpha}c_\alpha\otimes \phi_\alpha\right)
=\sum_{\alpha}s_i(c_\alpha)\otimes \phi_\alpha\quad
(c_\alpha\in \R,\;
\phi_\alpha\in G\Gamma),
$$
where $s_i(c_\alpha)$ is defined in \S \ref{subsec:WR}.
The action of $s_0$ is defined as follows.
For 
$\phi\in G\Gamma,$ define 
\begin{equation}
(s_0\phi)(x_1,x_2,\ldots)=
\phi(b_1,x_1,x_2,\ldots),\label{eq:s_0}
\end{equation}
which is a well-defined element in $\R\otimes_{\Z[\beta]}G\Gamma.$
In general, define
\begin{equation}
s_0\left(\sum_\alpha c_\alpha \phi_\alpha\right)=\sum_\alpha s_0(c_\alpha)\cdot s_0(\phi_\alpha)\quad
(c_\alpha\in \R,\;\phi_\alpha\in G\Gamma).
\end{equation}
We have $s_0^2=\mathrm{id}$ on $\R\otimes_{\Z[\beta]}G\Gamma$ 
by virtue of the $K$-theoretic $Q$-cancellation property (cf. \cite{DSP}).
Define also $s_{\hat{1}}=s_0s_1s_0.$

\begin{prop} The operators $s_i\;(i\in I)$ give 
an action of 
$W(X_\infty)$ on $G\Gamma_\R^X.$
\end{prop}
{\it Proof.} The proof is straightforward (cf. \cite{DSP}, Prop. 7.2). $\qed$

\subsection{Divided difference operators $\pi_i$ on $G\Gamma_\R^X$} The {\it divided difference operator}
$\pi_i$ on $G\Gamma_\R^X$ is defined by  
\begin{equation}
\pi_i F=
\frac{s_i F-(1+\beta e(\alpha_i))F}{e(\alpha_i)}\quad 
\mbox{for all}\;\;
F\in G\Gamma_\R^X.\label{eq:def_pi}
\end{equation}
One can check 
directly that $s_i F-F$ is divisible by $e(\alpha_i)$ so 
the right-hand side is 
an element in $G\Gamma_\R^X.$

\begin{rmk} Using the equation 
$(1+\beta e(\alpha_i))(1+\beta e(-\alpha_i))=1$, 
we can rewrite (\ref{eq:def_pi}) as 
$$\pi_i F=
\frac{F-(1+\beta e(-\alpha_i))s_i F}{e(-\alpha_i)}\quad 
\mbox{for all}\;\;
F\in G\Gamma_\R^X.
$$
\end{rmk}

\begin{thm}\label{thm:piGQ} We have
$$\pi_i GX_\lambda(x|b)=\begin{cases}
GX_{\lambda^{(i)}}(x|b)&\mbox{if}\quad 
s_i \lambda<\lambda \\
-\beta \;GX_{\lambda}(x|b) & \mbox{if}\quad s_i\lambda\geq \lambda. 
\end{cases}$$
\end{thm}

{\it Proof.} We first prove the type $C$ case. %
Let $i\geq 1.$
We will work with $n$ variables $x_1,\ldots,x_n.$
Since the operator $\pi_i$ acts only non-trivially 
on $b_i$ and $b_{i+1}$, 
we only have to show
$$\pi_i [[x|b]]^{\lambda}
=\begin{cases}[[x|b]]^{\lambda^{(i)}}&\mbox{if}\quad s_i\lambda<\lambda\\
-\beta [[x|b]]^{\lambda}&\mbox{if}\quad s_i\lambda\geq \lambda.\end{cases}
$$
 If $s_i\lambda\geq \lambda$ then one sees that 
 $[[x|b]]^{\lambda}$ 
 is symmetric 
with respect to $b_i$ and $b_{i+1}.$ 
 It follows that $\pi_i GQ_\lambda(x|b)=-\beta GQ_\lambda(x|b).$
If $s_i\lambda<\lambda$ then 
$\lambda_k=i+1$ for some $k.$
Then we note that 
$\prod_{j< k}[[x|b]]^{\lambda_j}$ is symmetric 
with respect to $b_i$ and $b_{i+1},$
and 
$\prod_{j>k}[[x|b]]^{\lambda_j}$ does not 
depend on $b_i$ and $b_{i+1}.$
So we only have to calculate $[[x_k|b]]^{\lambda_k}.$
It is easy to check the equation
 $$
x_k\oplus b_{i+1}-(1+\beta b_{i+1}\ominus b_i)(x_k\oplus b_i)=b_{i+1}\ominus b_{i},
 $$
 which is equivalent to $\pi_i(x_k\oplus b_i)=1.$
So we have
$\pi_i \left([[x_k|b]]^{i+1}\right)=[[x_k|b]]^i.$ 
  
Next we consider the operator $\pi_0.$ 
Define 
$S_{n,r}=\{w\in S_n\;|\;w(r+1)<\cdots<w(n)\}.$ Then the definition of
$GQ_\lambda(x_1,\ldots,x_n,x_{n+1}|\ominus b_1,b_2,\ldots)$ reads
\begin{equation}
\sum_{w\in S_{n+1,r}}
\prod_{i=1}^{r}
(x_{w(i)}\oplus x_{w(i)})
(x_{w(i)}\ominus b_1)
(x_{w(i)}\oplus b_2)\cdots
(x_{w(i)}\oplus b_{\lambda_i-1})\times 
\prod_{i=1}^r
\prod_{j=i+1}^{n+1}\frac{x_{w(i)}\oplus x_{w(j)}}{x_{w(i)}\ominus x_{w(j)}}.\label{eq:pi0-1}
\end{equation}

Suppose $s_0\lambda<\lambda$ then we have 
$\lambda_j\geq 2$ for $1\leq j\leq r-1,\lambda_r=1.$ 
We will prove 
\begin{eqnarray}
&&\frac{
GQ_\lambda(x_1,\ldots,x_n,b_1|\ominus b_1,b_2,\ldots)
-(1+\beta (b_1\oplus b_1))GQ_\lambda(x_1,\ldots,x_n|b_1,b_2,\ldots)}{b_1\oplus b_1}
\nonumber \\
&=&GQ_{\lambda^{(0)}}(x_1,\ldots,x_n|b_1,b_2,\ldots).\label{eq:pi0}
\end{eqnarray}
If $w(k)=n+1$ for some $k$ such that $1\leq k\leq r-1$ then the 
corresponding term in (\ref{eq:pi0-1}) vanishes when we substitute $x_{n+1}=b_1.$ 
Define $S_{n+1,r}^I$ to be the set of elements $w\in S_{n+1,r}$ such that 
 $w(r)=n+1.$ Then the corresponding part of (\ref{eq:pi0-1}), 
after substituting 
$x_{n+1}=b_1,$  
becomes
\begin{equation}
(b_1\oplus b_1)\sum_{w\in S_{n+1,r}^I}
\prod_{j=r+1}^{n+1}\frac{b_1 \oplus x_{w(j)}}{b_1\ominus x_{w(j)}}
\times
\prod_{i=1}^{r-1}
[[x_{w(i)}|b]]^{\la_i}
\prod_{i=1}^{r-1}\left(
\prod_{j=i+1,\; j\neq r}^{n+1}\frac{x_{w(i)}\oplus x_{w(j)}}{x_{w(i)}\ominus x_{w(j)}}
\right),\label{eq:b1+b1}
\end{equation}
where we used an obvious equation
\begin{equation}
(x_{w(i)}\ominus b_1)
\left(
\frac{x_{w(i)}\oplus b_1}{x_{w(i)}\ominus b_1}\right)\label{eq:obvious}
=x_{w(i)}\oplus b_1.
\end{equation}
Define $w'\in S_n$ by $w'(i)=w(i)\;(1\leq i\leq r-1),$
$w'(i)=w(i+1)\;(r\leq i\leq n).$
Then $w'\in S_{n,r}$ and this correspondence
gives a bijection $S_{n+1,r}^I\rightarrow S_{n,r-1}.$
Then (\ref{eq:b1+b1}) is written as follows:
\begin{equation}
(b_1\oplus b_1)
\sum_{w'\in S_{n,r-1}}
\prod_{j=r}^{n}\frac{b_1 \oplus x_{w'(j)}}{b_1\ominus x_{w'(j)}}
\times \prod_{i=1}^{r-1}
[[x_{w'(i)}|b]]^{\la_i}
\prod_{i=1}^{r-1}\left(
\prod_{j=i+1}^{n}\frac{x_{w'(i)}\oplus x_{w'(j)}}{x_{w'(i)}\ominus x_{w'(j)}}
\right).\label{eq:s0-1}
\end{equation}

Define $S_{n+1,r}^{I\!I}$ to be the set of elements
$w\in S_{n+1,r}$ such that $w(n+1)=n+1.$ 
We calculate the corresponding part in (\ref{eq:pi0-1}), after substitution $x_{n+1}=b_1$. 
This, using (\ref{eq:obvious}) again,
is equal to 
\begin{equation}
\sum_{w\in S_{n+1,r}^{I\!I}}
\frac{x_{w(r)}\oplus b_1}{x_{w(r)}\ominus b_1}\times
\prod_{i=1}^{r}
[[x_{w(i)}|b]]^{\la_i}
\times \left(
\prod_{i=1}^{r}
\prod_{j=i+1}^{n}\frac{x_{w(i)}\oplus x_{w(j)}}{x_{w(i)}\ominus x_{w(j)}}\right).
\end{equation}

Let $F_\lambda(x)$ denote the last function.
The 
natural embedding $S_n\subset S_{n+1} $ given by $w(n+1)=n+1$ 
gives a bijection $S_{n,r}\cong S_{n+1,r}^{I\!I}.$ 
Using the identity 
\begin{equation}
\dfrac{x_{w(r)}\oplus b_1}{x_{w(r)}\ominus b_1}
-(1+\beta (b_1\oplus b_1))
=\dfrac{b_1\oplus b_1}{x_{w(r)}\ominus b_1}\label{eq:alpha0}
\end{equation}
we have 
\begin{eqnarray}
&&F_\lambda(x)-(1+\beta (b_1\oplus b_1))GQ_{\lambda}(x_1,\ldots,x_n|b_1,b_2,\ldots)
\nonumber\\
&=&{}
\sum_{w\in S_{n,r}}
\frac{b_1\oplus b_1}{x_{w(r)}\ominus b_1}
\prod_{i=1}^{r}
[[x_{w(i)}|b]]^{\la_i}
\times \left(
\prod_{i=1}^{r}
\prod_{j=i+1}^{n}\frac{x_{w(i)}\oplus x_{w(j)}}{x_{w(i)}\ominus x_{w(j)}}\right)
\nonumber \\
&=&{}
(b_1\oplus b_1)\sum_{w\in S_{n,r}}
\frac{x_{w(r)}\oplus x_{w(r)}}{x_{w(r)}\ominus b_1}
\left(
\prod_{j=r+1}^{n}\frac{x_{w(r)}\oplus x_{w(j)}}{x_{w(r)}\ominus x_{w(j)}}\right)
\nonumber \\
&&\times
\prod_{i=1}^{r-1}
[[x_{w(i)}|b]]^{\la_i} \left(
\prod_{i=1}^{r-1}
\prod_{j=i+1}^{n}\frac{x_{w(i)}\oplus x_{w(j)}}{x_{w(i)}\ominus x_{w(j)}}\right).
\label{eq:s0-3}
\end{eqnarray}
Now using (\ref{eq:s0-1}), (\ref{eq:s0-3}), and  
(\ref{eq:C}) in Lemma \ref{lem:Id_C},
with $m=n-r+1,\;t=b_1,$
we deduce (\ref{eq:pi0}).

Next suppose $s_0\lambda\geq \lambda,$ which is equivalent to
$\lambda_r\geq 2,$ where $r$ is the length of $\lambda.$
 We show
$$
{GQ_\lambda}(x_1,\ldots,x_n,b_1|\ominus b_1,b_2,\ldots)
=GQ_\lambda(x_1,\ldots,x_n|b_1,b_2,\ldots).
$$
In (\ref{eq:pi0-1}), 
those $w\in S_{n+1}$ with $1\leq w(n+1)\leq r$ will vanish
after the specialization $x_{n+1}=b_1.$
Using (\ref{eq:obvious}), we see that 
the remaining terms, i.e. those elements
satisfying $w(n+1)=n+1,$ 
coincide with the summands
of $GQ_\lambda(x_1,\ldots,x_n|b_1,b_2,\ldots)$ in the definition.

Next we consider type $B.$ 
Calculations for 
equations for $\pi_i\;(i\geq 1)$ 
are almost the same as for the case of type $C$.
Indeed, we now use $x\prod_{i=1}^{k-1}(x\oplus b_i)$ instead of
$[[x|b]]^k$. For the equation with respect to the operator $\pi_{0}$, 
if $s_0\la<\la$, 
we use the identity 
\begin{equation}
\dfrac{x_{w(r)}\oplus b_1}{x_{w(r)}\ominus b_1}
-(1+\beta b_1)
=\dfrac{b_1(2+\beta x_{w(r)})}{x_{w(r)}\ominus b_1}\label{eq:alpha0B}
\end{equation}
instead of (\ref{eq:alpha0}). Noting that $x_{w(r)}(2+\beta x_{w(r)})
=x_{w(r)}\oplus x_{w(r)}$ the rest is quite similar
using the same equation (\ref{lem:Id_C})
in the same way. The case of $s_0\la\geq \la$ is the same as for type $C.$

Finally we consider the case of type $D.$ 
We also use finite $n$ variables $x_1,\ldots,x_n$,
but note that $n$ should be even here.
Calculations for 
equations for $\pi_i\;(i\geq 2)$ 
are quite similar to the case of type $C$. 
For the equation with respect to the operator $\pi_{\hat{1}}$, 
we need to calculate the function
\begin{equation}
GP_\lambda(x_1,\ldots,x_{n},b_1,b_2|\ominus b_2,\ominus b_1,b_3,\ldots).\label{eq:pi1hatGP}
\end{equation}
The case $s_{\hat{1}}\lambda\geq \lambda$ is easy.
Indeed (\ref{eq:pi1hatGP})
is equal to $GP_\lambda(x_1,\ldots,x_{n}| b_1,b_2,b_3,\ldots).$
For the case $s_{\hat{1}}\lambda< \lambda,$
we can deduce
\begin{eqnarray*}
&&\frac{GP_\lambda(x_1,\ldots,x_{n},b_1,b_2|\ominus b_2,\ominus b_1,b_3,\ldots)
-(1+\beta (b_1\oplus b_2))
GP_\lambda(x_1,\ldots,x_{n}| b_1,b_2,b_3,\ldots)}{b_1\oplus b_2}\\
&=&GP_{\lambda^{(\hat{1})}}(x_1,\ldots,x_{n}| b_1,b_2,b_3,\ldots)
\end{eqnarray*}
by using the following identities 
\begin{equation}
\frac{x_{w(r)}\oplus b_2}{x_{w(r)}\ominus b_1}
-(1+\beta(b_1\oplus b_2))
=\frac{b_1\oplus b_2}{x_{w(r)}\ominus b_1},
\end{equation}
\begin{equation}
\sum_{i=1}^m\frac{u_i\oplus t}{u_i\ominus t}
\prod_{j\neq i}\frac{u_i\oplus u_j}{u_i\ominus u_j}
+\prod_{i=1}^m\frac{t\oplus u_i}{t\ominus u_i}=1,
\label{eq:Id_D}
\end{equation}where
the latter holds for {\it even\/} integer $m.$
The last equation is exactly Lemma \ref{var0-Iv} for $k=0$ and the number of 
variables is odd (with $n=m+1$, and $\{x_i\}$ being $u_1,\ldots,u_m,t$).

Finally we prove the equation for $\pi_1$ in type $D$.
If $\la$ is $1$-removable, 
the calculation is reduced to $\pi_1(x_1\oplus b_1)=1$,
as in the proof for $\pi_i\;(i\geq 1)$ of type $C.$
Let us consider the case
$s_1\la=\la.$
This occurs when (i) $r$ (the length of $\la)$ is odd and $\la_{r}=1$,
or (ii) $r$ is even and $\la_{r}\geq 2.$
We have to show that $GP_\la(x|b)=GD_\la(x|b)$ is 
symmetric in $b_1$ and $b_2.$
The case (ii) is obvious. Let us consider the case (i).
In view of Prop. \ref{prop:GPexpand},
it is enough to check this 
for $GP_\la(x_1,\ldots,x_n|b)$
 with $n=r+1.$
By the factorization theorem, we have 
$$
GP_\la(x_1,\ldots,x_n|b)
=GP_{\rho_{n-1}}(x_1,\ldots,x_n)G_\mu(x_1\ldots,x_n|b),
$$
where $\mu$ is a partition in $\mathcal{P}_{n-2}.$
Thus it suffices to show that  
$G_\mu(x_1,\ldots,x_n|b)$ is symmetric 
in $b_1$ and $b_2.$
Let $a_{ij}=[x_i|b]^{\la_j+n-j}(1+\beta x_i)^{j-1}$ 
and consider $\det(a_{ij}).$ For $j\leq n-2,$ 
the entries $a_{ij}$ are symmetric
in $b_1$ and $b_2$
because $\la_j+n-j\geq 2$. 
%The last two columns have entries 
%$((x_i\oplus b_1)(1+\beta x_i)^{n-2}, \;(1+\beta x_i)^{n-1}).$
The entries of the last two columns in $\det(a_{ij})$
are $((x_i\oplus b_1)(1+\beta x_i)^{n-2}, \;(1+\beta x_i)^{n-1}).$
We can replace them with 
$(x_i(1+\beta x_i)^{n-2}, \;(1+\beta x_i)^{n-1})$
by a column operation. 
Thus $\det(a_{ij})$ is symmetric in $b_1$ and $b_2.$
Hence the proof is complete.
$\qed$

\setcounter{equation}{0}
\section{Localization map}\label{sec:localization}
In this section, we define localization map $\Phi$. 
This map gives an 
injective $\R$-algebra homomorphism 
from  $G\Gamma_\R^X$ 
%the ring of $K$-supersymmetric functions
into $\Psi.$
The image of $GX_\la(x|b)$'s under the map $\Phi$ is 
shown to be a (unique) family of Schubert classes.

\subsection{Vanishing property}\label{ssec:vanishing}

Let $v$ be the Grassmannian element in $W$
corresponding to a strict partition $\mu,$
i.e. $\lambda(v)=\mu$ ($\la_+(v)=\mu$ for the type $D$ case) in the notation of \S \ref{ssec:SPGr}. 
Define a sequence $b_\mu$ of elements in $\R$ by 
$$
(b_\mu)_i= \begin{cases}
 b_{v(i)}& \mbox{if}\; v(i)\;\mbox{is negative}\\
0 & \mbox{otherwise}.
\end{cases}
$$ 
Explicitly the sequences $b_\mu$ for types $B,C$ are as follows: 
$$b_\mu=(\ominus b_{\mu_1},\ldots,\ominus b_{\mu_r},0,\ldots),$$
where $r$ is the length of $\mu,$
and for type $D$: 
$$b_\mu=\begin{cases}
(\ominus b_{\mu_1+1},\ldots,\ominus b_{\mu_r+1},0,\ldots) &
\mbox{if}\; r \mbox{ is even} \\
(\ominus b_{\mu_1+1},\ldots,\ominus b_{\mu_r+1},\ominus b_1,0,\ldots)&
\mbox{if}\; r \mbox{ is odd.}
\end{cases}
$$
\begin{rmk}
We defined 
$b_\mu$ for $\mu\in P_n$ in \S \ref{subsec:Gro} (type A case). 
Notify that the above definition is for types $B,C$, and $D.$
\end{rmk}

\begin{prop}[Vanishing property] \label{prop:vanishing}
Let $\lambda,\mu$ be strict partitions.
Then 
$GX_\lambda(b_\mu|b)=0$ unless $\lambda\subset\mu$
and 
$GX_\lambda(b_\la|b)=\prod_{\alpha\in \mathrm{Inv}(\lambda)}e(-\alpha).$ 
\end{prop}

{\it Proof.} We only show this for $GP_\lambda(x|b)$ and 
the case when the length $r$ of $\mu$ is even. 
Other cases are similar. Using the stability property, $GP_\lambda(b_\mu|b)$ can be 
evaluated as
$GP_{\lambda}(\ominus b_{\mu_1+1},\ldots,\ominus b_{\mu_r+1}|b).$

Suppose $\lambda\not\subset\mu.$ So we have
$\mu_j<\lambda_j$ for some $1\leq j\leq r.$ 
Then we easily see that the polynomial
$\prod_{i=1}^r[x_{w(i)}|b]^{\lambda_i}$  
vanishes for any $w\in S_r$ when we specialize
$x_i$ to $\ominus b_{\mu_i+1}\;(1\leq i$
$\leq r).$ 
Note that the denominator of the factor 
$$\prod_{i=1}^r\prod_{j=i+1}^n\frac{x_{w(i)}\oplus x_{w(j)}}{x_{w(i)}\ominus x_{w(j)}}$$
does not vanish identically under the specialization.
Thus we have $GP_\lambda(b_\mu|b)=0.$

Next we consider the case 
$\lambda=\mu.$ We see that the terms other than 
the one comes from $w=e$ vanish, by a similar reason to the previous case.
The term corresponding to $w=e$ is evaluated by using the following 
equation:
$$
(\ominus b_{\lambda_i+1}\oplus b_1)\cdots (\ominus b_{\lambda_{i}+1}\oplus b_{\lambda_i})
\prod_{k=i+1}^r\frac{\ominus b_{\lambda_{i}+1} \ominus b_{\lambda_{k}+1}}{\ominus b_{\lambda_{i}+1} \oplus b_{\lambda_{k}+1}}
$$
$$
=\prod_{k\in \{1,\ldots,\lambda_i\}\setminus\cup_{k=i+1}^r
\{\lambda_{k}+1\}}(\ominus b_{\lambda_i+1}\oplus b_k)
\prod_{k=i+1}^r(\ominus b_{\lambda_i+1}\ominus b_{\lambda_k+1}),
$$
which is a simple consequence of cancellation. %?
This last expression is readily seen to be $\prod_{\alpha\in \mathrm{Inv}(\lambda)}e(-\alpha).$
$\qed$

\subsection{Algebraic localization map $\Phi$}
Let $\mu$ be a strict partition. 
Define an $\R$-algebra 
homomorphism 
$\Phi_\mu: G\Gamma_\R\rightarrow \R$ by 
``substitution'' $x=b_\mu.$
This is well-defined.
In fact, we know that arbitrary element $F(x)$ in $G\Gamma_\R$ can be written as a possibly
infinite $\R$-linear combination
$ F(x)=\sum_{\lambda\in \mathcal{SP}} c_\lambda \cdot GX_\lambda(x|b)$
of $GX_\lambda$'s. 
Since $GX_\lambda({b}_\mu|b)$ can be non-zero 
for only finitely many $\lambda$'s such that 
$\lambda\subset \mu,$ 
$F(b_\mu)=\sum_{\lambda\in \mathcal{SP}} c_\lambda \cdot GX_\lambda({b}_\mu|b)
$ is a well-defined element in $\R.$

\begin{Def}[Localization map] Define the homomorphism of $\R$-algebras as
$$\Phi: G\Gamma_\R^X\longrightarrow \mathrm{Fun}(\mathcal{SP},\R)\quad 
f\mapsto (\mu\mapsto \Phi_{\mu}(f)).$$
$\Phi$ is called the {\rm localization map} for $G\Gamma_\R^X.$
\end{Def}

We often call $\Phi$ the {\it algebraic} localization map
in order to form 
a contrast to the {\it geometric} one. See \S \ref{sec:geometry} for geometric background of the map $\Phi$.

\begin{prop}\label{InPsi}
We have $\mathrm{Im}(\Phi)\subset \Psi.$
\end{prop}

This is a consequence of the following 
more explicit statement.

\begin{lem}
Let 
$t_{ij}=s_{t_j-t_i}, \;s_{ij}=s_{t_j+t_i}\;(j>i\geq 1),\;s_{ii}=s_{t_i}=s_{2t_i}.$
If 
$F(x)$ is %a $K$-supersymmetric function 
in $G\Gamma$ then, for any 
$\mu\in \mathcal{SP}$, we have
\begin{eqnarray}
F(b_{t_{ij}\mu})-F(b_\mu)\in \langle b_j\ominus b_i\rangle
\quad \mbox{for}\;j>i\geq 1,\label{eq:tij}\\
F(b_{s_{ij}\mu})-F(b_\mu)\in \langle b_j\oplus b_i\rangle
\quad \mbox{for}\;j>i\geq 1,\label{eq:sij}\\
F(b_{s_{ii}\mu})-F(b_\mu)\in \langle b_i\rangle\quad \mbox{for}\;
i\geq 1.\label{eq:sii1}
\end{eqnarray}
If $F(x)$ is in $G\Gamma_{+}$ then we have
\begin{equation}
F(b_{s_{ii}\mu})-F(b_\mu)\in \langle b_i\oplus b_i\rangle
\quad \mbox{for}\;
i\geq 1.\label{eq:sii2}\end{equation}
\end{lem}
{\it Proof.} 
There are the following possibilities:  
(i) $i,j\in \lambda$, (ii) $i,j\notin \lambda$, $i\in \lambda$ and 
$j\notin\lambda$, (iv) $i\notin \lambda$ and 
$j\in\lambda.$
Actions of $t_{ij}$ and $s_{ij}$ are 
given as follows: 
$$
\begin{array}{|c|c|c|c|c|}\hline
 & i,j\in \lambda & i,j\notin \lambda
 & i\in \lambda, j\notin \lambda
 & i\notin \lambda, j\in \lambda\\\hline 
t_{ij}\cdot\lambda & \lambda & \lambda & (\lambda\setminus\{i\})\cup\{j\} & (\lambda\setminus\{j\})\cup\{i\} \\
\hline
s_{ij}\cdot\lambda & \lambda\setminus\{i,j\} & \lambda\cup\{i,j\} & \lambda & \lambda \\\hline
\end{array}
$$
We may concentrate on two variables $x_i,x_j.$ 
The check for (\ref{eq:tij})  
is easy in view of 
the fact that $F$ is symmetric (note that $\langle e(-\alpha)
\rangle
=\langle e(\alpha) \rangle).$
To show (\ref{eq:sij}), it suffices to consider the cases (i) and (ii)
(the cases (iii) and (iv) are obvious).
Let 
$G=F(\ominus b_i,\ominus b_j)-F(0,0).$
Then 
$G$ as a polynomial of $b_i$ vanishes at
$\ominus b_j,$ because $G\in G\Gamma$. This implies that 
$G$ is divisible by $b_j \oplus b_i$ (cf. Lemma \ref{lem:2var}). 
Hence we have (\ref{eq:sij}).

Next we show (\ref{eq:sii1}) and (\ref{eq:sii2}). It is obvious that 
$F(\ominus b_i)-F(0)$ is divisible by $\ominus b_i.$
The divisibility for $F(\ominus b_i)-F(0)$ by $b_i\oplus b_i$ is the very condition that 
$F$ is a member of $G\Gamma_{+}.$
This concludes the proof.
$\qed$

\bigskip

The following is the main result of this paper.
\begin{thm}\label{thm:main} 
$\{\Phi(GX_\lambda(x|b))\}_{\lambda\in \mathcal{SP}}$   
is a family of Schubert classes. 
\end{thm}
{\it Proof.} 
By the uniqueness of Schubert classes (Prop. \ref{prop:unique}),
it suffices to check the defining property. 
We know that $\Phi(GX_\lambda(x|b))\in \Psi$ by  Prop. \ref{InPsi}.
Note that $\Phi\circ\pi_i=\pi_i\circ\Phi$ holds (cf. \cite{DSP}, Prop. 7.4.).
By this fact and Thm.\:\ref{thm:piGQ}, 
the divided difference equation (\ref{eq:divdiffSchubert}) is satisfied.
The initial condition (\ref{eq:InitialSchubert}) is satisfied
because $GX_\lambda(x|b)$ vanishes at $x=b_\emptyset
=(0,0,\ldots)$ for $\lambda\neq \emptyset,$
and we have $GX_\emptyset(x|b)=1.$ $\qed$

\begin{cor}\label{cor:injective} 
$\Phi$ is injective.
\end{cor}

\begin{rmk}\label{inj} One can prove the injectivity of $\Phi$ in 
the same way as \cite{DSP}, Lemma 6.5.
\end{rmk}

\setcounter{equation}{0}
\section{Equivariant $K$-theory of the maximal isotropic Grassmannians}
\label{sec:geometry}
In this section, we show an application of 
the $K$-theoretic factorial $P$- and $Q$-functions
to the equivariant $K$-theory of the maximal isotropic 
Grassmannians.
\subsection{Maximal isotropic Grassmannians} Let $n$ be a positive integer.
Suppose $G$ is one of the groups $SO(2n+1,\C),Sp(2n,\C),$ or $SO(2n+2,\C).$ 
Let $T$ be a maximal torus of $G,$ and let 
$B$ be a Borel subgroup containing $T.$
Then we have the corresponding root system, 
which 
we identify with the one defined in \ref{ssec:roots} for 
types $B_n,C_n,D_{n+1}$ respectively.
Let $P$ be the standard maximal parabolic 
subgroups corresponding to the simple root $\alpha_0$ for $B_n$
and $C_n$, 
and $\alpha_{\hat{1}}$ for $D_{n+1}$.
Let $\mathcal{G}_n$ denote the homogeneous variety $G/P.$
Then $\mathcal{G}_n$ 
is $OG(n,2n+1),LG(n),OG(n+1,2n+2)$ 
respectively for type $B_n,C_n,$ and $D_{n+1},$ where 
$LG(n)$ denote 
the Lagrangian Grassmannian, and $OG(n,2n+1),\;OG(n+1,2n+2)$
the maximal odd and even orthogonal Grassmannians. 
It is known that $OG(n,2n+1)$ and
$OG(n+1,2n+2)$ are isomorphic as algebraic varieties, however, note that  
here we are considering actions of different tori.

We denote by $\!\mathcal{G}_n^T$ the set of $T$-fixed points in $\mathcal{G}_n.$
We have a natural bijection $\mathcal{G}_n^T\cong \mathcal{SP}(n),$
where $ \mathcal{SP}(n)$ is the set of strict partitions
$\lambda$ such that $\lambda\subset \rho_n.$
Let $e_\lambda\in \mathcal{G}_n^T$ denote the $T$-fixed point 
corresponding 
to $\lambda\in \mathcal{SP}(n).$
Let $B_{-}$ be the opposite Borel subgroup, so that we have $B\cap B_{-}=T.$
The {\it Schubert varieties} $\Omega_\lambda$ in $\mathcal{G}_n$ 
is defined to be the closure of the orbit $B_{-}e_\lambda\subset G/P.$
Then the 
co-dimension of $\Omega_\lambda$ is given by $|\lambda|.$
Note that $e_\mu\in \Omega_\lambda$ if and only of $\lambda\subset \mu.$

Let $K_{T}(\mathcal{G}_n)$ denote
the Grothendieck group of the abelian category of  
$T$-equivariant coherent sheaves on
$\mathcal{G}_n.$ It is known that 
$K_{T}(\mathcal{G}_n)$ has a natural ring structure.
Since $\Omega_\lambda$ is $T$-stable, its structure sheaf 
$\mathcal{O}_{\Omega_\lambda}$ defines a class $[\mathcal{O}_{\Omega_\lambda}]_T$
 in $K_{T}(\mathcal{G}_n).$ 
$K_{T}(\mathcal{G}_n)$ 
has a natural $R(T)$-algebra structure, where
$R(T)$ is the representation ring of $T$, which is also identified with $K_T(pt).$
As an $R(T)$-module 
$K_{T}(\mathcal{G}_n)$
is free and has basis consisting of the Schubert 
structure sheaves $[\mathcal{O}_{\Omega_\lambda}]_{T}$
($\lambda\in \mathcal{SP}(n)$).

\subsection{GKM ring}

Let $\Delta(n)^+$ denote the set of 
positive roots associated to the pair $(G,B).$
Then $\Delta(n)^+$ is $\Delta_n^+$
for types $B,C,$ and 
$\Delta_{n+1}^+$
for type $D$ (see \S \ref{ssec:roots} for the notation $\Delta_n^+$). We realize  
$R(T)$ as the ring of Laurent polynomials : 
$R(T)=\Z[e^{\pm t_1},\ldots,e^{\pm t_n}]$ for types
$B,C$ and 
$R(T)=\Z[e^{\pm t_1},\ldots,e^{\pm t_{n+1}}]$ for type $D.$ 
Let us denote by $\mathrm{Fun}(\mathcal{SP}(n),R(T))$ the set of 
all functions from $\mathcal{SP}(n)$ to $R(T).$
This is naturally an $R(T)$-algebra by pointwise multiplication.
\begin{Def} Let $\Psi_n$ be the $R(T)$-subalgebra of 
$\mathrm{Fun}(\mathcal{SP}(n),R(T))$
defined as follows: 
a map $\psi: \mathcal{SP}(n)\rightarrow R(T)$ is in $\Psi_n$ if and only if  
$$
\psi(s_\alpha \mu)-\psi(\mu)\in (1-e^{\alpha})R(T)\quad
\mbox{for all}\;
\mu \in \mathcal{SP}(n),\;
\alpha\in \Delta(n)^+.
$$
\end{Def}

Letting $\beta=-1$ we consider $R(T)$ as a subalgebra of $\Z\otimes_{\Z[\beta]}\R
=\Z[e^{\pm t_1},e^{\pm t_2},\ldots].$
If $\lambda,\mu$ are strict partitions in $\mathcal{SP}(n),$
we see that $\psi_\lambda(\mu)\in R(T).$  
Thus we have $\psi_\lambda^{(n)}:=\psi_\lambda|_{\mathcal{SP}(n)}
\in \mathrm{Fun}(\mathcal{SP}(n),R(T)),$ which 
is obviously an element in $\Psi_n.$ 
Thus we have the following.
\begin{prop} We have 
$\Psi_n=\bigoplus_{\la\in \mathcal{SP}(n)}R(T)\psi_\lambda^{(n)}.$
\end{prop}

\subsection{$K_T(\mathcal{G}_n)$ and its Schubert basis}\label{ssec:K_T}

Note that $K_T(\mathcal{G}_n^T)
\cong \prod_{e_\mu\in \mathcal{G}_n^T}
K_T(e_\mu)$  is naturally 
identified with $\mathrm{Fun}(\mathcal{SP}(n),R(T)),$
since we have $\mathcal{G}_n^T\cong \mathcal{SP}(n).$

\begin{thm}[\cite{KK}]\label{thm:KK} Let $i$ be the inclusion map $\mathcal{G}_n^T
\hookrightarrow \mathcal{G}_n.$
The induced map 
$$
i^*: K_T(\mathcal{G}_n)\longrightarrow K_T(\mathcal{G}_n^T)\cong 
\mathrm{Fun}(\mathcal{SP}(n),R(T))
$$
is injective, and its image is equal to the $R(T)$-subalgebra $\Psi_n.$
\end{thm}
{\it Proof.} For the proof, see \cite{KK}, Thm.\:3.13 and Corollary 3.20.
$\qed$

\begin{thm}[\cite{KK}]
The image 
$i^* [\mathcal{O}_{\Omega_\lambda}]_T,\;\lambda
\in\mathcal{SP}(n), $ is equal to 
$\psi_\lambda^{(n)}\in \Psi_n,$
where $\psi_\lambda^{(n)}$ is the restriction of 
$\psi_\lambda$ to $\mathcal{SP}(n).$
\end{thm}

An $R(T)$-basis $\{\tau^w\}_{w\in W}$ for $K_T(G/B)$ 
was constructed in \cite{KK}.
The class 
$\tau^w$ is closely related to 
$[\mathcal{O}_{\Omega_w}]_T$ but not-exactly coincides with it,
where $\Omega_w$ is the Schubert variety in $G/B$ 
corresponding to $w\in W.$
For the precise comparison
with the class $[\mathcal{O}_{\Omega_w}]_T$, 
see e.g. \cite{LSS}, \cite{GrKu}. 
Now we need the corresponding result for the parabolic case, i.e. for $K_T(G/P).$
The reader can consult \cite{LSS} for the result.

\bigskip

For $X=B,C,D$, we denote $G\Gamma_n^B=G\Gamma_n^D=G\Gamma_n$ and
$G\Gamma_n^C=G\Gamma_{n,+}.$
We consider $R(T)$-algebras 
$R(T)\otimes_{\Z[\beta]} G\Gamma_n^X,$
with the specialization $\beta=-1.$
For any strict partition $\lambda$ in $\mathcal{SP}_n,$ 
let $GX_\lambda^{(n)}(x|1-e^t)$ be functions
in $R(T)\otimes_{\Z[\beta]} G\Gamma_n^X,$
defined as follows:
$$
\begin{cases}
GB_\lambda^{(n)}(x|1-e^{t_1},\ldots,1-e^{t_n},0,\ldots)& \mbox{if}\;X=B,\\
GC_\lambda^{(n)}(x|1-e^{t_1},\ldots,1-e^{t_n},0,\ldots)& \mbox{if}\;X= C,\\
GD_\lambda^{(n)}(x|1-e^{t_1},\ldots,1-e^{t_n},1-e^{t_{n+1}},0,\ldots)& \mbox{if}\; X=D.
\end{cases}
$$
See \S \ref{ssec:GGamma_R} for the notation $GX_\la^{(n)}(x|b).$
Then $\{GX_\lambda^{(n)}(x|1-e^t)\;|\;\la\in \mathcal{SP}_n\}$
form an $R(T)$ basis of $R(T)\otimes_{\Z[\beta]} G\Gamma_n^X.$

\begin{thm}\label{thm:main-geom} There exists a surjective homomorphism of 
$R(T)$-algebras
$$\pi_n: R(T)\otimes_{\Z[\beta]} G\Gamma_n^X\longrightarrow K_{T}(\mathcal{G}_n),$$ 
which sends 
$GX_\lambda^{(n)}(x|1-e^t)$ to $[\mathcal{O}_{\Omega_\lambda}]_T$ if $\lambda\subset \rho_n$
and to $0$ if $\lambda\not\subset \rho_n.$
\end{thm}
{\it Proof.} The algebraic localization map
$\Phi$ naturally 
induces a homomorphism of $R(T)$-algebras
$\Phi_n: R(T)\otimes_{\Z[\beta]} G\Gamma_n^X
\rightarrow \Psi_n,$
given by $\Phi_n(F)=(\Phi_\mu(F))_{\mu\in \mathcal{SP}(n)}.$
By Thm.\:\ref{thm:main} we know that $\Phi_n$ maps
$GX_\lambda^{(n)}(x|1-e^t)$ to $\psi_\lambda^{(n)}$ 
if $\lambda\subset \rho_n$
and to $0$ if $\lambda\not\subset \rho_n.$
Then there is a unique $R(T)$-algebra homomorphism $\pi_n$ making 
the following diagram commutative:
$$
\begin{CD}
R(T) \otimes_{\Z[\beta]} G\Gamma_n^X @>\pi_n >>K_T(\mathcal{G}_n) \\
@V\Phi_nVV @ VVi^*V\\
\Psi_n @=\Psi_n
\end{CD}
$$
Since $i^*$ is an isomorphism (Thm.\:\ref{thm:KK}), $\pi_n$ is determined by the 
above diagram, and has the desired property.
$\qed$

\bigskip

Now we 
 derive the non-equivariant analogue of Thm.\:\ref{thm:main-geom}.
We first recall that two cases $B_n$ and $D_{n+1}$ coincide in non-equivariant setting.

Let  
$\ep: R(T)\rightarrow \Z$ 
be the ring homomorphism 
given by the evaluation at the identity element of $T$; it is given 
by $e^{t_i}\mapsto 1\;(1\leq i\leq r)$ or equivalently $b_i\mapsto 0$.
Let 
$GB_\lambda^{(n)}(x)=GD^{(n)}_\la(x)=GP_\la(x_1,\ldots,x_n)$
and $GC_\lambda^{(n)}(x)=GQ_\lambda^{(n)}(x_1,\ldots,x_n).$

\begin{cor}\label{thm:non-eq} There exists a surjective homomorphism of 
rings
$$\pi_n: G\Gamma_n^X\longrightarrow K(\mathcal{G}_n),$$ 
which sends 
$GX_\lambda^{(n)}(x)$ to $[\mathcal{O}_{\Omega_\lambda}]$ if $\lambda\subset \rho_n$
and to $0$ if $\lambda\not\subset \rho_n.$
\end{cor}
{\it Proof.} This follows from 
Thm.\:\ref{thm:main-geom}, since 
the canonical map   
\begin{equation}
\Z\otimes_{R(T)}K_T(\mathcal{G}_n)\rightarrow
K(\mathcal{G}_n),\label{eq:forget}
\end{equation}
is an isomorphism of rings, 
where $\Z$ is considered as an $R(T)$-algebra by $\ep$  
(see \cite{KK}, Proposition 3.25). $\qed$

\setcounter{equation}{0}
\section{Combinatorial expressions}\label{sec:comb} 

In this section, we present two combinatorial 
expressions for $GQ_\lambda(x|b)$ and $GP_\lambda(x|b)$
(Thm.\:\ref{thm:SVT} and Thm.\:\ref{thm:EYD}).
The proof will be given in \S \ref{sec:comb_proof}.

\subsection{Shifted set-valued tableaux}\label{ssec:SVT}

For each strict partition $\lambda$, 
let $\mathbb{D}(\lambda)$ be the shifted diagram of $\lambda,$
which is defined as   
$\mathbb{D}(\lambda)=\{(i,j)\in \Z\times \Z\;|\; i\leq j\leq \lambda_i+i-1,\;
1\leq i\leq \ell(\lambda)\}.$

\bigskip

Let $\A$ denote the ordered alphabet
$${1'<1<2'<2<\cdots<n'<n}.$$
Let $\X$ denote the set of all non-empty subsets of 
$\A.$

A {\it set-valued shifted tableau} $T$ of shape $\lambda$
is an assignment $T:\mathbb{D}(\lambda)\to \X$ such that

(1) $\max\, T(i,j)\leq \min\, T(i,j+1),\; \max\, T(i,j)\leq \min\, T(i,j+1),$

(2) Each $a\, (=1,2,\ldots,n)$ appears at most once in each 
column,

(3) Each $a'\,(=1',2',\ldots,n')$ appears at most once in each 
row.

\bigskip

We denote by $\T(\lambda)$ the set of all set-valued shifted tableaux of shape $\lambda.$
Define $\T'(\lambda)$ to be the subset of $\T(\lambda)$ consisting of 
$T\in  \T(\lambda)$ satisfying 
the condition that 
for each $i\;(1\leq i\leq n),$
$T(i,i)$ contains only $1,2',3,4',\ldots$ if $i$ is odd,
and $1',2,3',4,\ldots$ if $i$ is even.
Define $\T''(\lambda)$ to be the subset of $\T(\lambda)$ consisting of 
$T\in  \T(\lambda)$ satisfying 
the condition that 
for each $i\;(1\leq i\leq n),$
$T(i,i)$ contains no primed symbols $\{1',2',\ldots,n'\}.$

\begin{example}\label{ex:tab}
{\rm
The following
are examples of set-valued shifted tableaux
of shape $\lambda=(4,3,1)$ : 
\begin{align*}
T_1&=\Tableau{1'&12'&23&34'\\~&2&4'&6\\~&~&6'}, &
T_2&=\Tableau{1&12'&2&2\\~&23'&3&3\\~&~&4'},&
T_3&=\Tableau{1&12'&23&34'\\~&2&4'&6\\~&~&6}.
\end{align*}}
Of these three tableaux in $\mathcal{T}(\la),$ 
only $T_2$ is in $\mathcal{T}'(\la),$
and only $T_3$ is in $\mathcal{T}''(\la).$
\end{example}

For each integers $i,j$ such that $1\leq i\leq n,\; i\leq j$ and
$a\in \A$, we define 
$$
w_I(i,j;a)
=\begin{cases}
x_{a}\oplus b_{j-i}& \mbox{if}\quad a\in \{1,2,\ldots,n\}\\
x_{|a|}\ominus b_{j-i}& \mbox{if}\quad a\in \{1',2',\ldots,n'\},
\end{cases}
$$
where $b_0=0$ and
$$
w_{I\!I}(i,j;a)
=\begin{cases}
x_{a}\oplus b_{j-i+1}& \mbox{if}\quad a\in \{1,2,\ldots,n\}\\
x_{|a|}\ominus b_{j-i+1}& \mbox{if}\quad a\in \{1',2',\ldots,n'\}.
\end{cases}
$$
For $T\in \mathcal{T}(\la)$ we define 
\begin{equation}
(x|b)_I^T=\prod_{(i,j)\in \mathbb{D}(\lambda),\; a\in T(i,j)}w_I(i,j;a),
\quad
(x|b)_{I\!I}^T=\prod_{(i,j)\in \mathbb{D}(\lambda),\; a\in T(i,j)}w_{I\!I}(i,j;a).
\label{eq:(x|b)^T}
\end{equation}
Recall that $GX_\lambda^{(n)}(x|b)$ is defined in \S \ref{ssec:GGamma_R}

\begin{thm}[Combinatorial formula in terms of set-valued tableaux]\label{thm:SVT}
Let $\lambda$ be a strict partition of length $\leq n.$
Then we have
\begin{align}&\mbox{\rm{Type}}\; B:\quad
GB_\lambda^{(n)}(x|b)
=
\sum_{T\in \T''(\lambda)}
\beta^{|T|-|\lambda|} 
(x|b)_I^T,
\\
&\mbox{\rm{Type}}\; C:\quad GC_\lambda^{(n)}(x|b)
=
\sum_{T\in \T(\lambda)}
\beta^{|T|-|\lambda|} 
(x|b)_I^T,
\label{eq:TabGQ}
\\
&\mbox{\rm{Type}}\; D:\quad
GD_\lambda^{(n)}(x|b)
=
\sum_{T\in \T'(\lambda)}
\beta^{|T|-|\lambda|} 
(x|b)_{I\!I}^T.
\end{align}
\end{thm}

For example, let $T_1$ be the tableau in Example \ref{ex:tab}, then 
the corresponding summand on the right-hand side of (\ref{eq:TabGQ}) equals:
$$
x_1(x_1\oplus b_1)x_2(x_2\ominus b_1)(x_2\oplus b_2)
(x_3\oplus b_2)(x_3\oplus b_3)(x_4\ominus b_1)(x_4\ominus b_3)
x_6(x_6\oplus b_2).
$$

\subsection{Excited Young diagrams}\label{ssec:EYD}
Let $\mathcal{D}_n$ denote the subset of $\Z\times \Z$
consisting of $(i,j)$ satisfying $1\leq i\leq n,\; i\leq j.$
We call an element $(i,j)$ in $\mathcal{D}_n$ a {\it box}.
Let $D$ be a finite set of boxes, i.e. a finite subset of $\mathcal{D}_n.$
Suppose $(i,j)\in D$ and $(i+1,j),(i+1,j+1),(i,j+1)\notin D$
(if $i=j$, the condition $(i+1,j)\notin D$ is vacuous).
Then it is said that 
$D'=(D\setminus (i,j))\cup (i+1,j+1)$
is obtained from $D$ by 
an {\it elementary excitation}, and we denote $D\rightarrow D'.$
Let $\mathcal{E}^I_n(\lambda)$ denote the set of 
all subsets $D$ of $\mathcal{D}_n$ 
obtained from $D_0=\D(\lambda)$ by 
a sequence of successive elementary excitations 
$$
D_0\rightarrow D_1\rightarrow 
\cdots\rightarrow D_r=D.
$$

Let $D\in \mathcal{E}^I_n(\lambda).$
We denote by $B^I(D)$ the set of boxes
$(i,j)$ satisfying the following conditions:

(1) $(i,j)\notin D,$

(2) there is a positive integer $k$ such that $(i+k,\;j+k)\in D,$

(3) $(i+s,j+s)\notin D$ for $1\leq s<k,$

(4) $(i+s,j+s+1)\notin D$ for $0\leq s< k,$

(5) $(i+s-1,j+s)\notin D$ for $0\leq s< k$ (this condition is ignored 
if $i=j).$

\begin{example}\label{ex:EYD} Let
$n=4.$ The following are examples of elements 
of $\mathcal{E}^{I}_4(\lambda)$ for $\lambda=(4,2).$
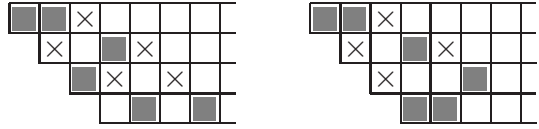
\begin{center}
\begin{picture}(205,45)
\put(0,40){\line(1,0){75}}
\put(0,30){\line(1,0){75}}
\put(10,20){\line(1,0){65}}
\put(20,10){\line(1,0){55}}
\put(30,0){\line(1,0){45}}
\put(0,30){\line(0,1){10}}
\put(10,20){\line(0,1){20}}
\put(20,10){\line(0,1){30}}
\put(30,0){\line(0,1){40}}
\put(40,0){\line(0,1){40}}
\put(50,0){\line(0,1){40}}
\put(60,0){\line(0,1){40}}
\put(70,0){\line(0,1){40}}
\put(1,31){$\cell$}
\put(11,31){$\cell$}
\put(31,21){$\cell$}
\put(21,11){$\cell$}
\put(41,1){$\cell$}
\put(61,1){$\cell$}
\put(51,12){$\times$}
\put(41,22){$\times$}
\put(21,32){$\times$}
\put(31,12){$\times$}
\put(11,22){$\times$}
\put(100,40){\line(1,0){75}}
\put(100,30){\line(1,0){75}}
\put(110,20){\line(1,0){65}}
\put(120,10){\line(1,0){55}}
\put(130,0){\line(1,0){45}}
\put(100,30){\line(0,1){10}}
\put(110,20){\line(0,1){20}}
\put(120,10){\line(0,1){30}}
\put(130,0){\line(0,1){40}}
\put(140,0){\line(0,1){40}}
\put(150,0){\line(0,1){40}}
\put(160,0){\line(0,1){40}}
\put(170,0){\line(0,1){40}}
\put(101,31){$\cell$}
\put(111,31){$\cell$}
\put(131,21){$\cell$}
\put(131,1){$\cell$}
\put(141,1){$\cell$}
\put(151,11){$\cell$}
\put(121,12){$\times$}
\put(111,22){$\times$}
\put(121,32){$\times$}
\put(141,22){$\times$}
\end{picture}
\end{center}
The boxes with symbol $\times$ are the elements of $B^I(D).$
\end{example}

Next we consider type $D.$
In the sequel, we assume $n$ is even, whenever we consider 
type $D$ case. In order to treat the type $D$ case,
we define  
$\mathcal{E}^{I\!I}_n(\lambda)$, which is a subset of 
$\mathcal{E}^{I}_n(\lambda)$ consisting of EYD 
whose way of excitation of each diagonal box $(i,i)$ is 
restricted in the following manner. 
%to be even steps. More precisely, 
The elementary excitation of $(i,i)\in D$ is defined as follows:  
if $(i,i)\in D$ and $(i,i+1),(i+1,i+1),(i+1,i+2),(i+2,i+2)\notin D$
then $D'=(D\setminus (i,i))\cup (i+2,i+2)$ is an excitation 
of type $I\!I.$
For off-diagonal boxes the definition of excitation is 
the same. For example the left diagram of Example \ref{ex:EYD} is in $\mathcal{E}^{I\!I}_n(\lambda),$
but the right one is not.

Let $D\in \mathcal{E}^{I\!I}_n(\lambda).$ We define the subset  
$B^{I\!I}(D)$ of $B^{I}(D)$ 
by restricting $k$ to be even in the definition of $B^I(D)$.
\begin{example}\label{ex:EYDs}
In the following figure, the boxes with symbol $\times$ are the elements of $B^{I\!I}(D).$
\begin{center}
\begin{picture}(105,40)
\put(0,40){\line(1,0){75}}
\put(0,30){\line(1,0){75}}
\put(10,20){\line(1,0){65}}
\put(20,10){\line(1,0){55}}
\put(30,0){\line(1,0){45}}
\put(0,30){\line(0,1){10}}
\put(10,20){\line(0,1){20}}
\put(20,10){\line(0,1){30}}
\put(30,0){\line(0,1){40}}
\put(40,0){\line(0,1){40}}
\put(50,0){\line(0,1){40}}
\put(60,0){\line(0,1){40}}
\put(70,0){\line(0,1){40}}
\put(1,31){$\cell$}
\put(11,31){$\cell$}
\put(31,21){$\cell$}
\put(31,1){$\cell$}
\put(41,1){$\cell$}
\put(51,11){$\cell$}
\put(11,22){$\times$}
\put(21,32){$\times$}
\put(41,22){$\times$}
\end{picture}
\end{center}
\end{example}

Define for $(i,j)\in \mathcal{D}_n$, the following weight of types $B,C,$ and $D:$
\begin{equation}
wt^B(i,j)=\begin{cases}
x_i & \mbox{if}\quad j\leq n\\
x_i\oplus b_{j-n} & \mbox{if}\quad j>n \end{cases},\quad
wt^C(i,j)=\begin{cases}
x_i\oplus x_j & \mbox{if}\quad j\leq n\\
x_i\oplus b_{j-n} & \mbox{if}\quad j>n \end{cases}.
\label{eq:weightC}
\end{equation}
For type $D$, note that we assumed $n$ is even.
We define
$$
wt^{D}(i,j)=\begin{cases}
x_i\oplus x_{j+1} &\mbox{if}\quad  j\leq n-1\\
x_i\oplus b_{j-n+1} & \mbox{if}\quad j\geq n
\end{cases}.
$$

\begin{example}
Let $n=4.$ The weight function $wt^C$ on $\mathcal{D}_4$ 
is given as follows: 
 
 \setlength{\unitlength}{0.7mm}
 \begin{picture}(300,55)
\put(0,50){\line(1,0){148}}
\put(0,40){\line(1,0){148}}
\put(20,30){\line(1,0){128}}
\put(40,20){\line(1,0){108}}
\put(0,40){\line(0,1){10}}
\put(20,30){\line(0,1){20}}
\put(40,20){\line(0,1){30}}
\put(60,20){\line(0,1){30}}
\put(60,10){\line(0,1){40}}
\put(80,10){\line(0,1){40}}
\put(100,10){\line(0,1){40}}
\put(120,10){\line(0,1){40}}
\put(140,10){\line(0,1){40}}
\put(60,10){\line(1,0){85}}
\put(2,43){$x_1\!\oplus\! x_1$}
\put(22,43){$x_1\!\oplus\! x_2$}
\put(42,43){$x_1\!\oplus\! x_3$}
\put(22,33){$x_2\!\oplus\! x_2$}
\put(42,33){$x_2\!\oplus\! x_3$}
\put(42,23){$x_3\!\oplus\! x_3$}
\put(62,43){$x_1\!\oplus\! x_4$}
\put(82,43){$x_1\!\oplus\! b_1$}
\put(102,43){$x_1\!\oplus\! b_2$}
\put(62,33){$x_2\!\oplus\! x_4$}
\put(82,33){$x_2\!\oplus\! b_1$}
\put(102,33){$x_2\!\oplus\! b_2$}
\put(62,23){$x_3\!\oplus\! x_4$}
\put(82,23){$x_3\!\oplus\! b_1$}
\put(102,23){$x_3\!\oplus\! b_2$}
\put(122,43){$x_1\!\oplus\! b_3$}
\put(122,33){$x_2\!\oplus\! b_3$}
\put(122,23){$x_3\!\oplus\! b_3$}
\put(62,13){$x_4\!\oplus\! x_4$}
\put(82,13){$x_4\!\oplus\! b_1$}
\put(102,13){$x_4\!\oplus\! b_2$}
\put(122,13){$x_4\!\oplus\! b_3$}
\put(143,13){$\cdots\;,$}
\put(143,23){$\cdots$}
\put(143,33){$\cdots$}
\put(143,43){$\cdots$}
\end{picture}
while $wt^D$ is as follows:

 \setlength{\unitlength}{0.7mm}
 \begin{picture}(300,55)
\put(0,50){\line(1,0){148}}
\put(0,40){\line(1,0){148}}
\put(20,30){\line(1,0){128}}
\put(40,20){\line(1,0){108}}
\put(0,40){\line(0,1){10}}
\put(20,30){\line(0,1){20}}
\put(40,20){\line(0,1){30}}
\put(60,20){\line(0,1){30}}
\put(60,10){\line(0,1){40}}
\put(80,10){\line(0,1){40}}
\put(100,10){\line(0,1){40}}
\put(120,10){\line(0,1){40}}
\put(140,10){\line(0,1){40}}
\put(60,10){\line(1,0){85}}
\put(2,43){$x_1\!\oplus\! x_2$}
\put(22,43){$x_1\!\oplus\! x_3$}
\put(42,43){$x_1\!\oplus\! x_4$}
\put(22,33){$x_2\!\oplus\! x_3$}
\put(42,33){$x_2\!\oplus\! x_4$}
\put(42,23){$x_3\!\oplus\! x_4$}
\put(62,43){$x_1\!\oplus\! b_1$}
\put(82,43){$x_1\!\oplus\! b_2$}
\put(102,43){$x_1\!\oplus\! b_3$}
\put(62,33){$x_2\!\oplus\! b_1$}
\put(82,33){$x_2\!\oplus\! b_2$}
\put(102,33){$x_2\!\oplus\! b_3$}
\put(62,23){$x_3\!\oplus\! b_1$}
\put(82,23){$x_3\!\oplus\! b_2$}
\put(102,23){$x_3\!\oplus\! b_3$}
\put(122,43){$x_1\!\oplus\! b_4$}
\put(122,33){$x_2\!\oplus\! b_4$}
\put(122,23){$x_3\!\oplus\! b_4$}
\put(62,13){$x_4\!\oplus\! b_1$}
\put(82,13){$x_4\!\oplus\! b_2$}
\put(102,13){$x_4\!\oplus\! b_3$}
\put(122,13){$x_4\!\oplus\! b_4$}
\put(143,13){$\cdots\;.$}
\put(143,23){$\cdots$}
\put(143,33){$\cdots$}
\put(143,43){$\cdots$}
\end{picture}
\end{example}
We denote 
$\mathcal{E}_n^B(\lambda)=\mathcal{E}_n^C(\lambda)=\mathcal{E}_n^I(\lambda)$ and 
$\mathcal{E}_n^D(\lambda)=\mathcal{E}_n^{I\!I}(\lambda).$
We also denote $B^B(D)=B^C(D)=B(D)^I(\lambda)$ and 
$B^D(D)=B(D)^{I\!I}.$
Define
$$
wt^X(D)=\prod_{(i,j)\in D}wt^X(i,j)\prod_{(i',j')\in B^X(D)}
(1+\beta wt^X(i',j')).
$$
Recall the definition (\ref{eq:defGXfinite}) of $GX^{(n)}(x|b)$.
\begin{thm}\label{thm:EYD}
Let $\lambda$ be a strict partition of length $\leq n.$
Then we have
\begin{equation}
GX_\lambda^{(n)}(x|b)=
\sum_{D\in \mathcal{E}_n^X(\lambda)}
wt^X(D).
\label{eq:EYDQ}
\end{equation}
\end{thm}

Consider the right diagram $D$ 
in Example \ref{ex:EYD}.
The corresponding summands of the right-hand side 
of (\ref{eq:EYDQ}) are equal to the following:
\begin{eqnarray*}
wt^B(D)&=&x_1(x_1\oplus x_2)(x_2\oplus x_4)
(1+\beta (x_1\oplus x_3))(x_3\oplus b_2)(1+\beta (x_2\oplus b_1))\\
&&
x_4(1+\beta (x_2\oplus x_2))
(1+\beta (x_3\oplus x_3))
(x_4\oplus b_1),\\
wt^C(D)&=&(x_1\oplus x_1)(x_1\oplus x_2)(x_2\oplus x_4)
(1+\beta (x_1\oplus x_3))(x_3\oplus b_2)(1+\beta (x_2\oplus b_1))\\
&&
(x_4\oplus x_4)(1+\beta (x_2\oplus x_2))
(1+\beta (x_3\oplus x_3))
(x_4\oplus b_1),\\
wt^D(D)&=&(x_1\oplus x_2)(x_1\oplus x_3)
(x_2\oplus b_1)(1+\beta (x_1\oplus x_4))
(x_3\oplus b_3)(1+\beta (x_2\oplus b_2))
\\
&&(x_4\oplus b_1)(1+\beta (x_2\oplus x_3))
(x_4\oplus b_2).
\end{eqnarray*}
Note that $B_{I\!I}(D)$ is given in Example \ref{ex:EYDs}.

\begin{cor}\label{cor:positive}
The coefficients of $GX_\lambda^{(n)}(x|b)$ 
as polynomials in $x_1,\ldots,x_n$ are 
in $\mathbb{N}[\beta,b_1,b_2,\ldots].$
\end{cor}

\subsection{Type A case}
There is an analogous formula for $G_\lambda(x_1,\ldots,x_n|b).$
Define $\mathcal{D}_n^A=
\{(i,j)\in \Z^2\;|\;1\leq i\leq n,\;j\geq 0\}.$
Define for $(i,j)\in \mathcal{D}_n^A$, the following weight of type $A:$
$$
wt^A(i,j)=x_i\oplus b_j.
$$

For a partition $\lambda\in \mathcal{P}_n$ 
the set of excited Young diagrams 
$\mathcal{E}_n^A(\la)$ is defined in a similar way
(see \cite{EYD}). For $D\in \mathcal{E}_n^A(\la)$ the set 
$B^A(D)$ is defined in the same way as in \S \ref{ssec:EYD}.
\begin{prop}[EYD formula for Type A]\label{prop:A}
Let $\lambda $ be a partition in $\mathcal{P}_n.$ Then 
$$G_\lambda(x_1,\ldots,x_n|b)
=\sum_{D\in \mathcal{E}_n^A(\lambda)}
\prod_{(i,j)\in D}wt_n^A(i,j)\prod_{(i',j')\in B^A(D)}(1+\beta wt_n^A(i',j')).
$$
\end{prop}

The proof of this proposition can be given in the same line as for Thm. \ref{thm:EYD}.
The arguments are left to the reader because 
these are quite similar and 
easier than types $B,C,$ and $D.$
Furthermore, we can show $G_\lambda(x_1,\ldots,x_n|b)$ is 
equal to a combinatorial expression in terms of set-valued tableaux
given by McNamara (cf. {\it Remark} \ref{rmk:McNa}).

\setcounter{equation}{0}
\section{Proofs of Thm.\:\ref{thm:SVT} and Thm. \ref{thm:EYD}}\label{sec:comb_proof}
In this section, we give proofs of Thm.\:\ref{thm:SVT} and Thm.\:\ref{thm:EYD}.
To make them more comprehensible, we consider 
type $C$. The modifications for other types are left to the reader.

\subsection{Locally equivalent weight functions}
Let $S$ be a subset of $\mathcal{D}_n.$ 
By a {\it weight function} on $S$ we mean a map $w$ from $S$ to a commutative ring $R.$
Suppose a weight function $w$ on $\mathcal{D}_n$ is given.
For a strict partition $\lambda\in \mathcal{SP}_n$, we denote 
$$
E_\lambda(w)=\sum_{D\in \mathcal{E}_n(\lambda)}w(D),\quad
w(D):=\prod_{c\in D}w(c)\prod_{c'\in B(D)}(1+\beta w(c')).
$$
Two weight functions $w$ and $w'$ on $\mathcal{D}_n$ are {\it equivalent\/} if 
the equality $E_\lambda(w)=E_\lambda(w')$ holds
for all $\lambda\in \mathcal{SP}_n.$
Two weight functions $w_1$ and $w_2$ on $S$ are {\it locally 
equivalent\/} on $S,$ 
 if for any weight function $w$ on $\mathcal{D}_n\setminus S,$
the following two weight functions $w\vee w_1$ and $w\vee w_2$ on $\mathcal{D}_n$
are equivalent: 
$$
(w\vee w_i)(\alpha)=\begin{cases}
w(\alpha) & \mbox{if}\;\alpha\in\mathcal{D}_n\setminus S\\
w_i(\alpha) & \mbox{if}\;\alpha\in  S
\end{cases}\quad (i=1,2).
$$

\setlength{\unitlength}{0.7mm}

\begin{example} 
We will show that the following two 
weight functions on $\mathcal{D}_2$
are equivalent for arbitrary $x,y,z$, etc. 
In other words, the two weight functions 
on $S$ 
are locally equivalent, where 
$S$ is the set of the gray boxes.
 
\begin{picture}(200,25)
\multiput(0,0)(90,0){2}{
\put(10,20){\line(1,0){70}}
\put(10,10){\line(1,0){70}}
\put(25,0){\line(1,0){55}}
\put(10,10){\line(0,1){10}}
\put(25,0){\line(0,1){20}}
\put(40,0){\line(0,1){20}}
\put(55,0){\line(0,1){20}}
\put(70,0){\line(0,1){20}}}
\put(131,1.5){$\Cell$}
\put(131,11.5){$\Cell$}
\put(116,11.5){$\Cell$}
\put(116,1.5){$\Cell$}
\put(41,1.5){$\Cell$}
\put(41,11.5){$\Cell$}
\put(26,11.5){$\Cell$}
\put(26,1.5){$\Cell$}
\put(15,13){$x$}
\put(25.8,13){$a\oplus t$}
\put(45,13){$a$}
\put(60,13){$y$}
\put(72,13){$\cdots$}
\put(60,3){$z$}
\put(72,3){$\cdots$}
\put(45,3){$0$}
\put(30,3){$0$}
\put(150,13){$y$}
\put(162,13){$\cdots$}
\put(150,3){$z$}
\put(162,3){$\cdots$}
\put(105,13){$x$}
\put(131,13){$a\oplus t$}
\put(121,13){$a$}
\put(120,3){$0$}
\put(137,3){$ t$}
\end{picture}

\noindent Let $w,w'$ denote the above weight functions.
If $\lambda=$ \begin{tiny}$\Tableau{&&}$\end{tiny} then the corresponding 
EYD sums are 
\begin{eqnarray*}
E_\lambda(w)&=&x(a\oplus t)a+x(a\oplus t)z(1+\beta a)\\
E_\lambda(w')&=&xa(a\oplus t)+xaz(1+\beta (a\oplus t))+xt(1+\beta a)z
\end{eqnarray*}
We can easily check the equality $E_\lambda(w)=E_\lambda(w').$
One sees that the only remaining $\lambda$ which we have to check is $\lambda=$ \begin{tiny}$\Tableau{&}$\end{tiny} .
Then we have the following two EYD sums which are equal:
\begin{equation*}
E_\lambda(w)=x(a\oplus t),\quad
E_\lambda(w')=xa+xt(1+\beta a).
\end{equation*}

\end{example}

\setlength{\unitlength}{0.6mm}

Next lemma is a generalization of this example.

\begin{lem}\label{lem:EYDsumA}   The
following weight functions on $S$
are locally equivalent :
\begin{center}
\begin{picture}(200,50)
\put(83,0){\multiput(2,0)(0,10){6}{\line(1,0){54}}
\multiput(2,0)(27,0){3}{\line(0,1){50}}
\put(33,3){$a_0\oplus t$}
\put(33,13){$a_1\oplus t$}
\put(33,23){$a_2\oplus t$}
\put(41,32){$\vdots$}
\put(33,43){$a_{n}\oplus t$}
\put(10,3){$0$}
\put(10,13){$a_1$}
\put(10,23){$a_2$}
\put(13,33){$\vdots$}
\put(10,43){$a_{n}$}
}
\put(0,0){\multiput(2,0)(0,10){6}{\line(1,0){54}}
\multiput(2,0)(27,0){3}{\line(0,1){50}}
\put(37,3){$a_0$}
\put(37,13){$a_1$}
\put(37,23){$a_2$}
\put(40,32){$\vdots$}
\put(37,43){$a_{n}$}
\put(10,3){$0$}
\put(7,13){$a_1\oplus t$}
\put(7,23){$a_2\oplus t$}
\put(15,33){$\vdots$}
\put(7,43){$a_{n}\!\oplus t$}
}
\end{picture}
\end{center}
where $S$ is 
the subset of $\mathcal{D}_{n+1}$ consisting of 
two columns ($i$-th and $(i+1)$-st) 
with $i\geq n+1.$
\end{lem}

\begin{lem}\label{lem:EYDsumC} The
following weight functions on $S$
are locally equivalent :

\begin{picture}(120,60)
\multiput(0,30)(0,10){3}{\line(1,0){60}}
\multiput(0,10)(0,10){2}{\line(1,0){155}}
\multiput(30,0)(30,0){2}{\line(0,1){50}}
\put(0,10){\line(0,1){40}}
\put(30,0){\line(1,0){125}}
\multiput(90,0)(30,0){3}{\line(0,1){20}}
\put(33,3){$c$}
\put(33,13){$c\oplus t$}
\put(35,23){$b_1$}
\put(38,32){$\vdots$}
\put(35,43){$b_{i-1}$}
%\put(35,53){$0$}
%
\put(2,13){$c\oplus t\oplus t$}
\put(3,23){$b_1\oplus t$}
\put(10,33){$\vdots$}
\put(3,43){$b_{i-1}\oplus t$}
%\put(2,53){$b_{h+1}\oplus t$}
%
\put(63,3){$a_1$}
\put(93,3){$a_2$}
\put(123,3){$a_3$}
\put(63,13){$a_1\oplus t$}
\put(93,13){$a_2\oplus t$}
\put(123,13){$a_3\oplus t$}
\end{picture}

\begin{picture}(120,70)
\multiput(0,30)(0,10){3}{\line(1,0){60}}
\multiput(0,10)(0,10){2}{\line(1,0){155}}
\multiput(30,0)(30,0){2}{\line(0,1){50}}
\put(0,10){\line(0,1){40}}
\put(30,0){\line(1,0){125}}
\multiput(90,0)(30,0){3}{\line(0,1){20}}
\put(33,3){$c\oplus t\oplus t$}
\put(33,13){$c\oplus t$}
\put(35,23){$b_1\oplus t$}
\put(38,32){$\vdots$}
\put(35,43){$b_{i-1}\oplus t$}
\put(2,13){$c$}
\put(3,23){$b_1$}
\put(10,33){$\vdots$}
\put(3,43){$b_{i-1}$}
\put(63,3){$a_1\oplus t$}
\put(93,3){$a_2\oplus t$}
\put(123,3){$a_3\oplus t$}
\put(63,13){$a_1$}
\put(93,13){$a_2$}
\put(123,13){$a_3$}
\end{picture}

\noindent where $S$ is 
the subset of $\mathcal{D}_{n}$ consisting of 
two columns and rows ($i$-th and $(i+1)$-st) 
with $i<n.$
\end{lem}

 Lemma \ref{lem:EYDsumA} and Lemma \ref{lem:EYDsumC} can be proved by 
elementary but quite tedious calculations. 
An alternative 
proof using the Yang--Baxter equation is available 
(the details will appear in \cite{N}). 

\subsection{Outline of the proof of Thm. \ref{thm:EYD}}
For $\lambda\in \mathcal{SP}_n$,
define   
$$E_\lambda(x_1,\ldots,x_n|b)=\sum_{D\in \mathcal{E}_n(\lambda)}wt(D).$$
\begin{prop}\label{prop:sym}
$E_\lambda(x_1,\ldots,x_n|b)
$ is in $G\Gamma_{n}.$
\end{prop}
{\it Proof.} Let $i$ be $1\leq i\leq n-1.$
We apply Lemma \ref{lem:EYDsumC}
to $i$-th and $(i+1)$-st columns and rows 
with $t=x_i\ominus x_{i+1}.$ It turns out that  $E_\lambda(x_1,\ldots,x_{n}|b)$
is invariant under the exchange of $x_i$ and $x_{i+1}.$
Hence $E_\lambda(x_1,\ldots,x_{n}|b)$ is symmetric.

Next we show the cancellation property of Def. \ref{def:Ksuper}.
Let $x_1=t,x_2=\ominus t.$
Then weight function is given as follows:

 \setlength{\unitlength}{0.7mm}
 \begin{picture}(300,45)
\put(0,40){\line(1,0){148}}
\put(0,30){\line(1,0){148}}
\put(20,20){\line(1,0){128}}
\put(40,10){\line(1,0){108}}
\put(0,30){\line(0,1){10}}
\put(20,20){\line(0,1){20}}
\put(40,10){\line(0,1){30}}
\put(60,10){\line(0,1){30}}
\put(60,0){\line(0,1){40}}
\put(80,0){\line(0,1){40}}
\put(100,0){\line(0,1){40}}
\put(120,0){\line(0,1){40}}
\put(140,0){\line(0,1){40}}
\put(60,0){\line(1,0){85}}
\put(5,33){$t\!\oplus\! t$}
\put(25.5,33){$t\!\ominus{t}$}
\put(45.5,33){$t\!\oplus\! x_3$}
\put(21,23){$\ominus{t}\!\ominus{t}$}
\put(41,23){$\ominus{t}\!\oplus\! x_3$}
\put(42,13){$x_3\!\oplus\! x_3$}
\put(65.5,33){$t\!\oplus\! x_4$}
\put(85.5,33){$t\!\oplus\! b_1$}
\put(105.5,33){$t\!\oplus\! b_2$}
\put(61,23){$\ominus{t}\!\oplus\! x_4$}
\put(81,23){$\ominus{t}\!\oplus\! b_1$}
\put(101,23){$\ominus{t}\!\oplus\! b_2$}
\put(62,13){$x_3\!\oplus\! x_4$}
\put(82,13){$x_3\!\oplus\! b_1$}
\put(102,13){$x_3\!\oplus\! b_2$}
\put(125.5,33){$t\!\oplus\! b_3$}
\put(121,23){$\ominus{t}\!\oplus\! b_3$}
\put(122,13){$x_3\!\oplus\! b_3$}
\put(62,3){$x_4\!\oplus\! x_4$}
\put(82,3){$x_4\!\oplus\! b_1$}
\put(102,3){$x_4\!\oplus\! b_2$}
\put(122,3){$x_4\!\oplus\! b_3$}
\put(143,3){$\cdots\;.$}
\put(143,13){$\cdots$}
\put(143,23){$\cdots$}
\put(143,33){$\cdots$}
\end{picture}

\noindent This is obviously equivalent to 

 \setlength{\unitlength}{0.7mm}
 \begin{picture}(300,45)
\put(0,40){\line(1,0){148}}
\put(0,30){\line(1,0){148}}
\put(20,20){\line(1,0){128}}
\put(40,10){\line(1,0){108}}
\put(0,30){\line(0,1){10}}
\put(20,20){\line(0,1){20}}
\put(40,10){\line(0,1){30}}
\put(60,10){\line(0,1){30}}
\put(60,0){\line(0,1){40}}
\put(80,0){\line(0,1){40}}
\put(100,0){\line(0,1){40}}
\put(120,0){\line(0,1){40}}
\put(140,0){\line(0,1){40}}
\put(60,0){\line(1,0){85}}
\put(8,33){$0$}
\put(28,33){$0$}
\put(45.5,33){$t\!\oplus\! x_3$}
\put(28,23){$0$}
\put(41,23){$\ominus{t}\!\oplus\! x_3$}
\put(42,13){$x_3\!\oplus\! x_3$}
\put(65.5,33){$t\!\oplus\! x_4$}
\put(85.5,33){$t\!\oplus\! b_1$}
\put(105.5,33){$t\!\oplus\! b_2$}
\put(61,23){$\ominus{t}\!\oplus\! x_4$}
\put(81,23){$\ominus{t}\!\oplus\! b_1$}
\put(101,23){$\ominus{t}\!\oplus\! b_2$}
\put(62,13){$x_3\!\oplus\! x_4$}
\put(82,13){$x_3\!\oplus\! b_1$}
\put(102,13){$x_3\!\oplus\! b_2$}
\put(125.5,33){$t\!\oplus\! b_3$}
\put(121,23){$\ominus{t}\!\oplus\! b_3$}
\put(122,13){$x_3\!\oplus\! b_3$}
\put(62,3){$x_4\!\oplus\! x_4$}
\put(82,3){$x_4\!\oplus\! b_1$}
\put(102,3){$x_4\!\oplus\! b_2$}
\put(122,3){$x_4\!\oplus\! b_3$}
\put(143,3){$\cdots\;.$}
\put(143,13){$\cdots$}
\put(143,23){$\cdots$}
\put(143,33){$\cdots$}
\end{picture}
Then the sum of weights over $\mathcal{D}_n$
reduces to the sum over 
the subset of $\mathcal{D}_n$
consisting of the boxes $(i,j)\in\mathcal{D}_n$ such that $i\geq 3.$  
Since this subset is identified with $\mathcal{D}_{n-2}$
of  
type $C_{n-2}$ weights with respect to 
$x_3,\ldots,x_n,$
and hence the sum is equal to
$E_\la(0,0,x_3,\ldots,x_n|b).$
This implies $E_\la(x_1,\ldots,x_n|b)$ is in $G\Gamma_n$.
$\qed$

\begin{prop}[Stability]\label{prop:EYD_stable}
We have $E_\lambda(x_1,\ldots,x_n,0|b)
=E_\lambda(x_1,\ldots,x_n|b).$
\end{prop}
{\it Proof.} 
To make the stability property clearer, 
it is useful to arrange the weight in the following way
by using symmetry. 

\setlength{\unitlength}{0.7mm}
\begin{picture}(300,35)
\put(0,30){\line(1,0){145}}
\put(0,20){\line(1,0){145}}
\put(20,10){\line(1,0){125}}
\put(40,0){\line(1,0){105}}
\put(0,20){\line(0,1){10}}
\put(20,10){\line(0,1){20}}
\put(40,0){\line(0,1){30}}
\put(60,0){\line(0,1){30}}
\put(60,0){\line(0,1){30}}
\put(80,0){\line(0,1){30}}
\put(100,0){\line(0,1){30}}
\put(120,0){\line(0,1){30}}
\put(140,0){\line(0,1){30}}
\put(2,23){$x_3\!\oplus\! x_3$}
\put(22,23){$x_3\!\oplus\! x_2$}
\put(42,23){$x_3\!\oplus\! x_1$}
\put(22,13){$x_2\!\oplus\! x_2$}
\put(42,13){$x_2\!\oplus\! x_1$}
\put(42,3){$x_1\!\oplus\! x_1$}
\put(62,23){$x_3\!\oplus\! b_1$}
\put(82,23){$x_3\!\oplus\! b_2$}
\put(102,23){$x_3\!\oplus\! b_3$}
\put(62,13){$x_2\!\oplus\! b_1$}
\put(82,13){$x_2\!\oplus\! b_2$}
\put(102,13){$x_2\!\oplus\! b_3$}
\put(62,3){$x_1\!\oplus\! b_1$}
\put(82,3){$x_1\!\oplus\! b_2$}
\put(102,3){$x_1\!\oplus\! b_3$}
\put(122,23){$x_3\!\oplus\! b_4$}
\put(122,13){$x_2\!\oplus\! b_4$}
\put(122,3){$x_1\!\oplus\! b_4$}
\end{picture}
\noindent Set
$x_{n+1}=0.$ Then the weight at the box $(1,1)$ is zero.
So if $(1,1)\in D$ for $D\in \mathcal{E}_{n+1}(\lambda)$, then 
$wt(D)=0.$
The set consists of $D\in \mathcal{E}_{n+1}(\lambda)$ such that $(1,1)\notin D$ is 
naturally in bijection with $\mathcal{E}_{n}(\lambda).$ 
This bijection yields the equation.
$\qed$

\bigskip

Let us denote by $E_\lambda(x|b)$ the inverse 
limit $\varprojlim E_\lambda^{(n)}(x_1,\ldots,x_n|b)$ (cf. Prop. \ref{prop:EYD_stable}).
By Prop. \ref{def:Ksuper} we have the following.
\begin{lem}\label{lem:Ess}
$E_\lambda(x|b)$ is in $G\Gamma_\R^C.$ 
\end{lem}

In order to prove Thm.\:\ref{thm:EYD}, 
we need the following lemma.
\begin{lem}\label{lem:EYD_div} We have
$\pi_i E_\lambda(x|b)=\begin{cases}
E_{\lambda^{(i)}}(x|b) &\mbox{if}\;s_i\lambda<\lambda\\
-\beta E_\lambda(x|b) &\mbox{if}\;s_i\lambda\geq \lambda.
\end{cases}$
\end{lem}

With this lemma at hand, the proof for Thm.\:\ref{thm:EYD} is
completed as follows.

\bigskip

{\it Proof of Thm.\:\ref{thm:EYD}.}
Lemma \ref{lem:Ess} enables us to send $E_\lambda(x|b)$
by $\Phi$ into $\Psi$. 
By the same argument as in the proof 
of Thm. \ref{thm:main}, Lemma \ref{lem:EYD_div} implies that 
$\{\Phi(E_\lambda(x|b))\}$ is a family of Schubert classes.
Since we know $\{\Phi(GQ_\lambda(x|b))\}$ is a family of Schubert classes 
(Thm. \ref{thm:main}), 
the injectivity of $\Phi$  
and the uniqueness of a family of Schubert classes 
implies $E_\lambda(x|b)=GQ_\lambda(x|b).$
$\qed$

\subsection{Proof of Lemma \ref{lem:EYD_div}}

Let  $i\geq 1.$ We calculate 
$s_i E_\lambda(x|b)$, which is equal to the weight 
sum over $\mathcal{E}_n(\lambda)$ of the weights 
obtained from the original one by replacing  
$(n+i)$-th and $(n+i+1)$-st columns with the following
left diagram:
\setlength{\unitlength}{0.6mm}
\begin{center}
\begin{picture}(150,55)
\multiput(0,0)(0,10){6}{\line(1,0){63}}
\multiput(2,0)(30,0){3}{\line(0,1){50}}
\put(49,3){$0$}
\put(37,13){$x_n\oplus b_i$}
\put(37,23){$\vdots$}
\put(37,32){$x_{2}\oplus b_i$}
\put(37,43){$x_1\oplus b_i$}
\put(13,3){$0$}
\put(4,13){$x_n\oplus b_{i+1}$}
\put(18,23){$\vdots$}
%\put(3,33){$x_{n-2}\oplus b_{i+1}$}
\put(4,32){$x_{2}\oplus b_{i+1}$}
\put(4,43){$x_1\oplus b_{i+1}$}
\put(73,30){$\rightarrow$}
\multiput(90,0)(0,10){6}{\line(1,0){63}}
\multiput(92,0)(30,0){3}{\line(0,1){50}}
\put(124,3){$b_{i+1}\ominus b_{i}$}
\put(125,13){$x_n\oplus b_{i+1}$}
\put(135,23){$\vdots$}
\put(125,32){$x_{2}\oplus b_{i+1}$}
\put(125,43){$x_{1}\oplus b_{i+1}$}
\put(100,3){$0$}
\put(95,13){$x_n\oplus b_{i}$}
\put(100,23){$\vdots$}
\put(95,32){$x_{2}\oplus b_{i}$}
\put(95,43){$x_1\oplus b_{i}$}
\end{picture}
\end{center}
We consider each $D\in \mathcal{E}_n(\lambda)$ 
as a subset in $\mathcal{D}_{n+1},$
on which the values of weight on $(n+1)$-st row are all zero.
Applying Lemma \ref{lem:EYDsumA}, 
we can change the weight  
with the right one without changing the weight sum in total.
There is an additional box at the bottom of $(i+n)$-th column,
 whose weight is $e(\alpha_i)=b_{i+1}\ominus b_i.$
 For example if $n=3,i=1$ then the two weight functions 
 look as follows:

 \setlength{\unitlength}{0.7mm}
 \begin{picture}(300,55)
\put(0,50){\line(1,0){145}}
\put(0,40){\line(1,0){145}}
\put(20,30){\line(1,0){125}}
\put(40,20){\line(1,0){105}}
\put(0,40){\line(0,1){10}}
\put(20,30){\line(0,1){20}}
\put(40,20){\line(0,1){30}}
\put(60,20){\line(0,1){30}}
\put(60,10){\line(0,1){40}}
\put(80,10){\line(0,1){40}}
\put(100,10){\line(0,1){40}}
\put(120,10){\line(0,1){40}}
\put(140,10){\line(0,1){40}}
\put(60,10){\line(1,0){85}}
\put(2,43){$x_3\!\oplus\! x_3$}
\put(22,43){$x_3\!\oplus\! x_2$}
\put(42,43){$x_3\!\oplus\! x_1$}
\put(22,33){$x_2\!\oplus\! x_2$}
\put(42,33){$x_2\!\oplus\! x_1$}
\put(42,23){$x_1\!\oplus\! x_1$}
\put(62,43){$x_3\!\oplus\! b_2$}
\put(82,43){$x_3\!\oplus\! b_1$}
\put(102,43){$x_3\!\oplus\! b_3$}
\put(62,33){$x_2\!\oplus\! b_2$}
\put(82,33){$x_2\!\oplus\! b_1$}
\put(102,33){$x_2\!\oplus\! b_3$}
\put(62,23){$x_1\!\oplus\! b_2$}
\put(82,23){$x_1\!\oplus\! b_1$}
\put(102,23){$x_1\!\oplus\! b_3$}
\put(122,43){$x_3\!\oplus\! b_4$}
\put(122,33){$x_2\!\oplus\! b_4$}
\put(122,23){$x_1\!\oplus\! b_4$}
\put(88,13){$0$}
\put(68,13){$0$}
\put(108,13){$0$}
\put(128,13){$0$}
\end{picture}

 \setlength{\unitlength}{0.7mm}
 \begin{picture}(300,55)
\put(0,50){\line(1,0){145}}
\put(0,40){\line(1,0){145}}
\put(20,30){\line(1,0){125}}
\put(40,20){\line(1,0){105}}
\put(0,40){\line(0,1){10}}
\put(20,30){\line(0,1){20}}
\put(40,20){\line(0,1){30}}
\put(60,20){\line(0,1){30}}
\put(60,10){\line(0,1){40}}
\put(80,10){\line(0,1){40}}
\put(100,10){\line(0,1){40}}
\put(120,10){\line(0,1){40}}
\put(140,10){\line(0,1){40}}
\put(60,10){\line(1,0){85}}
\put(2,43){$x_3\!\oplus\! x_3$}
\put(22,43){$x_3\!\oplus\! x_2$}
\put(42,43){$x_3\!\oplus\! x_1$}
\put(22,33){$x_2\!\oplus\! x_2$}
\put(42,33){$x_2\!\oplus\! x_1$}
\put(42,23){$x_1\!\oplus\! x_1$}
\put(62,43){$x_3\!\oplus\! b_1$}
\put(82,43){$x_3\!\oplus\! b_2$}
\put(102,43){$x_3\!\oplus\! b_3$}
\put(62,33){$x_2\!\oplus\! b_1$}
\put(82,33){$x_2\!\oplus\! b_2$}
\put(102,33){$x_2\!\oplus\! b_3$}
\put(62,23){$x_1\!\oplus\! b_1$}
\put(82,23){$x_1\!\oplus\! b_2$}
\put(102,23){$x_1\!\oplus\! b_3$}
\put(122,43){$x_3\!\oplus\! b_4$}
\put(122,33){$x_2\!\oplus\! b_4$}
\put(122,23){$x_1\!\oplus\! b_4$}
\put(82,13){$b_2\!\ominus\! b_1$}
\put(68,13){$0$}
\put(108,13){$0$}
\put(128,13){$0$}
\end{picture}
We denote the first one by $s_i(wt)$ (the function obtained from $wt$ by 
exchanging $b_i$ and $b_{i+1}$) and the second one by $wt_{+i}.$
Thus we have the following.
\begin{lem} 
For $i\geq 1$ we have 
$s_i E_\lambda(x_1,\ldots,x_n|b)=\sum_{D\in {\mathcal{E}}_{n+1}(\lambda)}{wt_{+i}}(D).$
\end{lem}
{\it Proof of Lemma \ref{lem:EYD_div} for $i\geq 1.$} 
Suppose $\lambda$ is $i$-removable.
We denote by ${\mathcal{E}}_n^{+i}(\lambda)$
the set of $D\in {\mathcal{E}}_{n+1}(\lambda)$ such that $wt_{+i}(D)\neq 0.$
We have a 
map $r_i: {\mathcal{E}}_n^{+i}(\lambda)
\rightarrow \mathcal{E}_n(\lambda^{(i)})$ defined by 
removing the box  of content $i$
having the largest row index. 
One readily sees that 
$r_i$ is surjective.
Fix $D_0\in \mathcal{E}_n(\lambda^{(i)}).$
It suffices to show the following equation:
\begin{equation}
\sum_{{D}\in {\mathcal{E}^{+i}_n}(\lambda),\;r_i({D})=D_0}
{wt}_{+i}(D)-(1+\beta e(\alpha_i))
\sum_{{D}\in {\mathcal{E}_n}(\lambda),\;r_i(D)=D_0}wt(D)
=e(\alpha_i)wt(D_0)\label{eq:DivEYD}
\end{equation}
Note that each term on the left-hand side of \ref{eq:DivEYD}
is clearly a multiple of $wt(D_0).$
The differences come from
 the contributions of the box to be removed by $r_i$. An 
 obvious equation 
\begin{equation}
\sum_{i=1}^r a_i\prod_{j=1}^{i-1}(1+\beta a_j)
-(1+\beta a_r)\sum_{i=1}^{r-1}a_i\prod_{j=1}^{i-1}(1+\beta a_j)=a_r.
\end{equation}
  depicted below gives (\ref{eq:DivEYD}).

%\vspace{1cm}
\setlength{\unitlength}{0.17mm}
\begin{picture}(350,65)
\put(20,41){$\scell$}
\put(90,31){$\scell$}
\put(160,21){$\scell$}
\put(230,11){$\scell$}
\put(360,41){$\scell$}
\put(419,41){\begin{tiny}$\times$\end{tiny}}
\put(430,31){$\scell$}
\put(489,31){\begin{tiny}$\times$\end{tiny}}
\put(500,21){$\scell$}
\put(79,41){\begin{tiny}$\times$\end{tiny}}
\put(139,41){\begin{tiny}$\times$\end{tiny}}
\put(199,41){\begin{tiny}$\times$\end{tiny}}
\put(149,31){\begin{tiny}$\times$\end{tiny}}
\put(209,31){\begin{tiny}$\times$\end{tiny}}
\put(219,21){\begin{tiny}$\times$\end{tiny}}
\put(479,41){\begin{tiny}$\times$\end{tiny}}
\put(319,11){\begin{tiny}$\times$\end{tiny}}
%\put(530,11){$\scell$}
%
\put(0,25){$\Bigl($}
\put(245,25){$\Bigr)$}
\put(265,25){$-$}
\put(60,25){$+$}
\put(120,25){$+$}
\put(180,25){$+$}
\put(400,25){$+$}
\put(460,25){$+$}
\put(340,25){$\Bigl($}
\put(530,25){$\Bigr)$}
\multiput(10,0)(60,0){4}{
\multiput(10,40)(10,-10){4}
{\put(0,0){\line(1,0){10}}
\put(0,0){\line(0,1){10}}
\put(10,0){\line(0,1){10}}
\put(0,10){\line(1,0){10}}}}
\multiput(360,0)(60,0){3}{
\multiput(0,40)(10,-10){4}
{\put(0,0){\line(1,0){10}}
\put(0,0){\line(0,1){10}}
\put(10,0){\line(0,1){10}}
\put(0,10){\line(1,0){10}}}}
\multiput(280,0)(60,0){1}{
\multiput(10,40)(10,-10){4}
{\put(0,0){\line(1,0){10}}
\put(0,0){\line(0,1){10}}
\put(10,0){\line(0,1){10}}
\put(0,10){\line(1,0){10}}}}
\put(550,25){$=$}
\put(610,11){$\scell$}
\multiput(570,0)(60,0){1}{
\multiput(10,40)(10,-10){4}
{\put(0,0){\line(1,0){10}}
\put(0,0){\line(0,1){10}}
\put(10,0){\line(0,1){10}}
\put(0,10){\line(1,0){10}}}}
\end{picture}

If $\lambda$ is not $i$-removable, 
one easily sees that $s_i E_\lambda(x_1,\ldots,x_n|b)
=E_\lambda(x_1,\ldots,x_n|b).$
Then we have $\pi_i E_\lambda(x|b)
=-\beta E_\lambda(x|b).$
$\qed$
\setlength{\unitlength}{0.7mm}

\bigskip

In order to prove Lemma \ref{lem:EYD_div} for $i=0$,
we need the following lemma.
Let 
$wt_{+0}$ denote the weight function on $\mathcal{D}_{n+1}$ defined 
similarly as $wt_{+i}\;(i\geq 1).$

\begin{lem}%\label{lem:EYD} 
We have 
$E_\lambda(x_1,\ldots,x_n,b_1|\ominus b_1,b_2,b_3,\ldots)=\sum_{D\in 
{\mathcal{E}}_{n+1}(\lambda)}{wt_{+0}}(D).$
\end{lem}

If $n=2$ then $E_\lambda(x_1,x_2,b_1|\ominus b_1,b_2,b_3,\ldots)$ is given by the following table.

\setlength{\unitlength}{0.7mm}
\begin{picture}(300,40)
\put(68,31){$*$}
\put(88,31){$*$}
\put(0,30){\line(1,0){145}}
\put(0,20){\line(1,0){145}}
\put(20,10){\line(1,0){125}}
\put(40,0){\line(1,0){105}}
\put(0,20){\line(0,1){10}}
\put(20,10){\line(0,1){20}}
\put(40,0){\line(0,1){30}}
\put(60,0){\line(0,1){30}}
\put(60,0){\line(0,1){30}}
\put(80,0){\line(0,1){30}}
\put(100,0){\line(0,1){30}}
\put(120,0){\line(0,1){30}}
\put(140,0){\line(0,1){30}}
\put(2,23){$x_2\!\oplus\! x_2$}
\put(22,23){$x_2\!\oplus\! x_1$}
\put(42,23){$x_2\!\oplus\! b_1$}
\put(22,13){$x_1\!\oplus\! x_2$}
\put(42,13){$x_1\!\oplus\! b_1$}
\put(42,3){$b_1\!\oplus\! b_1$}
\put(62,23){$x_2\!\ominus\! b_1$}
\put(82,23){$x_2\!\oplus\! b_2$}
\put(102,23){$x_2\!\oplus\! b_3$}
\put(62,13){$x_1\!\ominus\! b_1$}
\put(82,13){$x_1\!\oplus\! b_2$}
\put(102,13){$x_1\!\oplus\! b_3$}
\put(62,3){$b_1\!\ominus\! b_1$}
\put(82,3){$b_1\!\oplus\! b_2$}
\put(102,3){$b_1\!\oplus\! b_3$}
\put(122,23){$x_2\!\oplus\! b_4$}
\put(122,13){$x_1\!\oplus\! b_4$}
\put(122,3){$b_1\!\oplus\! b_4$}
\end{picture}

Note that $b_1\ominus b_1=0.$
By applying Lemma \ref{lem:EYDsumA}
to two columns indicated by $*$, we can deform 
this as follows:

\begin{picture}(300,40)
\put(0,30){\line(1,0){145}}
\put(0,20){\line(1,0){145}}
\put(20,10){\line(1,0){125}}
\put(40,0){\line(1,0){105}}
\put(0,20){\line(0,1){10}}
\put(20,10){\line(0,1){20}}
\put(40,0){\line(0,1){30}}
\put(60,0){\line(0,1){30}}
\put(60,0){\line(0,1){30}}
\put(80,0){\line(0,1){30}}
\put(100,0){\line(0,1){30}}
\put(120,0){\line(0,1){30}}
\put(140,0){\line(0,1){30}}
\put(2,23){$x_2\!\oplus\! x_2$}
\put(22,23){$x_2\!\oplus\! x_1$}
\put(42,23){$x_2\!\oplus\! b_1$}
\put(22,13){$x_1\!\oplus\! x_1$}
\put(42,13){$x_1\!\oplus\! b_1$}
\put(42,3){$b_1\!\oplus\! b_1$}
\put(62,23){$x_2\!\oplus\! b_2$}
\put(82,23){$x_2\!\ominus\! b_1$}
\put(102,23){$x_2\!\oplus\! b_3$}
\put(62,13){$x_1\!\oplus\! b_2$}
\put(82,13){$x_1\!\ominus\! b_1$}
\put(102,13){$x_1\!\oplus\! b_3$}
\put(62,3){$0$}
\put(82,3){$0$}
\put(102,3){$b_1\!\oplus\! b_3$}
\put(122,23){$x_2\!\oplus\! b_4$}
\put(122,13){$x_1\!\oplus\! b_4$}
\put(122,3){$b_1\!\oplus\! b_4$}
\end{picture}

Repeat this process sufficiently many times 
 we have

\begin{picture}(300,35)
\put(0,30){\line(1,0){145}}
\put(0,20){\line(1,0){145}}
\put(20,10){\line(1,0){125}}
\put(40,0){\line(1,0){105}}
\put(0,20){\line(0,1){10}}
\put(20,10){\line(0,1){20}}
\put(40,0){\line(0,1){30}}
\put(60,0){\line(0,1){30}}
\put(60,0){\line(0,1){30}}
\put(80,0){\line(0,1){30}}
\put(100,0){\line(0,1){30}}
\put(120,0){\line(0,1){30}}
\put(140,0){\line(0,1){30}}
\put(2,23){$x_2\!\oplus\! x_2$}
\put(22,23){$x_2\!\oplus\! x_1$}
\put(42,23){$x_2\!\oplus\! b_1$}
\put(22,13){$x_1\!\oplus\! x_1$}
\put(42,13){$x_1\!\oplus\! b_1$}
\put(42,3){$b_1\!\oplus\! b_1$}
\put(62,23){$x_2\!\oplus\! b_2$}
\put(82,23){$x_2\!\oplus\! b_3$}
\put(102,23){$x_2\!\oplus\! b_4$}
\put(62,13){$x_1\!\oplus\! b_2$}
\put(82,13){$x_1\!\oplus\! b_3$}
\put(102,13){$x_1\!\oplus\! b_4$}
\put(62,3){$0$}
\put(82,3){$0$}
\put(102,3){$0$}
\put(122,23){$x_2\!\oplus\! b_5$}
\put(122,13){$x_1\!\oplus\! b_5$}
\put(122,3){$0$}
\end{picture}

The proof of Lemma \ref{lem:EYD_div} for $i=0$
is 
completed.

\subsection{Excited Young diagrams and Set-valued Tableaux}
Let $E_\lambda(x|b)$ and $T_\lambda(x|b)$ be the functions 
on the right-hand sides of (\ref{eq:EYDQ}) and (\ref{eq:TabGQ}) respectively.
In this section, we prove the following.
\begin{prop}\label{prop:E=T}
We have $E_\lambda(x|b)=T_\lambda(x|b).$
\end{prop}

The
main idea of the proof is to use equivalence of weight functions
so that we have ``separation of variables''.
For example, let $n=2$, and consider the weight function $wt^C$ defined by (\ref{eq:weightC}).
One readily sees that the following weight function is equivalent to 
the original one.

\begin{small}
$$
\begin{array}{ccccccccccc}
 x_{1} & 0 & 0 & 0
 & 0&0 &  0 &   \cdots\\
 & x_1 & x_1\oplus x_2  & x_1\oplus b_1   & x_1\oplus b_2& x_1\oplus b_3&  
  x_1\oplus b_4& \cdots \\
%    & &    &  &  &  & 0 & x_{3'}\ominus b_2  & 0 & 0& \cdots\\
   &  & 0 &  0 & 0 & 0 & 0 &\cdots\\
     & & & x_2 & 0 & 0 &  0 & \cdots\\
            & & && x_{2} & x_{2}\oplus b_1 &  x_{2}\oplus b_2
                & \cdots
                \end{array}
$$
\end{small}

\noindent This can be further deformed to the following equivalent one:

\begin{small}
$$
\begin{array}{cccccccccccccccccc}
 x_{1} & 0 & 0 & 0
 & 0&0 &  0&   \cdots\\
 & x_1 & x_1\oplus b_1  & x_1\oplus x_2   & x_1\oplus b_2& x_1\oplus b_3
 &  x_1\oplus b_4&  \cdots \\
%    & &    &  &  &  & 0 & x_{3'}\ominus b_2  & 0 & 0& \cdots\\
   &  & 0 &  x_2\ominus b_1 & 0 & 0 & 0 & \cdots\\
     & & & x_2 & 0 & 0 &  0 & \cdots\\
 & & && x_{2} & x_{2}\oplus b_1 &  x_{2}\oplus b_2&\cdots
\end{array}
$$
\end{small}

\noindent 
We repeat this process.
For example, if $n=3$, 
the following weight function
is equivalent to the original one.
\begin{tiny}
$$
\begin{array}{cccccccccccccccccc}
 & &     & x_{1} & 0 & 0 & 0
 & 0&0 &  0 &0&  0&0 & 0 &0& \cdots\\
  & &    &  & x_1 & x_1\oplus b_1  & x_1\oplus b_2  & x_{2}\oplus x_{1} & x_1\oplus b_3& x_1\oplus b_4&  x_1\oplus b_5& x_1\oplus b_6& x_1\oplus b_7 &x_1\oplus b_8 & x_1\oplus b_9 &\cdots \\
     & &    &    &  & 0 &  0 & x_{2}\ominus b_2 & 0 & 0&  0&0&0 &0&0 &\cdots\\
        & &    &     & & & 0 & x_{2}\ominus b_1 & 0 &  0 & 0&0&0 &0&0 &\cdots\\
                         &    &  &&  & & &x_{2} &0&  0&  0& 0&0 &0 &0 & \cdots\\
                &     &  &  & & &&& x_{2} & x_{2}\oplus b_1 &  x_{2}\oplus b_2&x_{2}\oplus x_{3}&
                x_{2}\oplus b_3 &x_{2}\oplus b_4 & x_{2}\oplus b_5 &\cdots
                \\
                &     &  &  & & &&&& 0 & 0&x_{3}\ominus b_2&0&0 & 0& \cdots\\
                &     &  &&&  & &&&& 0 & x_{3}\ominus b_1&0&0 & 0&\cdots
\\
                &     &  & &&&&&&&& x_{3} & 0&0 & 0 & \cdots
\\
                &   &&&  &  &  & & &&&& x_{3} & x_3\oplus b_1&x_3\oplus b_2&\cdots
\end{array}
$$
\end{tiny}

\noindent If we only consider a strict partition $\lambda$ such that
$\lambda_1\leq 3$,
we can use the following weight of ``separated'' form:

\begin{tiny}
\begin{equation*}
\begin{array}{ccccc ccccc cccccc}
  &  & x_{1} & 0&0 &0&0&0&0&0&0&0&0&0&\cdots\\
  &&&{x_{1}} &x_{1}\oplus b_1&x_1\oplus b_2&0&0&0&0&0&0&0&0&\cdots\\
    &&& &0&0&x_{2}\ominus b_2&0&0&0&0&0&0&0&\cdots\\
    &&& &&0&{x_{2}\ominus b_1}&0&0&0&0&0&0&0&\cdots\\
        &&& &&&x_{2}&0&0&0&0&0&0&0&\cdots\\
        &&& &&&&x_2&{x_{2}}\oplus b_1&x_{2}\oplus b_2&0&0&0&0&\cdots\\
        &&& &&&&&0&0&x_{3}\ominus b_2&0&0&0&\cdots\\
                &&& &&&&&&0&x_{3}\ominus b_1&0&0&0&\cdots\\
                &&& &&&&&&&{x_{3}}&0&0&0&\cdots\\
                &&& &&&&&&&&x_3&x_3\oplus b_1&{x_{3}}\oplus b_2&\cdots
%                \\
%                &&& &&&&&&&&&0&0&\cdots\\
%                &&& &&&&&&&&&&0&\cdots\\
                  \end{array}
\end{equation*}
\end{tiny}

In general, we consider 
a weight function of the above separated form.
Precisely, for positive integer $k$,
we define the weight function $wt_{\mathrm{sep}}$ on 
$\mathcal{D}_{nk+n-k+1}$
as follows:
$wt_{\mathrm{sep}}(i,j)=0
$ unless $j-i\leq k$ and 
($i\equiv 2\;(\mathrm{mod}\; k+1)$ or 
 $j\equiv 1\;(\mathrm{mod} \;k+1)).$
On the $j$th column with $j=(k+1)s-k\;(1\leq s\leq n)$, 
 $wt_{\mathrm{sep}}(i,j)=x_s\ominus b_{j-i}\;(j-k+1\leq i\leq j).$ 
On the $i$th row with $i=(k+1)t-k+1\;(1\leq t\leq n)$, 
 $wt_{\mathrm{sep}}(i,j)=x_t\oplus b_{j-i}\;(i\leq j\leq i+k-1).$

Let $\lambda$ be a strict partition.
Suppose $k$ is large enough so that $k\geq \lambda_1,$
and let $wt_{\mathrm{sep}}$ be the  
weight
function defined by the preceding paragraph. 
Let $\mathcal{E}_{\mathrm{sep}}(\lambda)$ denote
the set of excited Young diagrams $D$ such that $wt_{\mathrm{sep}}(D)\neq 0.$
Then 
$$
E_\la(x|b)=\sum_{D\in\mathcal{E}_{\mathrm{sep}}(\lambda)}wt_{\mathrm{sep}}(D).
$$
We denote by $\mathcal{T}_0(\lambda)$ the set of ordinary (non-set-valued)
tableaux of shape $\lambda.$

\begin{lem}\label{lem:bijection} We have a bijection
$\mathcal{E}_{\mathrm{sep}}(\lambda)\rightarrow\mathcal{T}_0(\lambda).$
\end{lem}

{\it Proof.} We use a term 
{\it horizontal (resp. vertical) line} to call a row (resp. column) having 
nonzero weights in $wt_{\mathrm{sep}}$.
We put label $i$ to the horizontal line
having weights $x_i\oplus b_{j-i}$, and $j'$ to   
the vertical line having weight $x_j\ominus b_{j-i}.$
Let $D\in \mathcal{E}_{\mathrm{sep}}(\lambda).$
Since each box in $D$ is on some line, horizontal or vertical, 
we associate its label. 
Then we write the label in the original position of the box in $\D(\la)$. 
This gives a tableau of shape $\lambda.$
The fact that the resulting tableau is semistandard is 
readily seen from the configuration of the lines.
Conversely,
let $T\in \mathcal{T}_0(\lambda).$ 
Any box $\alpha$ in $\D(\lambda)$ 
has the right (excited) position, say $\alpha'$,  
on the line having label $T(\alpha)$.
Semistandardness ensures that
the set $\{\alpha'\}$ is  
 an excited Young diagram 
in $\mathcal{E}_{\mathrm{sep}}(\lambda).$
$\qed$

\begin{example}\label{exam:tab-EYD}
Tableau $T_0=\Tableau{1&2&3\\~&3'}$
corresponds to the following EYD:
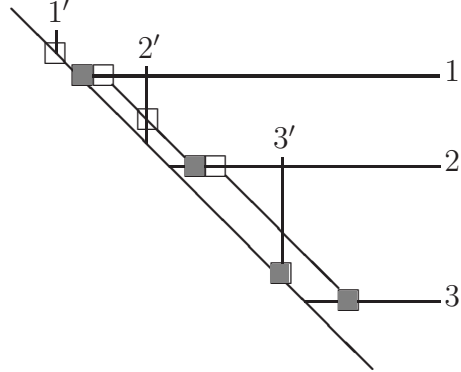
\begin{center}

\setlength{\unitlength}{0.6mm}
\begin{picture}(150,70)
\thicklines
\put(0,70){\line(1,-1){80}}
\put(10,60){\line(0,1){5}}
\put(15,55){\line(1,0){80}}
\put(8,67){$1'$}
\put(28,59){$2'$}
\put(58,39){$3'$}
\put(96,54){$1$}
\put(96,34){$2$}
\put(96,4){$3$}
%\put(120,2){$4$}
\put(7,57.7){$\square$}
%\put(6,87){$\times$}
\put(13,52.7){$\square$}
\put(17.7,52.7){$\square$}
\put(30,40){\line(0,1){17}}
\put(35,35){\line(1,0){60}}
\put(27.5,43){$\square$}
\put(38,32.6){$\square$}
\put(38.6,33){$\smallbox$}
\put(42.5,32.6){$\square$}
\put(57,8.8){$\square$}
%\put(62,33){$\square$}
%\put(67,33){$\square$}
\put(72,3){$\square$}
%\put(59,32){$\square$}
%\put(59,32){$\times$}
\put(60,10){\line(0,1){27}}
\put(65,5){\line(1,0){30}}
\put(47.5,33){\line(1,-1){25}}
\put(22.4,53){\line(1,-1){16}}
\put(13.7,53.2){$\smallbox$}
\put(72.5,3.3){$\smallbox$}
\put(57.4,9.3){$\smallbox$}
%\put(74.5,32){$\square$}
%\put(74.5,32){$\smallbox$}
%\put(90.5,1){$\square$}
%\put(90.5,1){$\smallbox$}
%\put(90,10){\line(0,1){100}}
%\put(94,5){\line(1,0){20}}
\end{picture}  %,\Tableau{2'&3'&3\\~&4},\Tableau{2&3'&3\\~&4'},\Tableau{2'&3'&3\\~&4'}$
\end{center}

\bigskip

Here $\smallbox$'s are the elements in $D$ whereas
 $\square$'s are the elements in $B(D).$
\end{example}

Let $T$ be in $\mathcal{T}(\lambda).$
The map $\D(\lambda)\rightarrow \A$ defined by 
 $(i,j)\mapsto \max\;T(i,j)$
gives an element in $\mathcal{T}_0(\lambda).$
Denote the resulting element by $\max(T).$
We denote $(x|b)^T_{I}$ in (\ref{eq:(x|b)^T}) simply by $(x|b)^T.$ 
\begin{lem}\label{lem:wsep=tab}
Let $T_0\in \mathcal{T}_0(\lambda)$ correspond to $D\in \mathcal{E}_{\rm sep}(\lambda)$  by the bijection of Lemma \ref{lem:bijection}. Then 
\begin{equation}
wt_{\mathrm{sep}}(D)=\sum_{T\in\mathcal{T}(\lambda), \;\max(T)=T_0}
\beta^{|T|-|\la|}(x|b)^T.\label{eq:wsep=tab}
\end{equation}
\end{lem}
If this 
lemma is true, the proof of Prop. \ref{prop:E=T} 
can be completed as follows:
\begin{eqnarray*}
E_\la(x|b)&=&
\sum_{D\in \mathcal{E}_{\mathrm{sep}}(\la)}
wt_{\mathrm{sep}}(D)\\
&=&
\sum_{T_0\in \mathcal{T}_0(\lambda)}\left(
\sum_{T\in\mathcal{T}(\lambda), \;\max(T)=T_0}
\beta^{|T|-|\la|}(x|b)^T\right)\\
&=&\sum_{T\in\mathcal{T}(\lambda)}\beta^{|T|-|\la|}(x|b)^T=T_\la(x|b),
\end{eqnarray*}
where
the second equality is
the consequence of bijection in Lemma \ref{lem:bijection}
and equation (\ref{eq:wsep=tab}).

\begin{example}\label{exam:E=T}
In order to see that equation (\ref{eq:wsep=tab}) holds, 
it is convenient to use weights 
$x_{j'}\ominus b_{j-i}$ on the vertical line
instead of $x_{j}\ominus b_{j-i}.$ 
This device makes clear that each term on the right-hand side 
of (\ref{eq:wsep=tab}) corresponds naturally to 
a term in suitably expanded form of $wt_{\rm sep}(D)$. 
If $D$ is the EYD in 
example \ref{exam:tab-EYD}, then $wt_{\rm sep}(D)$ is equal to 
the product of the weight in the following boxes: 

\setlength{\unitlength}{0.7mm}
\begin{picture}(180,25)
\put(0,20){\line(1,0){215}}
\put(0,10){\line(1,0){215}}
\put(34,0){\line(1,0){111}}
\put(0,10){\line(0,1){10}}
\put(34,0){\line(0,1){20}}
\put(145,0){\line(0,1){20}}
\put(215,10){\line(0,1){10}}
\put(2,13){$(1+\beta x_{1'})x_1$}
\put(37,13){$(1+\beta (x_{1}\oplus b_1))(1+\beta (x_{2'}\ominus b_1) )(x_2\oplus b_1)$}
\put(71,3){$x_{3'}$}
\put(150,13){$(1+\beta (x_2\oplus b_2))(x_{3}\oplus b_2)$}
\end{picture} 

\noindent First note that the factor $x_1(x_2\oplus b_1)(x_{3}\oplus b_2)x_{3'}$, when we ignore 
primes,  is
equal to $(x|b)^{T_0}.$ Expanding the remaining 
factors of the form $1+\beta wt_{\mathrm{sep}}(\alpha),$ 
which comes from $B(D)$, 
the weight $wt_{\mathrm{sep}}(D)$ can be expressed as 
a sum over tableaux
$$\Tableau{a1&b2&c3\\~&3'}$$
where $a,b,c$ are subsets of $\A$ such that  
$a\subset \{1'\},\;b\subset \{1,2'\},\;c\subset \{2\}.$
Note that 
these tableaux are exactly $T$'s in $\mathcal{T}(\lambda)$ such that 
$\max(T)=T_0.$
One sees that such $T$ gives a summand $\beta^{|T|-|\la|}(x|b)^T.$
For example, if $a=\{1'\},\,b=\{2'\},\,c=\emptyset$ then
for the corresponding $T$, we have the term 
$$
\beta^2 \cdot 
x_{1'}x_1
(x_{2'}\ominus b_1)(x_2\oplus b_1)
(x_{3}\oplus b_2)
x_{3'}.$$
We then ignore `primes' 
to have 
the right weight $\beta^2(x|b)^T.$
Thus we have (\ref{eq:wsep=tab}).
\end{example}

The argument in Example \ref{exam:E=T}
can be generalized to the following. 
 
 \bigskip
 
 {\it Proof of Lemma \ref{lem:wsep=tab}.}
For each $c=(i,j)\in \D(\lambda)$, let $R(c;T_0)$
denote the set of letters $a$ in $\mathcal{A}$ 
such that 
joining $a$ into $T_0$ at $c$ gives a semistandard 
set-valued 
tableau.
Now suppose, for all $c\in \D(\lambda),$
subsets $S_c$  in $R(c;T_0)$ are given. Then by joining $S_c$ into $T_0$ at $c$ for all $c\in \D(\lambda),$
we have a tableau $T\in \mathcal{T}(\la)$ such that $\max(T)=T_0.$
Conversely, all such tableaux are obviously given in this way.
This correspondence gives  
\begin{equation}
(x|b)^{T_0}
\sum_{(S_c)_{c\in \D(\lambda)} }
\prod_{c\in \D(\lambda)}
\beta^{|S_c|}\prod_{a\in S_c}w(c;a)=
\sum_{T\in \mathcal{T}(\la),\;\max(T)=T_0}\beta^{|T|-|T_0|}(x|b)^T,\label{eq:last}
\end{equation}
where $S_c$ runs for all the subsets in $R(c;T_0),$
and  
$w(c;a)$ is defined in \S \ref{ssec:SVT}.
Note that $(x|b)^{T_0}=\prod_{c\in D}wt_{\mathrm{sep}}(c).$
Finally as illustrated in Example \ref{exam:E=T}, we have
$$\sum_{(S_c)_{c\in \D(\lambda)}}
\prod_{c\in \D(\lambda)}
\beta^{|S_c|}\prod_{a\in S_c}w(c;a)
=\prod_{c\in B(D)}(1+\beta wt_{\mathrm{sep}}(c)).
$$
This together with (\ref{eq:last}) leads to (\ref{eq:wsep=tab}).
$\qed$

\begin{small} 
{\scshape Takeshi Ikeda, Department of Applied Mathematics,
Okayama University of Science,
Okayama 700-0005, JAPAN}
\end{small}
{\textit{email address}: \tt{ike@xmath.ous.ac.jp}}

\begin{small} 
{\scshape Hiroshi Naruse, Graduate School of Education,
Okayama University,
Okayama 700-8530, JAPAN}
\end{small}
{\textit{email address}: \tt{rdcv1654@okayama-u.ac.jp}}

\end{document}